\documentclass[journal,9pt]{IEEEtran}
\usepackage{graphicx}
\usepackage{hyperref}
\usepackage{mathrsfs}
\usepackage{amsmath}
\usepackage{amssymb}
\usepackage{url}
\usepackage{bm}
\usepackage{multirow}
\usepackage{diagbox}
\usepackage{multicol}
\usepackage{array}
\usepackage{booktabs}
\usepackage{amsthm}
\usepackage{amssymb}
\usepackage{color}
\usepackage[dvipsnames]{xcolor}
\usepackage[normal]{caption}
\usepackage{xcolor}

\graphicspath{{./},{./figures/}}

\newcommand{\black}{\color{black}}
\newcommand{\periwinkle}{\color{\periwinkle}}
\newcommand{\PID}{{\scriptscriptstyle {\it PID}}}
\newcommand{\PD}{{\scriptscriptstyle {\it PD}}}
\newcommand{\PI}{{\scriptscriptstyle {\it PI}}}
\newcommand{\p}{{\scriptscriptstyle {\it P}}}
\newcommand{\M}{{\scriptscriptstyle {\it M}}}

\begin{document}
\title{ Classical Stability Margins by PID Control
%Analytical Expressions of Maximal Gain and Phase Margins Attainable by PID Control
\thanks{This research was supported in part by the Natural Science Foundation of China under Grant 61876041, in part by Hong Kong RGC under the project CityU 11260016, in part by the City University of Hong Kong under Project 9380054, in part by NSF/USA (grants 1807664, 1839441), and in part by AFOSR/USA (grant FA9550-20-1-0029).}
}

\author{Qi~Mao,~\IEEEmembership
           \rm{Yong}~Xu,~\IEEEmembership
           \rm{Jianqi}~Chen,~\IEEEmembership
%            \rm{Islam}~Boussaada~\IEEEmembership
           \rm{Jie}~Chen,~\IEEEmembership{Fellow,~IEEE}
           Tryphon~Georgiou,~\IEEEmembership{Fellow,~IEEE}
\thanks{Y. Xu is with the School of Automation, Guangdong Provincial Key Laboratory of Intelligent Decision and Cooperative Control, Guangdong University of Technology, Guangzhou, China (e-mail: yxu@gdut.edu.cn.).}
\thanks{Q. Mao (e-mail: qimao3-c@my.cityu.edu.hk), Jianqi Chen (e-mail: jqchenlove@gmail.com), and Jie Chen (e-mail: jichen@cityu.edu.hk) are with the Department of Electrical Engineering, City University of Hong Kong, Hong Kong, China.}
%\thanks{I. Boussaada  is with the Department of Universit\'e Paris Saclay, CNRS-CentraleSupelec, Inria Saclay, Laboratoire %des Signaux et Syst\`emes (L2S), 3 rue Joliot-Curie 91192 Gif-sur-Yvette cedex (France)(e-mail: %Islam.Boussaada@l2s.centralesupelec.fr).}
% <-this % stops a space
%
%\thanks{Jianqi Chen is with the Department of Electronic Engineering, City University
%of Hong Kong, Hong Kong, China.}
%\thanks{Jie Chen is with the Department of Electrical Engineering, City University
%of Hong Kong, Kowloon, Hong Kong, China.}% <-this % stops a space
\thanks{T. Georgiou is with the Department of Mechanical $\&$ Aerospace Engineering, University of Calfornia, Irvine, CA 92697 USA (email: tryphon@uci.edu).}
% \href{mailto: jichen@cityu.edu.hk}{jichen@cityu.edu.hk}
}

\markboth{IEEE Transactions on Automatic Control}%
{Mao \MakeLowercase{\textit{et al.}}:
Classical Stability Margins by PID Control
%Analytical Expressions of Maximal Gain and Phase Margins Attainable by PID Control
}
%Maximal Gain and Phase Margin of Low-Order Systems Achievable by PID Controllers




\maketitle

\begin{abstract}
Proportional-Integral-Derivative (PID) control has been the workhorse of control technology for about a century. Yet to this day, designing and tuning PID controllers relies mostly on either tabulated rules (Ziegler-Nichols) or on classical graphical techniques (Bode). Our goal in this paper is to take a fresh look on PID control in the context of optimizing stability margins for low-order (first- and second-order) linear time-invariant systems. Specifically, we seek to derive explicit expressions for gain and phase margins that are achievable using PID control, and thereby gain insights on the role of unstable poles and nonminimum-phase zeros in attaining robust stability. In particular, stability margins attained by PID control for minimum-phase systems match those obtained by more general control, while for nonminimum-phase systems, PID control achieves margins that are no worse than those of general control modulo a predetermined factor. Furthermore, integral action does not contribute to robust stabilization beyond what can be achieved by PD control alone.
\black
\end{abstract}


\IEEEpeerreviewmaketitle

\section{Introduction}
\IEEEPARstart{T}{he} primary goal of feedback regulation is to maintain stability and
performance in the presence of modeling uncertainty and external disturbances.
Amongst the metrics that are commonly used to quantify robustness against such factors, traditionally, the most important
have been gain and phase margins, various types of induced norms ($H_\infty$, $L_\infty$), the gap metric, the structured singular value, and so forth. Each of these, and a few others, have been the subject of respective chapters in the modern robust control literature.
Herein, we focus on gain and phase margins that historically have been the first to be considered. Interestingly but perhaps unsurprisingly, these same metrics have also been the first to be tackled in the waning years of the 1970's with the modern tools of analytic function theory that gave rise to optimal designs \cite{tannenbaum1980feedback,tannenbaum1981invariance,tannenbaum1982modified}; see also \cite{khargonekar1985non,tannenbaum1986multivariable,feintuch1986gain} as well as \cite[Chapter 11]{doyle2013feedback}.

Gain and phase margins quantify tolerance of stability of a feedback system to perturbation in the respective features of its transfer function, i.e., the gain and phase. Historically these metrics proved central in feedback regulation of electronic amplifiers \cite{Bode45,Horowitz63} and are nowadays taught in every introductory course in control.
A vast volume of work is in existence, starting perhaps with Ziegler-Nichols \cite{ziegler1942optimum} and continuing to this day, aimed at analysis and synthesis techniques to ensure quantitative bounds on gain and phase for maintaining stability of linear time-invariant (LTI) systems (see \cite{maeda1986infinite,yan1989simultaneous,aastrom1995pid,ho1996performance,tantaris2003gain,aastrom2004revisiting,paraskevopoulos2006pid,wang2009synthesis} and the references therein).

Retaining stability against parametric and non-parametric uncertainty has been the subject of much debate and confusion in the years prior to the development of modern robust control theory in the early 1980's. Pivotal and influential in this journey have been works by Kalman \cite{kalman1964linear} who sought to quantify properties of optimally designed controllers,
exploration of the same in multivariable designs by Safonov and Athans \cite{safonov1977gain},
the influential attestation by Doyle on the absence of guaranteed margins in LQG regulators \cite{doyle1978guaranteed},
and several others.
%At the same time, curious differences between continuous and discrete-time with regard to robust stability were noted \cite{willems1978return}.
These and other works, to a large degree, set the stage for the urgency in the subsequent development of modern robust/$\mathcal{H}_\infty$ control methods \cite{zames1981feedback}.
It should however be noted that the earlier development by Tannenbaum \cite{tannenbaum1980feedback,tannenbaum1981invariance,tannenbaum1982modified} of optimal gain-margin designs for scalar LTI systems, introducing Nevanlinna-Pick interpolation in the context of robust control, was independent of the mainstream control literature. Tannenbaum referred to gain margin optimization as the ``blending problem'' -- an unconventional terminology which may have contributed to a delay in recognizing its connection to the fast developing $\mathcal{H}_\infty$-control literature at the time.

Interest in a deeper understanding of obstructions in achieving improved stability margins has persisted over the years, due to the practical significance of such metrics. Indicatively, we cite
\cite{yan1989simultaneous,%cockburn1992constructive,
cockburn1995stability} on multi-objective margins
and \cite{lee1987vibrational,francis1988stability,cusumano1988nonlinear} that highlight non-LTI controllers in allowing improved robustness margins while sacrificing other performance objectives of a feedback system.
Returning to the theme of our current paper,
undoubtedly, PID control has long been established for its simplicity, ease of implementation, and its cost-effectiveness \cite{aastrom1995pid}. Even in the post-modern and information-centric era of present time, it continues to dominate industrial control systems design as the most favored control technology. Indeed, the recent survey
\cite{samad2017survey} suggests that PID controllers continue to be most widespread in spite of extraordinary advances in control theory and design techniques.
Thus, while PID control has for a long time been by and large an {\em ad hoc} method, whose design and tuning drew heavily upon empirical results and trial-and-error studies, analytical studies of PID control began to emerge in the recent years, on a variety of problems ranging from performance analysis
\cite{aastrom1995pid,gundecs2020controller} to delay robustness \cite{Silva2002New,ma2018delay,chen2021delay},
pole placement \cite{michiels2002continuous,zitek2013dimensional,wang2017new},  and nonlinear regulation \cite{zhang2019theory, zhao2020control}. Notwithstanding these advances, the robustness qualities of PID control remain largely an open subject. Earlier results on the gain and phase margin enhancement by PID controllers include, e.g.,
\cite{ho1995tuning}, %\cite{ho1997self}, %\cite{ho1998optimal}, \cite{ho1999getting},
\cite{yaniv2004design}, where simple formulas and schemes were devised to tune PID controller parameters so that the required gain and phase margins can be met.
Hence, we herein explore more systematically, and compare PID control design, against the optimal achievable gain and phase margins by arbitrary LTI controllers.

\black
Specifically, in the present paper, we study the gain and phase margins attainable by PID control.
We consider first- and second-order unstable plants, with one or two unstable poles and possibly a nonminimum-phase zero; the gain and phase margin maximization problems are meaningful only for unstable plants. We note that, on one hand,
PID control is essentially pertinent for first- and second-order systemms \cite{aastrom1995pid,krstic2017applicability}, and on the other, industrial  control systems where PID control is prevalent rely heavily on such low-order  models. The simplicity in plant dynamics, albeit restrictive, renders our problems analytically tractable.
Our main contribution is in the derivation of explicit expressions of the maximal gain and phase
margins, tabulated in Table I. These results reveal how the
unstable poles and nonminimum-phase zeros limit achievable margins. In particular, our findings lead to the following conclusions:\\[0.03in]
i) PID and PD controllers achieve the same maximal stability margins;
(non-vanishing)
integral action can only reduce the margins.\\[0.03in]
ii) PID control achieves the same maximal gain and phase margins for first- and second-order  ``minimum-phase'' systems as general LTI controllers.\\[0.03in]
iii)
For first- and second-order nonminimum-phase plants, PID control attains {\em at least half of the maximal gain margin} achievable by general LTI controllers. While an analytical proof is currently unavailable in general, numerical results suggest that the same conclusion holds for the maximal phase margin as well. In fact, it can be shown that under suitable conditions,
the maximal phase margin achievable by an LTI controller is indeed {\em at most twice} that by PID controllers.\\[0.03in]
iv) For unstable and non-minimum phase systems, a fundamental limit on the maximal PID phase margin is $90^0$, while the maximal gain margin achievable can still be large, tending to infinity as the nonminimum-phase zero approaches the imaginary axis. Thus, while the statements i)-iii) serve as a testimony to the potency of PID control,
the statement iv) demonstrates a fundamental limitation of PID control, in spite of its other attributes.

\noindent Conceptually appealing, the results alluded to above establish the strong
robustness properties of PID controllers when appropriately
designed, and from the perspective of gain and phase margins,
lend a theoretical justification for the successes of PID control,
leading to a theory of ``explainable'' analytical design and
tuning of PID controllers. From a practical standpoint, the results provide insights and tuning rules for unstable systems,
to which the Ziegler-Nichols rules and their variants fail to
be applicable.  It is worth noting that in developing these results,
even for second-order systems, the problem is surprisingly challenging,
and arguably more so than its non-parametric $\mathcal{H}_\infty$ counterpart.
Thus, while our intention has been to be concise, the development is by necessity
intertwined and requires a delicate combination of tools drawing upon nonlinear programming and algebraic geometry.

The structure of the paper is as follows.
In Section II, we introduce the basic problem, which amounts to maximizing gain and phase margins using PID/PI/PD control.
Parametric nonlinear programming is brought in to tackle this problem in Section III for first-order systems, and in Section IV for second-order systems. More specifically, in Section III we employ the {\em Bilherz criterion} to solve the maximal
phase margin problem for first-order systems.
For second-order nonminimum-phase plants, the problem becomes substantially more involved. To this end we explore a nonlinear programming approach. This is undertaken in Section IV, where we exploit the {\em Karush-Kuhn-Tucker (KKT)} condition. The maximal phase margin is analytically determined by solving
the positive real roots of two third-order polynomials; this in turn can be obtained in closed-form via the {\em Cardano formula}. PI control is addressed
in Section V, where the maximal phase margin can also be obtained analytically. To streamline the presentation, we relegate the proofs of
the main results to appendices provided at the end of the paper, along with the necessary mathematical tools.
%, including the Bilherz criterion for stability of complex polynomials, the KKT condition for nonlinear programming problems, and the {\em Descartes Rule of Signs} %concerning the number of positive real roots of a polynomial.

Partial results of this paper were previously presented
in an abridged conference version \cite{mao2021maximal}, where due to space constraint all proofs were omitted.

\black

\section{Problem Formulation}

%%%%%%%%%%%%%%%%%%Notation
%\subsection{Notation}
%The notation used in this paper is fairly standard.

%%%%%%%%%%%%%%%%%%Notation
%\subsection{Gain and Phase Margin Maximization Problem}

%\subsection{Problem Formulation}
We consider the feedback system depicted in Fig.1, in which $P(s)$ represents the plant model, and $K(s)$ a finite-dimensional LTI controller.
Suppose that $P(s)$ is stabilized by $K(s)$.
\begin{figure}[!htb]
\centering
\includegraphics[width=5.50 cm]{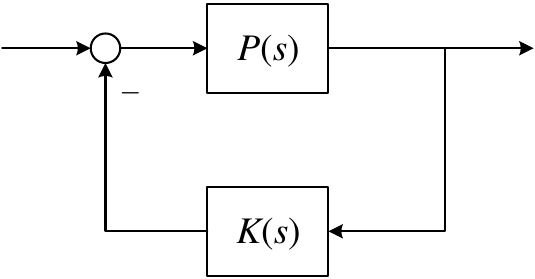}
\caption{Standard feedback control structure}
\label{Fig1}
\end{figure}
Then if $P(s)$ is perturbed by uncertainties in the form of gain or phase within a sufficiently small range, by continuity, the controller $K(s)$ still stabilizes the
perturbed plant. But how large a perturbation in gain or phase can be tolerated before the feedback system becomes unstable?

{Gain/phase margin maximization} provides a quantitative answer to the above question.
Following \cite[Chapter 11]{doyle2013feedback}, consider the family of plants
\begin{flalign} \label{2-1}
{\mathscr{P}_\mu } = \left\{ {\alpha P\left( s \right)\begin{array}{*{20}{c}}
%{\mathcal{P}_\alpha } = \left\{ {\alpha P\left( S \right)\begin{array}{*{20}{c}}
:&{1 \le \alpha  < \mu }
\end{array}} \right\}.
\end{flalign}
This represents uncertainty in the gain about the nominal plant model $P(s)$, and thereby, the {\em maximal gain margin} for $P(s)$,
\begin{equation}\nonumber
\begin{array}{l}
{k_\M} = \sup \{ \mu : \exists K\left( s \right) \mbox{ stabilizing } \alpha P\left( s \right), \; {\rm{ }} \; \forall \alpha  \in [1,\mu )\},
\end{array}
\end{equation}
specifies the largest family $\mathscr{P}_{k_\M}$ that can be stabilized by the same LTI controller.
Likewise, consider the family
\begin{equation} \label{2-2}
{\mathcal{Q}_\nu } = \left\{ {{e^{ - j\theta }}P\left( s \right)\begin{array}{*{20}{c}}
:&{\theta  \in ( - \nu ,\nu )}
\end{array}} \right\},
\end{equation}
and define the {\em maximal phase margin}
\begin{equation}\nonumber
\begin{array}{l}
{\theta _\M} = \sup \{ \nu : \exists K\left( s \right)\mbox{ stabilizing } {e^{ - j\theta }}P\left( s \right), {\rm{ }} \; \forall \theta  \in ( - \nu ,\nu )\},
\end{array}
\end{equation}
which specifies the largest family $\mathcal{Q}_{\theta_\M}$ stabilizable by the same suitably chosen LTI controller.

%Stated in words, the maximal gain margin ${k_\M} $ and the maximal phase margin ${\theta _\M}$ determine, respectively, the maximal ranges of gain and phase variations within which the families of ${\mathscr{P}_\alpha } $ and ${\mathcal{Q}_\theta }$ can be robustly stabilized by a LTI controller.

\textit{Remark 2.1}: The definition of gain margin above and in \cite[Chapter 11]{doyle2013feedback} serves some notational convenience, whereas a more practical definition would place the nominal plant in the (logarithmic) center of the parametric uncertainty, precisely as in \cite{tannenbaum1980feedback,tannenbaum1981invariance,khargonekar1985non}, requiring that a controller stabilizes a maximal family $\{\alpha P(s) :  \sqrt{\mu^{-1}}<\alpha <\sqrt{\mu}\}$. Indeed, optimizing for a one-sided perturbation of the gain necessitates that an $\epsilon$-suboptimal controller stabilizes gain perturbation in an interval $1-\epsilon\leq \alpha\leq k_\M-\epsilon$, placing the nominal plant in a precarious spot rendering the feedback system fragile to perturbations in one of the two direction. This very point became a source of controversy \cite{makila1998comments} in that earlier authors raised a debateable issue of fragility as endemic to gain margin maximization.
%Initially posed as a so-called ``blending'' problem, the maximization problem of the gain margin was solved in
%\cite{tannenbaum1980feedback,tannenbaum1981invariance} based on the Nevanlinna-Pick interpolation theory \cite{Ball90},
%which can be tackled more generally as an $\mathcal{H}_\infty$ optimal
%control problem \cite{khargonekar1985non,doyle2013feedback}. As in these formulations, we may consider $\alpha \in
%({\underline \mu},~ \bar \mu) $, where $0<{\underline \mu}<1<\bar \mu$. Suppose that $P\left( s \right)$
%can be stabilized by some $K(s)$. Define analogously
%\begin{equation}\nonumber
%\begin{array}{l}
%{\bar k_\M} = \sup \{ \bar \mu :{\rm{There \; exists \; some \; }}K\left( s \right){\rm{ \;  stabilizing}}\\
%{\kern 145pt} \alpha P\left( s \right), \; {\rm{ }} \; \forall \alpha  \in [1,\bar  \mu )\} ,
%\end{array}
%\end{equation}
%and
%\begin{equation}\nonumber
%\begin{array}{l}
%{\underline k_\M} = \inf \{ {\underline \mu} :{\rm{There \; exists \; some \; }}K\left( s \right){\rm{ \;  stabilizing}}\\
%{\kern 145pt} \alpha P\left( s \right), \; {\rm{ }} \; \forall \alpha  \in ({\underline \mu},1]\}.
%\end{array}
%\end{equation}
%Then a more general notion of the maximal gain margin can be defined as ${\bar k_\M}/{\underline k_\M}$. Similarly, the definition of
%$\theta_\M$ can be extended to an asymmetric interval.
$\hfill\blacksquare$

\textit{Remark 2.2}: At times, when no confusion
arises, we may measure the gain margin in decibels (dB), i.e., referring to $20\log_{10} k_\M$ as the maximal gain margin. $\hfill\blacksquare$

\black

%\subsection{Stability Margins}
Define the system's complementary sensitivity function by
$$
T\left( s \right) = P\left( s \right)K\left( s \right){\left[ {1 + P\left( s \right)K\left( s \right)} \right]^{ - 1}},
$$
and let
\begin{equation}\label{eqref:complementary}
{\gamma _{{\rm{opt}}}} = \inf \left\{ {{{\left\| T(s) \right\|}_\infty }: \; K\left( s \right) \; {\rm{ stabilizes }} \; P\left( s \right)} \right\},
\end{equation}
with $\|T(s) \|_\infty$ being the $\mathcal{H}_\infty$-norm of $T(s)$.
The following result due to Tannenbaum \cite{tannenbaum1980feedback,tannenbaum1982modified} (see also \cite{khargonekar1985non,doyle2013feedback} for a lucid exposition) provides analytical expressions of
$k_\M$ and $\theta_\M$ that can be determined by solving the standard $\mathcal{H}_\infty$ optimal control
problem in \eqref{eqref:complementary}.

\textit{Proposition 2.1}: (i). If $P$ is stable or minimum-phase, then the maximal gain margin ${k_\M} = \infty $. Otherwise,
\begin{equation} \label{2-6}
{k_\M} = {\left( {\frac{{{\gamma _{{\rm{opt}}}} + 1}}{{{\gamma _{{\rm{opt}}}} - 1}}} \right)^2}.
\end{equation}
 (ii). If $P$ is stable or minimum-phase, then the maximal phase margin is ${\theta_\M} = \pi$ (in radians). Otherwise,
\begin{equation} \label{2-7}
{\theta_\M} = 2{\sin ^{ - 1}}\frac{1}{{{\gamma _{{\rm{opt}}}}}}.
\end{equation}

The goal of the present paper is determine the maximal robustness margins that are achievable by PID controllers; that is, by controllers with transfer function of the form %$K\left( s \right) = {K_{\PID}}\left( s \right)$,
%where
\begin{equation} \label{2-3}
{K_{\PID}}\left( s \right) = {k_p} + \frac{{{k_i}}}{s} + {k_d}s.
\end{equation}
Thus, accordingly, we define
\begin{equation}\nonumber
\begin{array}{l}
k_\M^{\PID} = \sup \{ \mu :~\exists  {K_{\PID}}\left( s \right){\rm{\; stabilizing~}} \alpha P\left( s \right) \; \forall \alpha  \in [1,\mu )\} ,
\end{array}
\end{equation}
\begin{equation}\nonumber
\begin{array}{l}
\theta _\M^{\PID} = \sup \{ \nu :~\exists {K_{\PID}}\left( s \right){\rm{ \; stabilizing}} {e^{ - j\theta }}P\left( s \right)\; \forall \theta  \in ( - \nu ,\nu )\} .
\end{array}
\end{equation}
Likewise, it is of interest to consider subclasses of PID controllers, such as proportional, PI, and PD
controllers, given by
\begin{equation} \label{2-4}
\begin{array}{*{20}{cc}}
\displaystyle {{K_\p}\left( s \right) = {k_p},}& \displaystyle K_{\PI}(s)=k_p+\frac{k_i}{s}, & \displaystyle {{K_{\PD}}\left( s \right) = {k_p} + {k_d}s},
\end{array}
\end{equation}
respectively. We shall denote the corresponding maximal gain and phase margins with the corresponding superscripts P, PI, and PD.
%
%{\color{red} Not needed here?
%Consider the open-loop transfer function
%\begin{equation} \label{2-4+1}
%L\left( s \right) = P\left( s \right){K_{\PID}}\left( s \right).
%\end{equation}}
%

In our analysis, we only consider first- and second-order plants. This is partly due to the  fact that industrial processes are often modeled by low-order, and in fact, mostly first-order and second-order systems, and partly due to the limitation of PID controllers in controlling high-order dynamics as, in general, PID controllers may not be able to stabilize certain third- and higher-order unstable systems \cite{krstic2017applicability,zhao2020control}. From a technical standpoint, unlike $k_\M$ and $\theta_\M$, computing maximal PID gain and phase margins
requires solving a parametric optimization problem.

Our main results are summarized in Table 1, which provides analytical expressions of maximal gain and phase margins for
first- and second-order plants. The rest of the paper is devoted to establishing these results and to the discussion of our findings.
\black

\begin{table*}[!htbp] \label{Stabiliby Margin Table1}
\centering
%\begin{tabular}{c|c|c|c}
\begin{tabular}{ccccc}
\specialrule{0em}{5pt}{5pt}
%\toprule
\hline
\specialrule{0em}{1.5pt}{1.5pt}
%\multicolumn{2}{c|}{Plant $P\left( s \right)$}   &Maximal Gain Margin $k_\M^{\PID} $   &Maximal Phase Margin $\theta_\M^{\PID} $ \\
\multicolumn{2}{c}{Plants $P\left( s \right)$}   &Gain Margin $k_\M^{\PID} $   &Phase Margin $\theta_\M^{\PID}$ &Relation with ${k_\M},{\theta _\M}$ \\
\hline
\specialrule{0em}{1.5pt}{1.5pt}
\multirow{2}*{First-Order Systems}                    &minimum-phase Plants &$\infty$         &$\pi$         &$\begin{array}{l}
{k_\M} = k_\M^{\PID}\\
{\theta _\M} = \theta _\M^{\PID}
\end{array}$\\
\cline{2-5}
\specialrule{0em}{1.5pt}{1.5pt}
&$\frac{{s - z}}{{s - p}},~p,z > 0$    &${\max \left\{ \frac{z}{p},~\frac{p}{z} \right\}}$        &${\cos ^{ - 1}}\left( {\frac{{2\sqrt {z/p} }}{{1 + \left( {z/p} \right)}}} \right)$        &$\begin{array}{l}
{k_\M} = {(k_\M^{\PI})^2}\\
{\theta _\M} = 2\theta _\M^{\PI}
\end{array}$\\
\cline{1-5}
\specialrule{0em}{1.5pt}{1.5pt}
\multirow{2}*{Second-Order Systems}               &minimum-phase Plants                     &$\infty$          &$\pi$                &$\begin{array}{l}
{k_\M} = k_\M^{\PID}\\
{\theta _\M} = \theta _\M^{\PID}
\end{array}$\\
\cline{2-5}
\specialrule{0em}{1.5pt}{1.5pt}
&$\frac{{s - z}}{{\left( {s - {p_1}} \right)\left( {s - {p_2}} \right)}},~{p_1},{p_2},z > 0$                     &explicit expression  &analytical solution               &$k_\M^{\PID} \le {k_\M} \le {(k_\M^{\PID})^2}$\\
\hline
%\bottomrule
\end{tabular}
\setlength{\abovecaptionskip}{2pt}
\caption{Stability Margins: Gain and Phase.}
\end{table*}

\section{First-Order Unstable Systems}

In this section, we provide  explicit expressions of maximal gain and phase margins of first-order systems achievable by PID controllers.
The results  show the dependence of these measures on the system's unstable pole and nonminimum-phase zero.

%main results for first-order systems, which consist of two cases, minimum-phase system and nonminimum-phase system. In particular, the explicit expressions of maximal gain margin and phase margin are presented for systems with nonminimum-phase zero, which shows the dependence of gain margin and phase margin on the system's unstable poles and nonminimum-phase zeros.

%%%%%%%%%%%%%%%%%%
%\subsection{PI Control}
We first consider first-order unstable plants described by
\begin{equation} \label{1-8}
P\left( s \right) = \frac{{{\beta_0}s + {\beta_1}}}{{s - p}}, \;\; p > 0.
\end{equation}
Without loss of generality, it is assumed that ${\beta_0} \ge 0$ and ${\beta_1} \ne 0$. Note that in the case ${\beta_0} \ne 0$, derivative control will result in an improper system. For this reason, in this section we shall focus on PI controllers.

%In the case ${b_0} \ne 0$, however, we cannot employ derivative control to stabilize $P\left( s \right)$ due to the introduction of an improper system. Under this circumstance, we implement the PI controllers. On the other hand, it is typical to use full PID controllers in the case ${b_0} > 0$. The maximal gain and phase margins are provided in the following theorem.

\textit{Theorem 3.1}: Let $P(s)$ be given by (\ref{1-8}). Then %under Assumption 2.1,
 the following statements hold:

(i). For ${\beta_0} > 0$,
\begin{equation} \label{2-9}
\begin{array}{*{20}{l}}
{\kern -8pt}{k_\M^{\PI} = k_\M^\p = \left\{ {\begin{array}{*{20}{l}}
\displaystyle {{\max \left\{ {\frac{{\left| {{\beta_1}} \right|/{\beta_0}}}{p},\frac{p}{{\left| {{\beta_1}} \right|/{\beta_0}}}} \right\},\;\;\;{\beta_1} < 0,}}\\
\displaystyle {\infty,\;\;\;\;\;\;\;\;\;\;\;\;\;\;\;\;\;\;\;\;\;\;\;\;\;\;\;\;\;\;\;\;\;\;\;\;\;\;\;\,{\beta_1} > 0.}
\end{array}} \right.}\\
{\kern -8pt}{\theta _\M^{\PI} = \theta _\M^\p = \left\{ {\begin{array}{*{20}{l}}
\displaystyle{{{\cos }^{{\rm{ - }}1}}\frac{{2\sqrt {{\beta_0}\left| {{\beta_1}} \right|/p} }}{{{\beta_0} + (\left| {{\beta_1}} \right|/p)}},\;\;\;\;\;\;\;\;\;\;\;\;\;\,{\beta_1} < 0,}\\
\displaystyle {\pi, \;\;\;\;\;\;\;\;\;\;\;\;\;\;\;\;\;\;\;\;\;\;\;\;\;\;\;\;\;\;\;\;\;\;\;\;\;\;\;\;\;\,{\beta_1} > 0.}
\end{array}} \right.}
\end{array}
\end{equation}
%where $\left| {{b_1}} \right|/{b_0} > p$ for ${b_1} < 0$.

(ii). For ${\beta_0} = 0$,

\begin{equation} \label{2-10}
\begin{array}{l}
{\kern -126pt} k_\M^{\PID} = k_\M^{\PD} = k_\M^\p = \infty ,\\
{\kern -126pt} \theta _\M^{\PID} = \theta _\M^{\PD} = 2\theta _\M^\p = \pi .
\end{array}
\end{equation}

\textit{Proof}: See Appendix \ref{Theorem 3.1-Proof}. $\hfill\blacksquare$

The expressions given in Theorem 3.1 lead us to a number of pertinent observations. First, it is clear that integral control has no effect on maximizing either the gain or the phase margin. This is consistent with one's intuition; integral control has its essential utility in tracking reference signals.
%, which is seen as a conflict with a system's stability robustness and henceforth with increasing the system's gain and phase margins.
Secondly, we note that for minimum-phase plants (${\beta_0} > 0$ and ${\beta_1}  >  0$) of relative degree zero, proportional control suffices to achieve the maximum possible infinite  gain margin and a phase margin of $\pm {180^ \circ }$. For nonminimum-phase plants (${\beta_0} > 0$ and ${\beta_1}  <  0$), it is instructive to consider
\begin{equation} \label{new3-1}
P\left( s \right) = \frac{{s - z}}{{s - p}},\;\;\;\; p > 0, \;\;z > 0.
\end{equation}
In this case, it follows from Theorem 3.1 that
\begin{equation} \label{new3-2}
\begin{array}{l}
\displaystyle  {k_\M^{\PI} = k_\M^\p = \max \left\{ {\frac{z}{p},~\frac{p}{z}} \right\}}.\\
\displaystyle  {\theta _\M^{\PI} = \theta _\M^\p = {\cos ^{{\rm{ - }}1}}\frac{{2\sqrt {z/p} }}{{1 + \left( {z/p} \right)}}}.
\end{array}
\end{equation}
It is interesting to see that
\begin{equation} \label{PI3-1}
\theta _\M^{\PI} = \theta _\M^\p = {\cos ^{ - 1}}\frac{{2\sqrt {k_\M^\p} }}{{1 + k_\M^\p}} = {\cos ^{ - 1}}\frac{{2\sqrt {k_\M^{\PI}} }}{{1 + k_\M^{\PI}}}.
\end{equation}
%A subsequent comparison of (\ref{new3-2}) with Proposition 2.1 (cf. Corollary 3.1) shows that for a first-order nonminimum-phase plant, the maximal phase margin achievable by proportional control is half that achievable by general LTI controllers, and when measured on the logarithmic scale, the maximal gain margin is also half that by  LTI controllers.

%%%%%%%%%%%%%%%%%%

%%
We conclude this section with a comparison of the maximal gain and phase margins herein with those achievable by the general LTI controllers. The following corollary states %the comparison for both the PI controller ${K_{PI}}\left( s \right)$. %and the filtered PID controller $K_{\PID}^{T_f}\left( s \right)$. The results show
that for a first-order nonminimum-phase plant, the maximal gain and phase margins achievable by a proportional controller are half as good as those by general LTI controllers.
%With derivative control added, an appropriately designed filtered PD controller can restore the maximal margins to those by LTI controllers.

%%

\textit{Corollary 3.1}: Let $P(s)$ be given by (\ref{new3-1}). %Denoted by $k_\M$, $\theta _M$ the maximal gain and phase margins by the general LTI controllers $K\left( s \right)$.
Then%for $K\left( s \right)$to stabilizes $P(s)$, the following statements holds,
\begin{equation} \label{Corollary 3.1-1}
\begin{array}{l}
%{k_\M} = {\left( {k_\M^{\PI}} \right)^2} = {\left( {k_\M^\p} \right)^2}.
{\log _{10}}{k_\M}~~ = ~~2{\log _{10}}k_\M^{\PI}~~ =~~ 2{\log _{10}}k_\M^\p.\\
{\theta _\M} ~~~~~\;\,~~= ~~2\theta _\M^{\PI} ~~~~~~\,~~= ~~2\theta _\M^\p.
\end{array}
\end{equation}

%(ii)
%\begin{equation} \label{Corollary 3.1-2}
%\theta _M=  2 \theta _\M^{\PI} = 2 \theta _\M^\p.
%\end{equation}

%(ii). For ${T_f} = {1}/{p}$,
%\begin{equation} \label{Corollary 3.1-3}
%k_\M = k_{M,{T_f}}^{\PID} = k_{M,{T_f}}^{PD},
%\end{equation}
%where ${T_f} = \frac{1}{p}$ when $z > p$ and ${T_f} = \frac{p}{{{z^2}}}$ when $p > z$.

\textit{Proof}: With $P(s)$ given by (\ref{new3-1}), it follows from \cite{chen2000logarithmic} (see also \cite{khargonekar1985non}) that
$${\gamma _{{\rm{opt}}}} = \left| {\frac{{z + p}}{{z - p}}} \right|.$$
Invoking Proposition 2.1 yields
\begin{equation}\nonumber
{k_\M} = \max \left\{\left( {\frac{z}{p}} \right)^2,~\left( {\frac{p}{z}} \right)^2 \right\},
\end{equation}
%and
\begin{equation}\nonumber
{\theta _\M} = 2{\cos ^{{\rm{ - }}1}}\frac{{2\sqrt {z/p} }}{{1 + \left( {z/p} \right)}}.
\end{equation}
In view of (\ref{new3-2}), the corollary follows.
%The statement (ii) follows from  (\ref{1-18}) and  (\ref{1-19}) by evaluating $k_{M,{T_f}}^{\PID}, k_{M,{T_f}}^{PD}$ at ${T_f} = {1}/{p}$.%it is not difficult to establish (\ref{Corollary 3.1-1}). With ${T_f} = \frac{1}{p}$, we have $k_{M,{T_f}}^{\PID} = k_{M,{T_f}}^{PD} = k_\M $.
$\hfill\blacksquare$

\section{Second-Order Unstable Systems}
%We now study the second-order unstable systems.
This section presents results for second-order unstable plants, which include minimum-phase as well as nonminimum-phase systems. In general, the computation of the gain and phase margins for second-order plants poses a more difficult problem. Throughout this section, our development seeks to recast the gain and phase maximization problems as one of nonlinear programming.
%Since the derivations are lengthy, we omit the proofs of the subsequent results.
For their essential flavor, we employ the KKT condition to obtain value of the integral gain as the necessary solution to the nonlinear programming problems under consideration,
which show that unequivocally in all cases, the optimal integral gain is $k_i=0$. % and the optimal derivative gain $k_d$ (or proportional gain $k_p$) is a certain boundary value.
The maximal gain and phase margins are then obtained by solving a univariate optimization problem, defined in terms of a function of the proportional gain $k_p$ or the derivative gain $k_d$ alone.
%For their essential flavor, we employ the KKT condition to obtain values of the integral and derivative gains as the necessary solution to the nonlinear programming problems under consideration, which show that unequivocally in all cases, the optimal integral gain is $k_i=0$, and the optimal derivative gain $k_d$ is a certain boundary value. The maximal gain and phase margins are then obtained by solving a univariate optimization problem, defined in terms of a function of the proportional gain $k_p$ alone.

%%%%%%%%%%%%%%%%%%
\subsection{Minimum-phase Systems}
We begin with minimum-phase plants that contain a pair of unstable poles ${p_1}$, ${p_2}$, described by
\begin{equation} \label{3-1}
P\left( s \right) = \frac{{{\beta _0}s + {\beta _1}}}{{\left( {s - {p_1}} \right)\left( {s - {p_2}} \right)}},\;\;{\mathop{\rm Re}\nolimits} \left( {{p_1}} \right) > 0,\;\;{\mathop{\rm Re}\nolimits} \left( {{p_2}} \right) > 0,
\end{equation}
where ${\beta _0} \ge 0,{\beta _1} > 0$. To ensure that $P\left( s \right)$ is a real rational plant, we assume that ${p_1}$ and ${p_2}$ are either real poles or form a complex conjugate pair.

\textit{Theorem 4.1}: Let $P(s)$ be given by (\ref{3-1}). Then %under Assumption 2.1,
the following statements hold:

(i). \begin{equation} \label{2-20}
\begin{array}{l}
k_\M^\p = \left\{ {\begin{array}{*{20}{l}}
{{\rm{none}},\;\;\;{\rm{when}}\;{\beta_0} = 0},\\
{\infty, \;\;\;\;\;\;\;{\rm{when}}\;{\beta_0} \ne 0.}
\end{array}} \right.\\
k_\M^{\PID} = k_\M^{\PD} = \infty.
\end{array}
\end{equation}

(ii). \begin{equation} \label{2-21}
\begin{array}{l}
\theta _\M^\p = \left\{ {\begin{array}{*{20}{l}}
{{\rm{none}},\;\;\;{\rm{when}}\;{\beta_0} = 0,}\\
{\pi/2,\;\;\;\;\;\;\;\;{\rm{when}}\;{\beta_0} \ne 0.}
\end{array}} \right.\\
\theta _\M^{\PID} = \theta _\M^{\PD} = \left\{ {\begin{array}{*{20}{l}}
{\pi/2,\;\;\;{\rm{when}}\;{\beta_0} = 0,}\\
{\pi ,\;\;\;\;{\rm{when}}\;{\beta_0} \ne 0.}
\end{array}} \right.
\end{array}
\end{equation}

%Specifically, Define the set
%\begin{array}{l}
%\Omega  = \left\{ {\left( {{k_p},{k_i},{k_d}} \right):} \right.{k_p} >  - \frac{{{b_1}}}{{{b_0}}}{k_d} + \frac{{{p_1} + {p_2}}}{{{b_0}}},\\
%{\kern 75pt}\left. {{k_p} >  - \frac{{{b_0}}}{{{b_1}}}{k_i} - \frac{{{p_1}{p_2}}}{{{b_1}}},{k_i} > 0,{k_d} >  - \frac{1}{{{b_0}}}} \right\}.
%\end{array}
%\end{equation}
%It is necessary and sufficient that ${k_p} >  - \frac{{{p_1}{p_2}}}{{{b_1}}}$, ${k_i} > 0$, ${k_d} > \frac{{{p_1} + {p_2}}}{{{b_0}}}$ for ${K_{\PID}}\left( s \right)$ to stabilize $P(s)$ if ${{b_0} = 0}$, and $\left( {{k_p},{k_i},{k_d}} \right) \in \Omega $ for ${K_{\PID}}\left( s \right)$ to stabilize $P(s)$ if ${\;{b_0} \ne 0}$.

\textit{Proof}: The proof follows analogously as in that of Theorem 3.1, and hence is omitted. $\hfill\blacksquare$

It is clear from Theorem 4.1 that for second-order systems the maximal gain and phase margins achievable by PID control are the same as those by general LTI controllers, provided that the plant is minimum-phase and has a relative degree no greater than one. On the other hand, if the plant does have a relative degree greater than one, then the maximal phase margin is reduced to ${\pi }/{2}$, despite the fact that the maximal gain margin remains unchanged.

%%%

%We can find from Theorem 4.1 that the maximal gain and phase margins achievable by PID control is in line with that by LTI control. In some particular cases when the system's relative degree is two, i.e., ${{\beta_0} = 0}$, we show that the bound on the phase margin is ${\pi }/{2}$. On the contrary, the maximal phase margin can achieve $\pi$ if the system's relative degree is one, in which ${{\beta_0} \ne 0}$. The case also implies that the minimum-phase zero $ - {\beta_1}/{\beta_0}$ is capable of furnishing extra phase for the system.

%%%%%%%%%%%%%%%%%%
\subsection{Nonminimum-phase Systems}
One of our primary results in this paper pertains to second-order unstable, nonminimum-phase plants. We consider specifically the plant described by
%Now we extend our results to second-order nonminimum-phase plants. We consider first the plant described by
%We first seek to determine in parallel the maximal gain and phase margins for plants that contain with two unstable poles, denoted by
%%%%%%%%%%%%%%%%%%
%\subsection{Complex Conjugate Poles}
%We now consider the plant with a pair of complex conjugate poles given below
\begin{equation} \label{2-49}
{\kern -5pt} P\left( s \right) \! = \!  \frac{{s - z}}{{\left( {s \! - \! {p_1}} \right)\left( {s \! - \!  {p_2}} \right)}},\;\;  {\mathop{\rm Re}\nolimits} \left( {{p_1}} \right) > 0, \;\; {\mathop{\rm Re}\nolimits} \left( {{p_2}} \right) > 0,
\end{equation}
where likewise, $z > 0$, ${p_1}$ and ${p_2}$ are either real or are a complex conjugate pair. To facilitate our development, we introduce the functions
%
%We first present the theorem for the plants with two real poles, which begins by defining
\begin{equation}\label{hattheta1}
{\kern 0pt}  \hat \theta \left( {{k_p}} \right) = {\tan ^{ - 1}}\frac{{( {{p_1} + {p_2}})\hat{\omega} ({k_p})}}{{{p_1}{p_2} - {\hat{\omega} ^2}({k_p})}}-
{\tan ^{ - 1}}\frac{{\hat{\omega} ({k_p})}}{z} - {\tan ^{ - 1}}\frac{{\hat{\omega} ({k_p})}}{{{k_p}}},
%{\kern -4 pt} \hat \theta \left( {{k_p}} \right)=
%\sum\limits_{i = 1}^2 {{{\tan }^{ - 1}}\frac{{{\omega(k_p)}}}{{{p_i}}}}  - {\tan ^{ - 1}}%\frac{{{\omega(k_p)}}}{z} - {\tan ^{ - 1}}\frac{\omega(k_p)}{k_p},
\end{equation}
and
\begin{equation}\label{hattheta-kd1}
\begin{array}{l}
{\kern -15pt}  \displaystyle \tilde \theta \left( {{k_d}} \right) = \pi  - {\tan ^{ - 1}}\frac{{( {{p_1} + {p_2}} )\tilde{\omega} ({k_d})}}{{{\tilde{\omega} ^2}({k_d}) - {p_1}{p_2}}}\\
{\kern 60pt}  \displaystyle - {\tan ^{ - 1}}\frac{{\tilde{\omega} ({k_d})}}{z} + {\tan ^{ - 1}}\frac{{z{k_d}\tilde{\omega} ({k_d})}}{{{p_1}{p_2}}},
\end{array}
%{\kern -4pt} \tilde \theta \left( {{k_d}} \right) = \pi   -  {\tan ^{ - 1}}\frac{{\left( {{p_1} \!+\! {p_2}} \right)\omega (\!{k_d}\!)}}{{{\omega ^2}(\!{k_d}\!) \! - \!{p_1}{p_2}}} - {\tan ^{ - 1}}\frac{{\omega (\!{k_d}\!)}}{z} + {\tan ^{ - 1}}\frac{{{k_d}\omega (\!{k_d}\!)}}{{{p_1}{p_2}/z}},
\end{equation}
where
\begin{equation}\label{hattheta2}
\hat{\omega}(k_p)=\sqrt {\frac{{p_1^2p_2^2 - {z^2}k_p^2}}{k_p^2 + {z^2} - p_1^2 - p_2^2}},
\end{equation}
\begin{equation}\label{hattheta-kd2}
\tilde{\omega} ({k_d}) = \sqrt {\frac{{{z^2}k_d^2 + {{\left( {{p_1}{p_2}/z} \right)}^2} - p_1^2 - p_2^2}}{{1 - k_d^2}}} .
\end{equation}
We first present the following theorem for plants with two real unstable poles.
%In the same light, we develop the maximal gain and phase margins in the following theorem. %The proof and the interpretations of this result share the spirit of those for Theorem 4.2 and hence are omitted.

\textit{Theorem 4.2} {\it(Real Poles)}: Let $P(s)$ be given by (\ref{2-49}), and $p_1,~p_2$ be real poles.  %and ${p_{\min }} = \min {\rm{\{ }}{p_1},{p_2}{\rm{\} }}$, ${p_{\max }} = \max {\rm{\{ }}{p_1},{p_2}{\rm{\} }}$.
For ${K_{\PD}}\left( s \right)$ and ${K_{\PID}}\left( s \right)$ to stabilize $P(s)$, it is necessary that $({p_1} - z)({p_2} - z) \ne 0$. Under this condition, the following statements hold:

(i). {\bf Gain Margin}: For $({p_1} - z)({p_2} - z)> 0$,
\begin{equation} \label{2-45-1}
{\kern -4 pt} \displaystyle k_\M^{\PID} = k_\M^{\PD} = \left\{ %{\begin{array}{*{20}{l}}
\begin{array}{ll}
\displaystyle {\kern -8 pt}  \frac{p_1p_2}{z({p_1} + {p_2})-z^2}, & \mbox{if~} z < \min \{ {p_1},{p_2}\},\\
%\displaystyle {\kern -11 pt}  \frac{{z({p_1} \!+\! {p_2}) \!-\! {z^2}}}{{{p_1}{p_2}}},\;\,{\rm{if}}\;\min \{ {p_1},{p_2}\}  %< z \le \sqrt {{p_1}{p_2}},\\
%\displaystyle {\kern -13 pt}  \frac{{z({p_1} \!+\! {p_2})\! - \!{p_1}{p_2}}}{{{z^2}}},{\rm{if}}\;\sqrt {{p_1}{p_2}}  < z < %\max \{ {p_1},{p_2}\},
%\end{array}\\
\displaystyle {\kern -8 pt}  \frac{z^2}{z({p_1} + {p_2}) - {p_1}{p_2}} & \mbox{if~} z > \max \{ {p_1},{p_2}\}.
\end{array}
%\end{array}}
\right.
\end{equation}
Otherwise, for $({p_1} - z)({p_2} - z)< 0$,
\begin{equation} \label{2-45-1a}
{\kern -2 pt} \displaystyle k_\M^{\PID} = k_\M^{\PD} = \left\{ {\begin{array}{*{20}{l}}
\displaystyle {\kern -8 pt}  \frac{{z({p_1} \!+\! {p_2}) \!-\! {z^2}}}{{{p_1}{p_2}}},\;\;\;\,{\rm{if}}\;\min \{ {p_1},{p_2}\}  < z \le \sqrt {{p_1}{p_2}},\\
\displaystyle {\kern -6 pt}  \frac{{z({p_1} \!+\! {p_2})\! - \!{p_1}{p_2}}}{{{z^2}}},{\rm{if}}\;\sqrt {{p_1}{p_2}}  < z < \max \{ {p_1},{p_2}\}.
\end{array}} \right.
\end{equation}
%(ii). For ${{p_1} + {p_2} > ({{{p_1}{p_2}}}/{z})  + z}$,
%\begin{equation} \label{2-45-2}
%{\kern -128pt}  k_\M^{\PID} = k_\M^{\PD} =  \displaystyle \frac{{z\left( {{p_1} + {p_2}} \right)}}{{{z^2} + {p_1}{p_2}}}.
%\end{equation}
%(ii) Whenever ${p_1} + {p_2} \ne \left( {{p_1}{p_2}/z} \right)  + z$,
%\begin{equation} \label{2-45+1}
%{\kern -180pt} \theta _\M^{\PID} = \theta _\M^{\PD}.
%\end{equation}

(ii). {\bf Phase Margin}: For $({p_1} - z)({p_2} - z)> 0$,
\begin{equation} \label{phase-realpoles}
\theta _\M^{\PID} = \theta _\M^{\PD} = \left\{
\begin{array}{ll}
|{\hat \theta ({{\hat k}_p})}|, & {\rm{if}}\;z < \min \{ {p_1},{p_2}\}, \\
\tilde \theta ({{\tilde k}_d}), & {\rm{if}}\;z > \max \{ {p_1},{p_2}\},
%{\kern -9 pt} \max \{ \hat \theta ({{\hat k}_{p0}}),| {\tilde \theta ({{\tilde k}_{d0}})} |\},\,{\rm{otherwise}}.
%\end{array}
\end{array} \right.
\end{equation}
where ${\hat k_p} \in ({p_1} + {p_2} - z,\;{p_1}{p_2}/z)$ is the unique positive solution to the polynomial equation
\begin{equation}\label{poly1}
%\begin{array}{l}
%{\kern -8 pt} zk_p^5 - {z^2}k_p^4 + z\left( {({p_1} + {p_2})(z - {p_1} - {p_2} - {p_1}{p_2}/z) + 2{p_1}{p_2}} \right)k_p^3\\
%{\kern 18 pt}  + 2p_1^2p_2^2k_p^2 + \left( {{p_1}{p_2}({p_1} + {p_2})(p_1^2 + p_2^2 - {p_1}{p_2}) + } \right.z{p_1}{p_2} \cdot \\
%{\kern 32 pt} \left. {({p_1}{p_2} - z({p_1} + {p_2}))} \right){k_p} + p_1^2p_2^2\left( {{z^2} - p_1^2 - p_2^2} \right) = 0,
%\end{array}
%{a_{15}}k_p^5 + {a_{14}}k_p^4 + {a_{13}}k_p^3 + {a_{12}}k_p^2 + {a_{11}}{k_p} + {a_{10}} = 0,
k_p^3 + {c_{2}}k_p^2 + {c_{1}}{k_p} + {c_{0}} = 0
\end{equation}
with coefficients
%${c_2} = {p_1} + {p_2} - z,{c_1} =  - ({p_1}{p_2}/z){c_2},{c_0} = ({z^2} - p_1^2 - p_2^2){p_1}{p_2}/z.$
\begin{align}
{\kern 0 pt}  &{c_{2}} ={p_1} + {p_2} - z,\notag \\
{\kern 0 pt}  &{c_{1}} = (z - {p_1} - {p_2}){p_1}{p_2}/z, \label{poly-kp-1} \\
{\kern 0 pt}  &{c_{0}} = ({z^2} - p_1^2 - p_2^2){p_1}{p_2}/z,\notag
\end{align}
%\begin{equation}\label{poly-kp-1}
%\begin{array}{l}
%{c_2} = {p_1} + {p_2} - z,~~{c_1} =  - ({p_1}{p_2}/z){c_2},\\
%{c_0} = ({z^2} - p_1^2 - p_2^2){p_1}{p_2}/z.
%\end{array}
%\end{equation}
%\begin{equation}\label{poly1}
%\begin{array}{*{20}{l}}
%{\kern -3 pt}  { zk^5_p- z^2k_p^4 + \left(z(p_1 + p_2)(z-p_1-p_2-p_1p_2)+2zp_1p_2\right) k_p^3 }\\
%{\kern 48 pt}  +2 p_1^2 p_2^2 k^2_p \\
%{\kern 48 pt} + \left(p_1p_2(p_1+p_2)(p^2_1+p^2_2-{p_1}{p_2})+zp_1p_2(p_1p_2-z(p_1+p_2))\right)k_p\\
%{\kern 92 pt}  + p^2_1p^2_2\left( z^2 - p_1^2 - p_2^2 \right)= 0.
%\end{array}
%\end{equation}
and ${{\tilde k}_d} \in ( - 1,\;({p_1}{p_2}/{z^2}) - (({p_1} + {p_2})/z))$ is the  unique negative solution to the polynomial equation
\begin{equation}\label{kd-poly1}
k_d^3 + {d_{2}}k_d^2 + {d_{1}}{k_d} + {d_{0}} = 0
%\begin{array}{l}
%{\kern -9 pt}  -\! {z^2}k_d^5 \!- \!{p_1}{p_2}k_d^4 \!+\! (p_1^2\! + \!p_2^2 \!+\! ({p_1} \!+\! {p_2})(z \!-\! {p_1}{p_2}/z))k_d^3 \!+ \! 2{p_1}{p_2}k_d^2 \\
%{\kern 24 pt} + \!(({p_1} \!+ \!{p_2})({({p_1}{p_2}/z)^2} \!- \!p_1^2 \!-\! p_2^2) \!+\! ({p_1} + {p_2} - {p_1}{p_2}/z) \cdot \\
%{\kern 45 pt} \left. {{p_1}{p_2}/z} \right){k_d} + {p_1}{p_2}({({p_1}{p_2}/z)^2} - p_1^2 - p_2^2)/{z^2} = 0.
%\end{array}
\end{equation}
with coefficients
\begin{equation}\label{poly-kd-1}
\begin{array}{l}
{d_{2}} = {d_{1}} = ({p_1}{p_2} - ({p_1} + {p_2})z)/{z^2},\\
{d_{0}} = (p_1^2 + p_2^2 - {({p_1}{p_2}/z)^2})/{z^2}.
\end{array}
\end{equation}
Otherwise, for $({p_1} - z)({p_2} - z)< 0$,
\begin{equation} \label{phase-realpoles1a}
{\kern -6 pt} \theta _\M^{\PID} = \theta _\M^{\PD} = \max \{ \hat \theta ({{\hat k}_{p}}),| {\tilde \theta ({{\tilde k}_{d}})} |\},
\end{equation}
where ${{\hat k}_{p}} \in ({p_1}{p_2}/z,\;{p_1} + {p_2} - z)$ is the unique solution to the equation (\ref{poly1}), and ${{\tilde k}_{d}} \in (({p_1}{p_2}/{z^2}) - (({p_1} + {p_2})/z),\; - 1)$ is the unique solution to the equation (\ref{kd-poly1}).

\textit{Proof}: See Appendix \ref{Theorem 4.2-Proof}. $\hfill\blacksquare$

The following theorem is a counterpart to Theorem 4.2 for plants with complex
conjugate poles. While the result is essentially similar, subtle differences
do exist. For example, for a pair of complex conjugate poles
$p_1=p$ and $p_2=p^*$, it always holds that $({p_1} - z)({p_2} - z)=|p-z|^2>0$,
while on the contrary the case $({p_1} - z)({p_2} - z)<0$ is precluded.
For this reason, we state the result separately.

\textit{Theorem 4.3} {\it(Complex Poles)}: Let $P(s)$ be given by (\ref{2-49}),
 and ${p_1} = p = \sigma  + j\nu$, ${p_2} = {p^*} = \sigma  - j\nu$ with $\sigma  > 0$.
 %For ${K_{\PD}}\left( s \right)$ and ${K_{\PID}}\left( s \right)$ to stabilize $P(s)$, it is necessary that $|p-z|>0$. Under %this condition,
 Then, the following statements hold:

(i). {\bf Gain Margin}: For $\left| {p - z} \right| < z$,
\begin{equation} \label{gain-complexpoles1}
{\kern -9 pt}  \displaystyle k_\M^{\PID} = k_\M^{\PD} = \left\{ {\begin{array}{*{20}{l}}
 \displaystyle {\frac{{{p_1}{p_2}}}{{z({p_1} \!+\! {p_2}) \!-\! {z^2}}},\;\;\;\;\;\;\;\;\;\;{\rm{if}}\;z \le \left| p \right|,}\\
 \displaystyle {\frac{{{z^2}}}{{z({p_1} \!+\! {p_2}) \!-\! {p_1}{p_2}}},\;\;\;\;\;\;\;{\rm{if}}\;z > \left| p \right|.}
\end{array}} \right.
\end{equation}
Otherwise, for $\left| {p - z} \right| \ge z$,
\begin{equation} \label{gain-complexpoles2}
{\kern -9 pt} \displaystyle  k_\M^{\PID} = k_\M^{\PD} = \left\{ {\begin{array}{*{20}{l}}
\displaystyle  {\frac{{{p_1}{p_2}}}{{z({p_1} + {p_2}) - {z^2}}},\;\;\,{\rm{if}}\;z < {p_1} + {p_2},}\\
\displaystyle  {\frac{{{p_1}{p_2} + {z^2}}}{{z({p_1} + {p_2})}},\;\;\;\;\;\;\;\;\;\;\,{\rm{if}}\;z \ge {p_1} + {p_2}.}
\end{array}} \right.
\end{equation}

(ii). {\bf Phase Margin}: For $\left| {p - z} \right| < z$,
\begin{equation} \label{phase-conplexpoles1}
{\kern -12 pt}\theta _\M^{\PID} = \theta _\M^{\PD} = \left\{ {\begin{array}{*{20}{l}}
{|\hat \theta ({{\hat k}_p})|,\;\;\;\;\;\;\;\;\;\;\;\;\;{\rm{if}}\;z \le \left| p \right|,}\\
{\tilde \theta ({{\tilde k}_d}),\;\;\;\;\;\;\;\;\;\;\;\;\;\;\;{\rm{if}}\;z > \left| p \right|.}
\end{array}} \right.
\end{equation}
Otherwise, for $\left| {p - z} \right| \ge z$,
\begin{equation} \label{phase-conplexpoles2}
\theta _\M^{\PID} = \theta _\M^{\PD} = |\hat \theta ({{\hat k}_p})|,
\end{equation}
where ${\hat k_p} \in ({p_1} + {p_2} - z,\;{p_1}{p_2}/z)$ is the unique solution to the polynomial equation (\ref{poly1}), and ${{\tilde k}_d} \in ( - 1,\;({p_1}{p_2}/{z^2}) - (({p_1} + {p_2})/z))$ is the unique solution to the polynomial equation (\ref{kd-poly1}).
%
%In particular,
%\begin{equation} \label{2-45+3}
%{{\hat \theta }_1}\left( {{{\bar k}_p}} \right) = {\tan ^{ - 1}}\frac{{\left( {{p_1} + {p_2}} \right){{\tilde \omega }^*}}}{{{p_1}{p_2} - {{\tilde \omega }^{*2}}}} - {\tan %^{ - 1}}\frac{{\left( {{{\bar k}_p} + z} \right){{\tilde \omega }^*}}}{{z{{\bar k}_p} - {{\tilde \omega }^{*2}}}},
%\end{equation}
%and
%\begin{equation} \label{2-45+4}
%{{\hat \theta }_2}\left( {{{\bar k}_p}} \right) = \pi  - {\tan ^{ - 1}}\frac{{\left( {{p_1} + {p_2}} \right){{\tilde \omega }^*}}}{{{{\tilde \omega }^{*2}} - {p_1}{p_2}}} - %{\tan ^{ - 1}}\frac{{\left( {{{\bar k}_p} + z} \right){{\tilde \omega }^*}}}{{z{{\bar k}_p} - {{\tilde \omega }^{*2}}}},
%\end{equation}
%solution of the following equation

\textit{Proof}: See Appendix \ref{Theorem 4.3-Proof}. $\hfill\blacksquare$

\textit{Remark 4.1}: It is worth noting that, while given in terms of the positive real roots of polynomials, both $\hat{k}_p$ and $\tilde{k}_d$ can be calculated explicitly using the classical {\em Cardano formula} \cite{borwein1995polynomials}. In this sense, the expressions for the phase margin
(i.e., (\ref{phase-realpoles})-(\ref{phase-realpoles1a}), and (\ref{phase-conplexpoles1})-(\ref{phase-conplexpoles1}))
 provide in fact analytical solutions to the phase maximization problem. Furthermore,
 note also that because of the uniqueness of the solutions, the functions $\hat{\theta}(k_p)$ and $\tilde{\theta}(k_d)$ define, respectively, unimodal convex and concave functions on the corresponding intervals, a fact further evidenced by the following example.

\textit{Example 4.1}: In this example, we compare the maximal margins achievable by PID controllers and those by general LTI controllers. We consider plants with both real poles and complex conjugate poles specified as follows.

%The result thus highlights the essential role of the derivative control in stabilization. Note also that similar to (\ref{new3-2}), the maximal gain and phase margins exhibit the relation
%\begin{equation} \label{PID4-1}
%\theta _\M^\p = {\cos ^{ - 1}}\frac{{2\sqrt {k_\M^\p} }}{{1 + k_\M^\p}}.
%\end{equation}

%In light of Corollary 4.1 and Theorem 4.3, Theorem 4.2, we assert that the involvement of the integral term, i.e., ${k_i} \ne 0$ will make the gain and phase margins smaller. It can be drawn analogously the conclusion that the maximal gain and phase margins by PID control is equivalent to that by PD control. This highlights the essential role of the PD control in stabilization.

%\textit{Example 4.1}: In this example, similarly, we draw a companion between the Towards this end, it was found in \cite{chen2000logarithmic} that
%\begin{equation} \nonumber
%{\gamma _{{\rm{opt}}}} = \prod\limits_{i = 1}^2 {\left| {\frac{{{p_i} + z}}{{{p_i} - z}}} \right|} .
%\end{equation}
%We consider the second-order system (\ref{2-49}) with two real poles and a pair of unstable complex conjugate poles, respectively.

\begin{figure}[!htb]
\centering
\includegraphics[width=5.80cm]{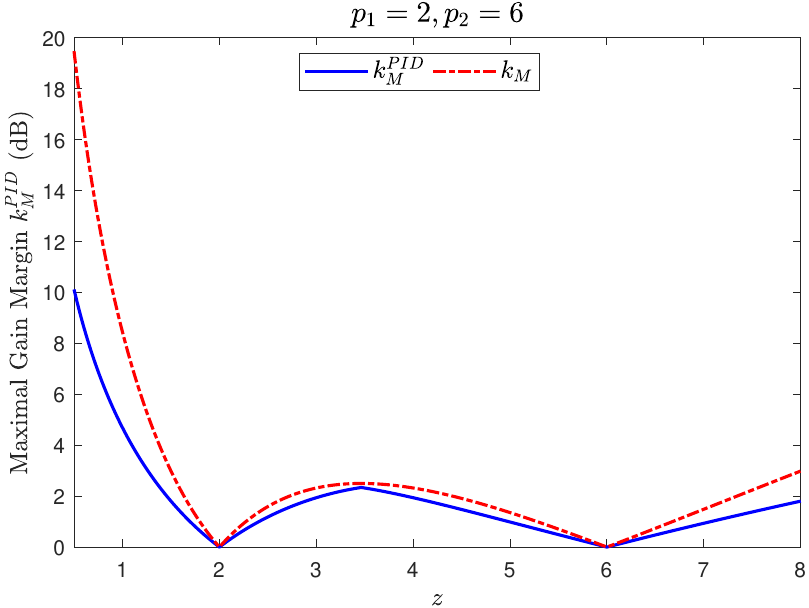}
\caption{Maximal gain margin $k_\M^{\PID}$ of the system (\ref{2-49}) vs. ${k_\M}$ : Real poles.}
\label{Fig4-secondorder_gain_real pole}
\end{figure}

\begin{figure}[!htb]
\centering
\includegraphics[width=5.80cm]{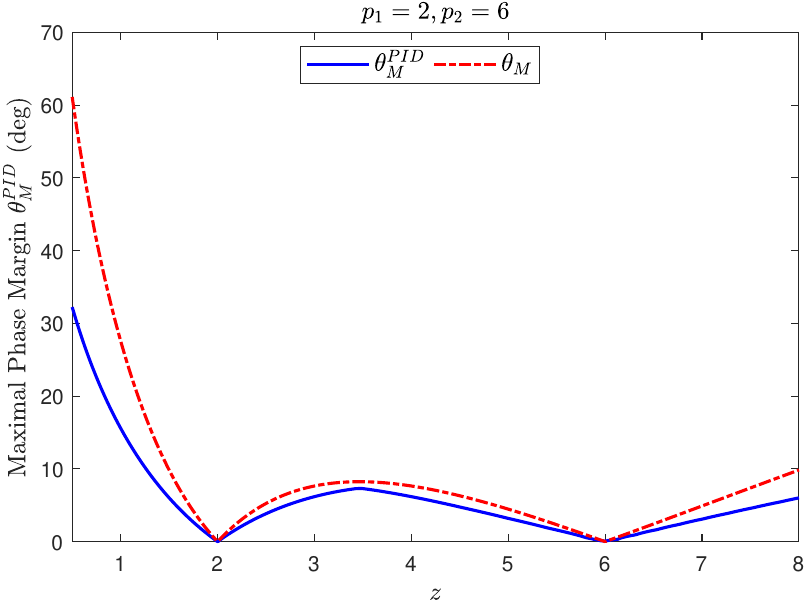}
\caption{Maximal phase margin $\theta_\M^{\PID}$ of the system (\ref{2-49}) vs. ${\theta_\M}$ : Real poles.}
\label{Fig4-secondorder_phase_real pole}
\end{figure}

\textit{Real unstable poles} ${p_1} > 0, {p_2} > 0$: In this case we take ${p_1} = 2$, ${p_2} = 6$ and let $z$ vary in the interval $\left[ {0.5,{\rm{ }}8} \right]$. Fig. \ref{Fig4-secondorder_gain_real pole} plots the gain margins $20{\log _{10}}k_\M^{\PID}$ and $20{\log _{10}}{k_\M}$ (both in db) as functions of $z$, and Fig. \ref{Fig4-secondorder_phase_real pole} presents the phase margins $\theta_\M^{\PID}$ and $\theta_\M$ (in deg) as functions of $z$. Note that, when $z= p_1$ or $z= p_2$, the gain and phase margins both vanish.

\begin{figure}[!htb]
\centering
\includegraphics[width=5.80cm]{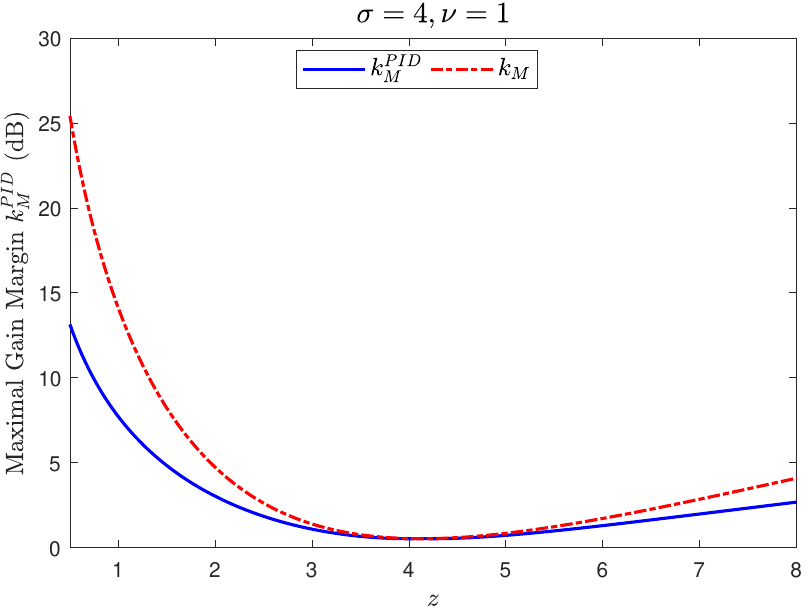}
\caption{Maximal gain margin $k_\M^{\PID}$ of the system (\ref{2-49}) vs. ${k_\M}$ : Complex poles.}
\label{Fig4-secondorder_gain_complex pole}
\end{figure}

\begin{figure}[!htb]
\centering
\includegraphics[width=5.80cm]{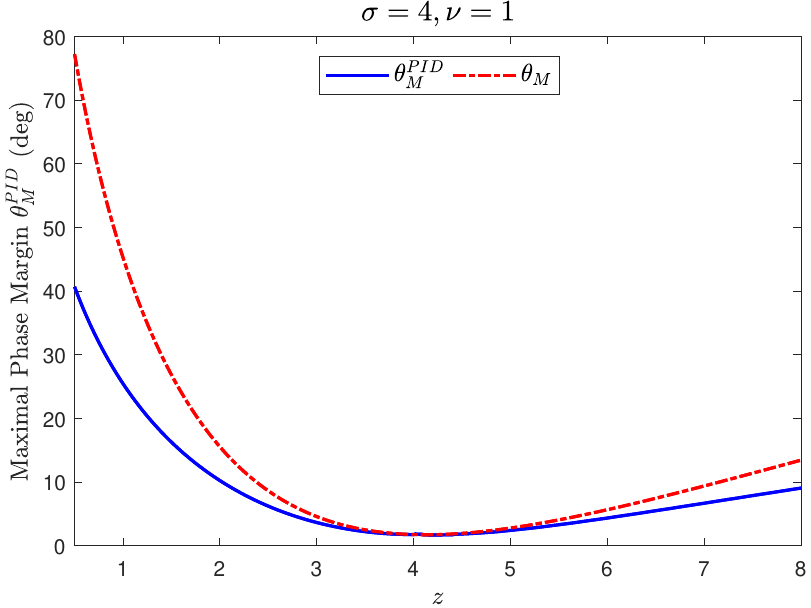}
\caption{Maximal phase margin $\theta_\M^{\PID}$ of the system (\ref{2-49}) vs. ${\theta_\M}$ : Complex poles.}
\label{Fig4-secondorder_phase_complex pole}
\end{figure}

\textit{Complex conjugate unstable poles} ${p_1} = \sigma  + j\nu ,{p_2} = \sigma  - j\nu$:
%$\sigma > 0$: %In an anlogous vein
We take $\sigma = 4$, $\nu = 1$ and let $z$ vary in the interval $\left[ {0.5,{\rm{ }}8} \right]$. Fig.   \ref{Fig4-secondorder_gain_complex pole} shows the gain margins $20{\log _{10}}k_\M^{\PID}$ and $20{\log _{10}}{k_\M}$, %(both in dB) %$k_\M^{\PID}$ and $k_\M$
%as a function of $z$,
while Fig. \ref{Fig4-secondorder_phase_complex pole} describes the phase margins $\theta_\M^{\PID}$ and $\theta_\M$.

\begin{figure}[!htb]
\centering
\includegraphics[width=5.80cm]{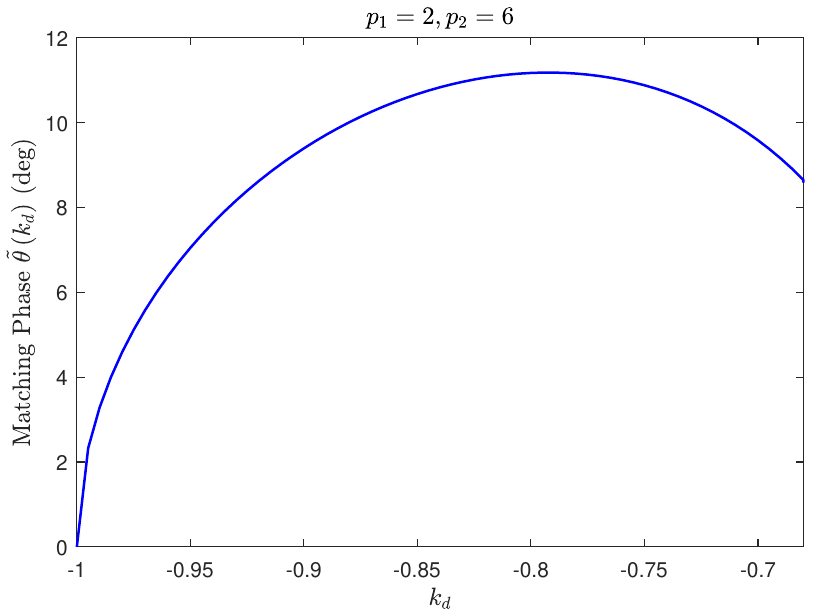}
\caption{The phase $\tilde \theta \left( {{k_d}} \right)$ of the system (\ref{2-49}): $z > \max {\rm{\{ }}{p_1},{p_2}{\rm{\} }}$.
%${{p_1} \!+\! {p_2} > {p_1}{p_2}/z}$.
}
\label{Fig4-matching-phase-realpoles}
\end{figure}

\begin{figure}[!htb]
\centering
\includegraphics[width=5.80cm]{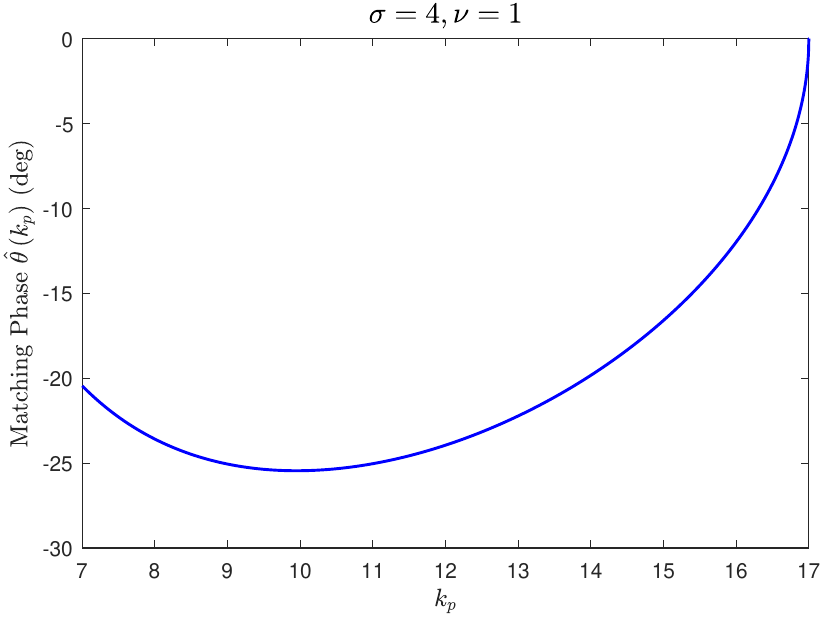}
\caption{The phase $\hat \theta \left( {{k_p}} \right)$ of the system (\ref{2-49}): $\left| {p - z} \right| \ge z$.
%${{p_1}\! +\! {p_2} \le {p_1}{p_2}/z}$.
}
\label{Fig4-matching-phase-complexpoles}
\end{figure}

To facilitate the understanding of Theorem 4.2-4.3 and their proofs, we also plot $\hat \theta \left( {{k_p}} \right)$ and $\tilde \theta \left( {{k_d}} \right)$. %where for the case of two real poles ${p_1} + {p_2} \le z$, and for the case of complex poles, ${p_1} + {p_2} > z$.
Fig. \ref{Fig4-matching-phase-realpoles} corresponds to the case $z > \max {\rm{\{ }}{p_1},{p_2}{\rm{\} }}$
%${p_1} \!+ \!{p_2} \!-({p_1}{p_2}/z) = 6.8 > 0$
(herein $z = 10$), where we observe that $\tilde \theta \left( {{k_d}} \right) > 0$ in $\left( { - 1,-0.68} \right)$, where $\tilde{\theta}(k_d)$ is concave and has a unique maximum. Fig. \ref{Fig4-matching-phase-complexpoles} demonstrates that for $\left| {p - z} \right| \ge z$
%${p_1} \!+ \!{p_2} \! - ({p_1}{p_2}/z) =  - 9 < 0$
(herein $z = 1$), $\hat \theta \left( {{k_p}} \right) < 0$ in $\left( {7,17} \right)$, indicating that in this interval, $\hat \theta \left( {{k_p}} \right)$ is convex and admits a unique minimum. %the polynomial (\ref{poly1}) admits a unique solution.
$\hfill\blacksquare$

Interestingly, by a careful inspection, we observe further from the gain margin plots that the ratio ${\log _{10}}{k_\M}/{\log _{10}}k_\M^{\PID}$ seemingly is within a factor of two. This, in fact, turns out not to be a coincidence. The following corollary shows that, as with first-order plants,  ${\log _{10}}{k_\M}$ is always within a factor of two of ${\log _{10}}k_\M^{\PID}$.

\textit{Corollary 4.1}: Let $P(s)$ be given by (\ref{2-49}). Then,
\begin{equation} \label{corollary 4.1}
%\begin{array}{l}
{\log _{10}}k_\M^{\PID} \le {\log _{10}}{k_\M} \le 2{\log _{10}}k_\M^{\PID}.
%k_\M^{\PID} \le {k_\M} \le {\left( {k_\M^{\PID}} \right)^2}.
%\\
%\theta _\M^{\PID} \le {\theta _\M} \le 2\theta _\M^{\PID}.
%\end{array}
\end{equation}

\textit{Proof}: We consider the case of real poles. The proof for the case of complex poles follows
analogously and hence is omitted. Note from \cite{chen2000logarithmic} that for $P(s)$ given by (\ref{2-49}),
\begin{equation} \nonumber
{\gamma _{{\rm{opt}}}} = \prod\limits_{i = 1}^2 {\left| {\frac{{{p_i} + z}}{{{p_i} - z}}} \right|},
\end{equation}
which can be written explicitly as
$$
\gamma_{opt}=\displaystyle{\left\{\begin{array}{ll}
\frac{(p_1+z)(p_2+z)}{(z-p_1)(p_2-z)} & \mbox{~if ~} p_1<z<p_2, \\
\frac{(p_1+z)(p_2+z)}{(p_1-z)(p_2-z)} & \mbox{~otherwise.}\end{array}\right.
}
$$
%Note that the second case in this expression includes that of complex poles.
Accordingly,
$$
\sqrt{k_\M} =  \frac{\gamma_{opt}+1}{\gamma_{opt}-1}=\displaystyle{\left\{\begin{array}{ll}
\frac{z(p_1+p_2)}{z^2+p_1p_2} & \mbox{~if ~} p_1<z<p_2, \\
\frac{z^2+p_1p_2}{z(p_1+p_2)} & \mbox{~otherwise.}\end{array}\right.
}
$$
Note also that under the condition $({p_1} - z)({p_2} - z) < 0$, it holds that
$$
\frac{{z({p_1} + {p_2})}}{{{z^2} + {p_1}{p_2}}} \le \frac{{z({p_1} + {p_2}) - {z^2}}}{{{p_1}{p_2}}},
$$
$$
\frac{{z({p_1} + {p_2})}}{{{z^2} + {p_1}{p_2}}} \le \frac{{z({p_1} + {p_2}) - {p_1}{p_2}}}{{{z^2}}}.
$$
%$$
%\frac{{z({p_1} + {p_2})}}{{{z^2} + {p_1}{p_2}}} \le \frac{{{z^2} + {p_1}{p_2}}}{{z({p_1} + {p_2})}}.
%$$
On the other hand, it is straightforward to verify that
%Under the condition $p_1<z<p_2$, it holds that $p_1+p_2>p_1p_2/z$.
%It is then straightforward to verify that
%$$
%\sqrt{k_\M}=\frac{z(p_1+p_2)}{z^2+p_1p_2}\leq \frac{z^2}{z(p_1+p_2)-p_1p_2}=k^{\PID}_M.
%$$
%Otherwise, it can also be readily shown that
$$
\frac{z^2+p_1p_2}{z(p_1+p_2)}\leq \frac{z^2}{z(p_1+p_2)-p_1p_2},
$$
and
$$
\frac{z^2+p_1p_2}{z(p_1+p_2)}\leq \frac{p_1p_2}{z(p_1+p_2)-z^2},
$$
provided that $({p_1} - z)({p_2} - z) > 0$, or equivalently, $p_1+p_2-z<p_1p_2/z$.
Thus, in all cases, it holds that $\sqrt{k_\M}\leq k^{\PID}_\M$.
This concludes the proof for the inequalities in (\ref{corollary 4.1}).
$\hfill\blacksquare$

\textit{Remark 4.2}: From the numerical results in Fig. 3 and Fig. 5, it appears plausible to contend that a similar relationship holds for the phase margins $\theta _\M^{\PID}$ and $\theta _M$; that is,
%In fact, it follows from the numerical computations that similar conclusion may hold for the maximal phase margin as well, that is
\begin{equation}\label{Phase2Times}
\theta _\M^{\PID} \le {\theta _\M} \le 2\theta _\M^{\PID}.
\end{equation}
While a general proof is currently unavailable, it can be shown that under
%However, this statement cannot be ascertained explicitly in its full interval in general. Under
certain circumstances, the inequality does hold.
%it can be determined by analysing $\hat \theta ({k_p})$ or $\tilde \theta ({k_d})$ at its boundary value.
Consider for example, the case $z <\min \{ {p_1},{p_2}\} $. %${p_1} + {p_2} \le {p_1}{p_2}/z$.
%In this case, $p_1>z$, $p_2>z$.
It follows instantly from
Proposition 2.1 that
$$
{\tan}\frac{{{\theta _\M}}}{2} = \frac{{({p_1} - z)({p_2} - z)}}{{2\sqrt {z({p_1} + {p_2})({p_1}{p_2} - {z^2})} }}.
$$
Hence, to establish the inequality ${\theta _\M} \le 2\theta _\M^{\PID}$, it is both necessary and sufficient to show that
$$
{\tan}\frac{\theta _\M}{2} \le {\tan}|\hat \theta ({\hat k_p})|.
$$
We evaluate instead $\hat{\theta}(p_1+p_2-z)$, and find conditions such that $\theta_\M\leq 2|\hat{\theta}(p_1+p_2-z)|$. By algebraic manipulation, we obtain
%It is clear that ${\omega ^2}({k_p}) = ({p_1}{p_2} + z({p_1} + {p_2} - z))/2$ at ${k_p} = {p_1} + {p_2} - z$. Noting that
%$$
%{\omega ^2}({k_p}) - {k_p}z = {p_1}{p_2} - {\omega ^2}({k_p}) = ({p_1}{p_2} - z({p_1} + {p_2} - z))/2 > 0
%$$
\begin{equation}\nonumber
\begin{array}{l}
\displaystyle {\kern -8 pt}  |\hat \theta ({k_p})| = \pi  - {\tan ^{ - 1}}\frac{{({p_1} + {p_2})\hat{\omega} ({k_p})}}{{{p_1}{p_2} - {\hat{\omega} ^2}({k_p})}} - {\tan ^{ - 1}}\frac{{({k_p} + z)\hat{\omega} ({k_p})}}{{{\hat{\omega} ^2}({k_p}) - {k_p}z}}.
\end{array}
\end{equation}
At ${k_p} = {p_1} + {p_2} - z$, we have
$$
|\hat \theta (p_1+p_2-z)| ={\tan ^{ - 1}}\frac{({p_1} + {p_2})({p_1}{p_2} - {k_p}z)\bar{\omega}}{({p_1} + {p_2})^2\bar{\omega}^2 - ({p_1}{p_2} - \bar{\omega}^2)^2},
$$
where
\begin{equation}\label{omegabar}
\bar{\omega}=\hat{\omega}(p_1+p_2-z)=\sqrt{\frac{p_1p_2+z(p_1+p_2-z)}{2}~},
\end{equation}
and ${p_1}{p_2} - {k_p}z=(p_1-z)(p_2-z)$.
To show that $\theta_\M\leq 2 |\hat{\theta}(p_1+p_2-z)|$, it suffices to prove that $\tan(\theta_\M/2)\leq \tan |\hat{\theta}(p_1+p_2-z)|$, which is
equivalent to the condition
\begin{equation}\label{PID-LTI_Phase1}
\begin{array}{l}
{\kern -25 pt}  2({p_1} + {p_2})\bar{\omega}\sqrt {z({p_1} + {p_2})({p_1}{p_2} - {z^2})}  \\
{\kern  25 pt}  \geq ({p_1} + {p_2})^2 \bar{\omega} ^2 - ({p_1}{p_2} - \bar{\omega}^2)^2.
%{({p_1}{p_2} - z({p_1} + {p_2} - z))^2}/4.
\end{array}
\end{equation}
It can be readily verified, however, that this inequality holds whenever
%Indeed, we claim that the inequality above holds for
\begin{equation}\label{PID-LTI_Phase2}
\frac{{{p_1} + {p_2}}}{z} \le 4\frac{{{p_1}{p_2} + {z^2}}}{{{p_1}{p_2} + z({p_1} + {p_2} - z)}}.
\end{equation}
Consequently, in the case $z < \min \{ {p_1},{p_2}\} $, under the condition (\ref{PID-LTI_Phase2}), we have found that
$$
\theta_\M \leq 2 |\hat \theta (p_1+p_2-z)|\leq 2 \theta_\M^{\PID}.
$$
In fact, since $|\hat \theta (p_1+p_2-z)|$ can be significantly smaller than $\theta_\M^{\PID}$, the analysis reveals, even more strikingly, that
$\theta_\M^{\PID}$ can be much closer to $\theta_\M$ than within a factor of two.
%Toward this end, we denote
%$$
%A({p_i},z) = \sqrt {z({p_1} + {p_2})({p_1}{p_2} - {z^2})},
%$$
%$$
%B({p_i},z) = ({p_1}{p_2} - z({p_1} + {p_2} - z))/2,
%$$
%$$
%C({p_i},z) = ({p_1} + {p_2})\sqrt {({p_1}{p_2} + z({p_1} + {p_2} - z))}.
%$$
%Noting also that $C({p_i},z) > A({p_i},z)$, $C({p_i},z) \gg  B({p_i},z)$ for ${p_1} + {p_2} \le {p_1}{p_2}/z$. As such, the inequality (\ref{PID-LTI_Phase1}) becomes
%\begin{equation}\label{PID-LTI_Phase3}
%0 < {C^2}({p_i},z) - {B^2}({p_i},z) \le 2A({p_i},z)C({p_i},z).
%\end{equation}
%Clearly, this inequality holds for $A({p_i},z) < C({p_i},z) \le \sqrt 2 A({p_i},z)$, which subsequently results in the inequality (\ref{PID-LTI_Phase2}). This allows us to %conclude that when the condition in (\ref{PID-LTI_Phase2}) is fulfilled, we have $\theta _\M^{\PID} \le {\theta _\M} \le 2\theta _\M^{\PID}.$
$\hfill\blacksquare$

Equally of interest, the following corollary suggests that in the presence of a
nonminimum-phase zero, the maximal phase margin attainable by PID control is
limited to $90^0$, thus also highlighting the limitation of PID controllers.

\textit{Corollary 4.2}: Let $P(s)$ be given by (\ref{2-49}). Then,
\begin{equation} \label{corollary 4.2}
\theta^{\PID}_\M \leq \frac{\pi}{2}.
\end{equation}

\textit{Proof}: We prove the corollary for the case $z <\min \{ {p_1},{p_2}\} $; the
proof for other cases is analogous and hence omitted. In this vein,
%$p_1+p_2\leq p_1p_2/z$,
we note that if $\hat{\omega}^2(k_p)\leq k_p z$, then
\begin{equation}\nonumber
|\hat \theta \left( {{k_p}} \right)| = {\tan ^{ - 1}}\frac{{\left( {{k_p} + z} \right)\hat \omega }}{{{k_p}z - {{\hat \omega }^2}}} - {\tan ^{ - 1}}\frac{{\left( {{p_1} + {p_2}} \right)\hat \omega }}{{{p_1}{p_2} - {{\hat \omega }^2}}} \le \frac{\pi }{2}.
\end{equation}
Otherwise, if $\hat{\omega}^2(k_p) > {k_p}z$, then
\begin{equation}\nonumber
|\hat \theta \left( {{k_p}} \right)| =\pi  - {\tan ^{ - 1}}\frac{{({p_1} + {p_2})\hat\omega }}{{{p_1}{p_2} - {\hat\omega ^2}}} - {\tan ^{ - 1}}\frac{{({k_p} + z)\hat\omega }}{{{\hat\omega ^2} - {k_p}z}}.
\end{equation}
For $k_p>p_1+p_2-z$, it is straightforward to see that $k_p(p_1+p_2-z)+z(p_1+p_2)-p_1p_2>0$.
Thus, for any $\hat{\omega}(k_p)>0$,
$$
\hat{\omega}^4+\left(k_p(p_1+p_2-z)+z(p_1+p_2)-p_1p_2\right)\hat{\omega}^2+k_p zp_1p_2>0,
$$
which is equivalent to the inequality
$$
\left( {\frac{{\left( {{k_p} + z} \right)\hat \omega }}{{{{\hat \omega }^2} - {k_p}z}}} \right)\left( {\frac{{({p_1} + {p_2})\hat \omega }}{{{p_1}{p_2} - {{\hat \omega }^2}}}} \right) > 1.
$$
In view of Lemma A.4, we may rewrite $|\hat \theta \left( {{k_p}} \right)|$ as
$$
|\hat \theta \left( {{k_p}} \right)|=\tan^{-1} \displaystyle{
\frac{\left(\frac{{\left( k_p+z \right){{ \hat\omega }}}}{{\hat\omega^2-{k_p}z }}\right)
+\left(\frac{{({p_1} + {p_2})\hat\omega }}{{{p_1}{p_2} - {\hat\omega ^2}}}\right)}{\left(\frac{{\left( k_p+z \right){{ \hat\omega }}}}{{\hat\omega^2}-{k_p}z }\right)
\left(\frac{{({p_1} + {p_2})\hat\omega }}{{{p_1}{p_2} - {\hat\omega ^2}}}\right)-1}\leq \frac{\pi}{2}.
}
$$
Hence, whenever $z <\min \{ {p_1},{p_2}\} $,
%$p_1+p_2\leq p_1p_2/z$,
we have $|\hat \theta \left( {{k_p}} \right)|\leq \pi/2$.
%The proof for the case of $z > \max \{ {p_1},{p_2}\} $,
%$p_1+p_2> p_1p_2/z$
%follows similarly and hence is omitted.
$\hfill\blacksquare$
%For any
%$k_d<0$,
%$$
%\tilde \theta \left( {{k_d}} \right) \leq \pi  - {\tan ^{ - 1}}\frac{{\left( {{p_1} + {p_2}} \right)\tilde{\omega} ({k_d})}}{{{\tilde{\omega} ^2}({k_d}) - {p_1}{p_2}}}.
%$$
%Note also that $\tilde{\omega}(k_d)$ is a monotonically decreasing function of $k_d$ on $(-1,~0)$, and $\tilde{\omega}(0)
%From the proof of Theorem 5.1, we know that the upper bound is a monotonically increasing function of
%\vspace{-0.25 cm}
%\vspace{-0.25 cm}
%\begin{figure}[!htb]
%\centering
%\includegraphics[width=6.80cm]{secondorder_gain_real and complex pole.pdf}
%\caption{Maximal gain margin $k_\M^{\PID}$ of the system (\ref{2-49}) vs. ${k_\M}$ : a). Real poles, b). Complex poles.}
%\label{Fig3-secondorder_gain_real and complex pole}
%\end{figure}
%
%\vspace{-0.25 cm}
%\begin{figure}[!htb]
%\centering
%\includegraphics[width=6.80cm]{secondorder_phase_real and complex pole.pdf}
%\caption{Maximal phase margin $\theta_\M^{\PID}$ of the system (\ref{2-49}) vs. ${\theta_\M}$ : a). Real poles, b). Complex poles.}
%\label{Fig3-secondorder_phase_real and complex pole}
%\end{figure}

%\textit{Example 4.2}: In this example, herein we show the matching phases $\hat \theta \left( {{k_p}} \right)$ in two cases, i.e., ${p_1} + {p_2} > z$, and ${p_1} + {p_2} \le z$ for $p_1 > z, p_2 > z$. We take ${p_1} = 6$, ${p_2} = 8$, and we let $z = 12$,  $z =32$, respectively. Fig.\ref{Fig4-matching-phase} gives the corresponding matching phases $\hat \theta \left( {{k_p}} \right)$, which illustrates the properties of $\hat \theta \left( {{k_p}} \right)$ in these two cases.
%\vspace{-0.25 cm}

\section{PI Control}
%\section{Extensions and Discussions: PI Control}
%\section{PI Control and Special Cases}

In the preceding section, we have shown that the maximal phase margin can be determined by solving two third-order polynomials. In this section, we
derive a simple closed-form analytical expression of the maximal phase margin, requiring no solution of polynomial equations.
By virtue of this development, it becomes clear that the derivative
control provides a strict improvement in increasing the gain and phase margins. On the other hand, it is also seen that this improvement is at best
limited to a factor of two.

%which solutions of polynomial equations are no longer needed. The results too reinforce the proceeding conclusion that the integral gain has no effect  in improving the gain and phase margins.

\textit{Theorem 5.1}: Let $P(s)$ be given by (\ref{2-49}). % For ${K_{\PID}}\left( s \right)$ to stabilize $P(s)$, it is necessary and sufficient that ${k_p} < \frac{{{k_i} + {p_1}{p_2}}}{z}$, ${k_i} < 0$, $- 1 < {k_d} < \frac{{{k_p} - {p_1} - {p_2}}}{z}$.
For ${K_{\PI}}\left( s \right)$ to stabilize $P(s)$, it is necessary that ${p_1} + {p_2} < {p_1}{p_2}/z$.  Under this condition, the following statements hold:

(i).
\begin{equation} \label{2-44}
{k_\M^{\PI} = k_\M^\p = \frac{{{p_1}{p_2}}}{{z\left( {{p_1} + {p_2}} \right)}}.}
\end{equation}

(ii).
\begin{equation} \label{2-44-phase}
%{\theta _\M^{\PI} = \theta _\M^\p = {{\tan }^{ - 1}}\frac{{{\omega _0}}}{z} - {{\tan }^{ - 1}}\frac{{{\omega _0}}}{{{p_1}}} - {{\tan }^{ - 1}}\frac{{{\omega _0}}}{{{p_2}}},}
\theta _\M^{\PI} = \theta _\M^\p = \bar{\theta}(\omega_0),
%
%{\tan ^{ - 1}}\frac{{{\omega _0}}}{z} - {\tan ^{ - 1}}\frac{{\left( {{p_1} + {p_2}} \right){\omega _0}}}{{{p_1}{p_2} - \omega _0^2}},
\end{equation}
where
$$
{\bar \theta (\omega ) = {{\tan }^{ - 1}}\frac{\omega }{z} - {{\tan }^{ - 1}}\frac{{({p_1} + {p_2})\omega }}{{{p_1}{p_2} - {\omega ^2}}},}
$$
%\begin{equation}
\begin{equation}
\begin{array}{l}
{\kern -4pt}  \displaystyle  {\omega_0} = \left\{ {\frac{{z({p_1} + {p_2})}}{{2({p_1} + {p_2} \!-\! z)}}\left( {\frac{{p_1^2 + p_2^2}}{{{p_1} + {p_2}}} \!-\! z \!-\! \frac{{{p_1}{p_2}}}{z} + \left( {{{(z \!- \!\frac{{{p_1}{p_2}}}{z})}^2}} \right.} \right.} \right. \\ \\
{\kern 20pt}  \displaystyle  {\left. {\left. {{{\left. { + {{({p_1} \!-\! {p_2})}^2}(1 - 2\frac{{z + {p_1}{p_2}/z}}{{{p_1} + {p_2}}})} \right)}^{\frac{1}{2}}}} \right)} \right\}^{\frac{1}{2}}}.\label{PI-phase-analytical}
\end{array}
\end{equation}
%\begin{array}{l}
%{\kern -6pt}  \displaystyle  {\omega _0} = \left\{ {\frac{{z({p_1} + {p_2})}}{{2({p_1} + {p_2} \!-\! z)}}\left( {\frac{{p_1^2 + p_2^2}}{{{p_1} + {p_2}}} - z - %\frac{{{p_1}{p_2}}}{z} + \left( {{{(z - \frac{{{p_1}{p_2}}}{z})}^2}} \right.} \right.} \right. \\ \\
%{\kern 18pt}  \displaystyle  {\left. {\left. {{{\left. { + {{({p_1} \!-\! {p_2})}^2}(1 - 2\frac{{z + {p_1}{p_2}/z}}{{{p_1} + {p_2}}})} \right)}^{\frac{1}{2}}}} \right)} %\right\}^{\frac{1}{2}}}.\label{PI-phase-analytical}
%\end{array}
%$$
%\end{equation}
The optimal PI coefficients are given by $k^\ast_i=0$, and
\begin{equation}\label{optimal-kp}
k^\ast_p = \sqrt{\frac{p_1+p_2}{z}(\omega^2_0+p_1p_2)}.
\end{equation}
%\begin{equation}\label{optimal-kp}
%k^\ast_p = \sqrt{\frac{p_1+p_2}{z}(\omega^2_0+p_1p_2)}.
%\end{equation}

%(ii).
%\begin{equation} \label{2-44+1}
%\begin{array}{l}
%{\kern -10pt} \displaystyle {\tan ^{ - 1}}\frac{{{p_i} - z}}{{2\sqrt {{p_i}z} }} - {\tan ^{ - 1}}\frac{{\sqrt {{p_i}z} }}{{{p_j}}} \le \theta _M^{P}\\
%{\kern 106pt} \displaystyle  \le {\tan ^{ - 1}}\frac{{{p_1}{p_2} - \left( {{p_1} + {p_2}} \right)z}}{{2\sqrt {{p_1}{p_2}z\left( {{p_i} - z} \right)} }},
%\end{array}
%\end{equation}
%where $i,j = 1,2$ and $i \ne j$.

\textit{Proof}: See Appendix \ref{Theorem 5.1-Proof}. $\hfill\blacksquare$

%\textit{Remark 5.1}:  In fact, it is not difficult to show that $\theta_2\left( {{\omega _2}} \right)$ is convex %function of $\omega_2$ in the interval $\omega^2_2 \in %[0,~p_1p_2 - z(p_1+p_2)]$.
%in the neighbourhood of $\omega_0$.
%To see this, we take the second-order derivative
%\begin{equation} \nonumber
%\frac{{{d^2}{\theta_2}}}{{d\omega^2}} = 2{\omega}\left( {\frac{z}{{{{\left( {\omega^2 + {z^2}} \right)}^2}}} - \frac{{{p_1}}}{{{{\left( {\omega^2 + p_1^2} \right)}^2}}} - %\frac{{{p_2}}}{{{{\left( {\omega^2 + p_2^2} \right)}^2}}}} \right).
%\end{equation}
%Noting from (\ref{derivative2}), however,
%\begin{equation} \nonumber
%\frac{{{z^2}}}{{{{\left( {\omega^2_0 + {z^2}} \right)}^2}}} > \frac{{p_1^2}}{{{{\left( {\omega^2_0 + p_1^2} \right)}^2}}} + \frac{{p_2^2}}{{{{\left( {\omega^2_0 + p_2^2} %\right)}^2}}},
%\end{equation}
%or equivalently,
%\begin{eqnarray*}
%\frac{{{z}}}{{{{\left( {\omega^2_0 + {z^2}} \right)}^2}}} & > & \left(\frac{p_1}{z}\right)\frac{{p_1}}{{{{\left( {\omega^2_0 + p_1^2} \right)}^2}}}
%+ \left(\frac{p_2}{z}\right)\frac{{p_2}}{{{{\left( {\omega^2 + p_2^2} \right)}^2}}}, \\
%\hfill & > & \frac{{p_1}}{{{{\left( {\omega^2_0 + p_1^2} \right)}^2}}}+\frac{{p_2}}{{{{\left( {\omega^2 + p_2^2} \right)}^2}}}.
%\end{eqnarray*}
%Here the second inequality follows from the fact ${p_1} > z,~{p_2} > z$. This leads to the conclusion that
%${d^2}{\theta }/d\omega_0^2 > 0$, and hence  $\theta_2\left( {{\omega _2}} \right)$ is convex in the neighbourhood of $\omega_0$.

%$\theta_2$ is a convex function of $\omega_2$.

%\textit{Remark 5.2}:
Again, similar to the case of first-order plants, it is clear from Theorem 4.2, Theorem 4.3 and Theorem 5.1 that with both PID and PI controllers, the integral control has no effect in improving the gain and phase margins; in fact, from the proofs of both theorems, one can see that a non-zero integral control gain $k_i$ will render the gain and phase margins smaller. On the other hand, these theorems also reveal that the derivative control improves the margins. The following result establishes this fact for the gain margin, and meanwhile it shows that the improvement is at best within a factor of two.
 %Consider, for example, the case $p_1+p_2\leq p_1p_2/z$ in (\ref{2-45-1}). In this case, it follows from Theorem 4.2 and (\ref{2-44}) that
%\begin{equation} \label{PID-PD-Gain}
%k_\M^{\PD} = \frac{{{p_1} + {p_2}}}{{{p_1} + {p_2} - z}}k_\M^\p.
%\end{equation}

%In the sequel, we end this section with a comparison of the maximal margins herein with those attainable by the PID controllers given in Theorem 4.2. The following corollary presents the comparison for both second-order unstable plants described by  (\ref{2-49}).

{\textit{Corollary 5.1}}: Let $P(s)$ be given by (\ref{2-49}). Under the condition ${p_1} + {p_2} < {p_1}{p_2}/z$,

\begin{equation} \label{corollary 5.1-gain}
k_\M^{\PI} < k_\M^{\PID} < 2k_\M^{\PI}.
%k_\M^{\PI} \le k_\M^{\PID} \le 2k_\M^{\PI}.
%k_\M^{\PI} \le k_\M^{\PID} \le 2k_\M^{\PI}. \\
\end{equation}
%In particular, the lower bound holds strictly if $z > 0$, $z \ne {p_1}$, $z \ne {p_2}$.

%(ii). For ${k_p} \ge {p_1} + {p_2} + z$,
%\begin{equation} \label{corollary 5.1-phase}
%\theta _\M^{\PI} < \theta _\M^{\PID} < 2\theta _\M^{\PI}.
%\theta _\M^{\PI} \le \theta _\M^{\PID} \le 2\theta _\M^{\PI}.
%\end{equation}

\textit{Proof}: Clearly, the lower bound of (\ref{corollary 5.1-gain}) holds strictly for $z > 0$. An inspection of (\ref{2-45-1}), (\ref{gain-complexpoles2}) and (\ref{2-44}) shows that under the condition
${p_1} + {p_2} < {p_1}{p_2}/z$,
%For the plant (\ref{2-49}) that can be both stabilized by PID controller ${K_{\PID}}\left( s \right)$ and PI controller ${K_{PI}}\left( s \right)$, clearly, it is necessary that ${p_1} + {p_2} \le {p_1}{p_2}/z$. This gives rise to (\ref{PID-PD-Gain}) thence suggesting that $k_\M^{\PID} \ge  k_\M^{\PI}.$
%$$
%k_\M^{\PID} = \frac{{{p_1} + {p_2}}}{{{p_1} + {p_2} - z}}k_\M^{\PI}  \ge  k_\M^{\PI}.
%$$
%Noting also from Theorems 4.2 and 5.2, however, that necessarily ${p_1} \ge z$, ${p_2} \ge z$.
% in this vein.
%As such, we are then led to
%It follows that under the condition ${p_1} + {p_2} < {p_1}{p_2}/z$, it is necessary that ${z < {p_1}{p_2}/\left( {{p_1} + {p_2}} \right) < {p_1},{p_2}}$. Noting also that ${p_1}{p_2}/\left( {{p_1} + {p_2}} \right) < \left( {{p_1} + {p_2}} \right)/2$.
$$
k_\M^{\PID} = \left( {1 + \frac{z}{{{p_1} + {p_2} - z}}} \right)k_\M^{\PI} < 2k_\M^{\PI},
$$
where the last inequality of (\ref{corollary 5.1-gain}) follows from the fact that under the condition ${p_1} + {p_2} < {p_1}{p_2}/z$, it is necessary that $z < {p_1}{p_2}/({p_1} + {p_2})$, and additionally, ${p_1} > z$, ${p_2} > z$ for real
poles and $|p-z|>z$ for complex poles.
$\hfill\blacksquare$

\textit{Remark 5.1}: Similarly, it is plausible to conjecture that
\begin{equation} \label{corollary 5.1-phase}
\theta _\M^{\PI} < \theta _\M^{\PID} < 2\theta _\M^{\PI}.
%\theta _\M^{\PI} \le \theta _\M^{\PID} \le 2\theta _\M^{\PI}.
\end{equation}
Here the first inequality is obvious, but the second inequality is nontrivial. We show below, however, that it holds under the
condition $p_1+p_2+z\leq p_1p_2/z$. For this purpose, we note from Theorem 4.2 that for any $k_p\in (p_1+p_2-z, ~p_1p_2/z)$,
$\hat{\omega}(k_p)\in (0,~\bar{\omega})$, with $\bar{\omega}$ given by (\ref{omegabar}), and that
%It follows from  (\ref{2-45+2})  and (\ref{2-44-phase}) that the lower bound of (\ref{corollary 5.1-phase}) holds. Noting from Theorem 4.2 that
$$
\begin{aligned}
 \theta_\M^{\PID} &= \tan^{ - 1}\frac{\hat{\omega} ({\hat k}_p)}{{\hat k}_p}
+ \tan^{ - 1}\frac{\hat{\omega} ({\hat k}_p)}{z} - \tan^{ - 1}\frac{({p_1} + {p_2})\hat{\omega} ({{\hat k}_p})}{{p_1}{p_2} - \hat{\omega}^2({{\hat k}_p})}
\\
& = {\tan ^{ - 1}}\frac{{\hat{\omega} ({{\hat k}_p})}}{{{{\hat k}_p}}} -\theta_2(\hat{\omega}(\hat{k}_p)) \\
&\leq {\tan ^{ - 1}}\frac{{\hat{\omega} ({{\hat k}_p})}}{{{{\hat k}_p}}} - \inf_{\omega\in (0,~\bar{\omega})}\theta_2(\omega).
\end{aligned}
$$
Here $\theta_2(\omega)$ is defined in the proof of Theorem 5.1, where it is also shown that
$ - \inf_{\omega\in (0,~\bar{\omega})}\theta_2(\omega)=\theta^-_{\PI}\leq \theta_\M^{\PI}$.
%However, from the proof of Theorem 5.1, we have $\sup\limits_{\omega\in (0,~\bar{\omega})}\bar{\theta}(\omega)=\bar{\theta}(\omega_0)$.
Therefore,
$$
\theta_\M^{\PID}\leq \theta_\M^{PI}+\tan^{ - 1}\frac{\hat{\omega} (\hat k_p)}{\hat k_p}.
$$
Under the condition $p_1+p_2+z\leq p_1p_2/z$, we have
%
%
%{\kern - 9 pt} \displaystyle  \le {\tan ^{ - 1}}\frac{{\omega ({{\bar k}_p})}}{{{{\bar k}_p}}} \!+\! \! \mathop {\sup }\limits_{{k_p} > {p_1} \!+ \!{p_2}} ({\tan ^{ - %1}}\frac{{\omega ({{\bar k}_p})}}{z} \!-\! {\tan ^{ - 1}}\frac{{({p_1} \!+\! {p_2})\omega ({{\bar k}_p})}}{{{p_1}{p_2} \!-\! {\omega ^2}({{\bar k}_p})}})
%\end{array}
%$$
%which leads to
%$$
%\begin{array}{l}
%\displaystyle  \theta _\M^{\PID} \le {\tan ^{ - 1}}\frac{{\omega ({{\bar k}_p})}}{{{{\bar k}_p}}} + {\tan ^{ - 1}}\frac{{{\omega _0}}}{z} - {\tan ^{ - 1}}\frac{{\left( {{p_1} %+ {p_2}} \right){\omega _0}}}{{{p_1}{p_2} - \omega _0^2}}\\
%{\kern 24 pt} \displaystyle   = {\tan ^{ - 1}}\frac{{\omega ({{\bar k}_p})}}{{{{\bar k}_p}}} + \theta _\M^{\PI}.
%\end{array}
%$$
%It follows analogously from the proof of Theorem 4.2 that
$$
{\tan ^{ - 1}}\frac{{\hat\omega ({{\hat k}_p})}}{z} - {\tan ^{ - 1}}\frac{{\hat\omega ({{\hat k}_p})}}{{{{\hat k}_p}}} \ge {\tan ^{ - 1}}\frac{{({p_1} + {p_2})\hat\omega ({{\hat k}_p})}}{{{p_1}{p_2} - {\hat\omega ^2}({{\hat k}_p})}}.
$$
%since ${\omega ^2}({{\bar k}_p}) < z{{\bar k}_p} < {p_1}{p_2}$, ${{\bar k}_p} - z \ge {p_1} + {p_2}$.
This gives rise to
$$
\begin{aligned}
 {\tan ^{ - 1}}\frac{{\hat{\omega} ({{\hat k}_p})}}{{{{\hat k}_p}}}&\leq
{\tan ^{ - 1}}\frac{{\hat\omega ({{\hat k}_p})}}{z} -
{\tan ^{ - 1}}\frac{{({p_1} + {p_2})\hat\omega ({{\hat k}_p})}}{{{p_1}{p_2} - {\hat\omega ^2}({{\hat k}_p})}} \\
&\leq \theta _\M^{\PI},
\end{aligned}
$$
and henceforth (\ref{corollary 5.1-phase}).
$\hfill\blacksquare$

\section{Conclusion}
In this paper we have studied the gain and phase margins of LTI systems attainable using PID controllers. The issue is to seek
%the fundamental margins of gain and phase, or equivalently, the
intrinsic limits under which a PID controller may exist to stabilize the family of plants with different gain and phase values. We derived analytical expressions for the maximal gain and phase margins for first- and second-order plants. For minimum-phase systems, we proved that the maximal gain and phase margins  achievable by PID controllers coincide with those by LTI controllers. For nonminimum-phase systems,  the results show that the maximal gain margin achievable by general LTI controllers is at most twice that by PID controllers. %, when measured on the logarithmic scale.
 Furthermore, in all cases, the results reveal, in a way consistent with one's intuition, that the maximal gain and phase margins achievable by PID and PD controllers coincide. This in turn implies that integral control plays no positive role in enlarging gain and phase margins.

In practical implementation of PID controllers, the derivative control is often accompanied by a lowpass filter. Whether this implementation scheme
will restrict the gain and phase margins achievable is of significant practical interest but is yet to be investigated. In particular, it is of
interest to see whether a possible tradeoff exists between the margins desired and the bandwidth of the filter. Furthermore, while the present paper
is fully focused on the robustness margins, it will be useful to extend the study to incorporate performance specifications, including, e.g., specification
on asymptotic reference tracking. In this latter scenario, the integral control becomes necessary, thus ushering in a conflict between tracking
performance and the system's stability robustness.

\appendices
%%%%%%%%%%%%%%%%%%%%%%%%%%%Proof of Theorem 3.2
%%%%%%%%%%%%%%%%%%%%%%%%%%%Proof of Theorem 3.2 \displaystyle
\section{Mathematical Preliminaries}
\setcounter{equation}{0}
\setcounter{subsection}{0}
\renewcommand{\theequation}{A.\arabic{equation}}
\renewcommand{\thesubsection}{A.\arabic{subsection}}
The gain and phase margin optimization problems are reformulated in this paper as nonlinear programming problems, and the well-known KKT
condition is employed to tackle the problems. The KKT %Karush-Kuhn-Tucker
condition is  a first-order necessary condition of optimality for constrained nonlinear programming problems whose objective function and constraints are differentiable, and the constraints may be equalities or inequalities. This class of problems is stated as
\begin{equation}
\begin{aligned}
&\min \quad f(x)\\
\mathrm{s.t.} \quad g_i(x)\leq 0,\quad &i\in \mathcal{I}_1= \left\lbrace 1,2,\cdots,m \right\rbrace, \\
h_j(x)=0,\quad &j\in \mathcal{I}_2=\left\lbrace m+1,m+2,\cdots,k \right\rbrace,
\end{aligned}
\label{eq:programming problem}
\end{equation}
where $f\!:\!\mathbb{R}^n\rightarrow \mathbb{R},\,g_{i}\!:\!\mathbb{R}^n\rightarrow \mathbb{R}$ and  $h_{j}\!:\!\mathbb{R}^n\rightarrow \mathbb{R}$ all have continuous first-order partial derivatives in $\mathbb{R}^n$.

%\begin{prop}\label{lemma:Fritz John condition}
%{\textit{The Fritz John Condition}}:
\textit{Lemma A.1 (Karush–Kuhn–Tucker Condition)}: If $x^{*}$ is an optimal solution of $f(x)$ in \eqref{eq:programming problem}, then there exists a row vector $\lambda=[\lambda _1, \lambda _2,\dots,\lambda _k]$ such that:
\begin{align}
\nabla  f(x^{\ast})+\sum_{i\in \mathcal{I}_1} \lambda_i \nabla  g_i(x^{\ast})+\sum_{j\in \mathcal{I}_2} \lambda_j  \nabla h_j(x^{\ast})=0&,  \notag
%\label{eq:KKT first}
\\
g_i(x^{\ast}) \le 0,  i\in \mathcal{I}_1,\; h_j(x^{\ast}) = 0,  j\in \mathcal{I}_2&,  \notag
%\label{eq:KKT second}
%\\
%h_j(x^{\ast}) = 0,  j\in \mathcal{I}_2&,
%\label{eq:KKT third}
\\
\sum_{i\in \mathcal{I}_1} \lambda_i g_i(x^{\ast})=0&,  \notag
%\label{eq:KKT forth}
\\
\lambda_i\geq 0, i\in \mathcal{I}_1&, \notag
%\label{eq:KKT fifth}
%\\
%\lambda\neq 0&,
%\label{eq:KKT sixth}
\end{align}
where $\nabla \theta(x)=[\partial \theta(x)/\partial x_1,\,\cdots,\,\partial \theta(x)/\partial x_n]^{\top}$ denotes the gradient of $\theta(x)$.
%\end{prop}

%Consider the polynomial with complex coefficients
%\begin{equation}\nonumber

%P\left( z \right) = {z^n} + {q_1}{z^{n - 1}} + {q_2}{z^{n - 2}} + ... + {q_n},
%\end{equation}
%where ${q_k} = {a_k} + j{b_k},k = 1,2,...,n$ is the complex coefficient. Then the polynomial $P\left( z \right) $ has all its zeros in the half-plane ${\mathop{\rm Re}\nolimits} \left( z \right) > 0$ if and only if the determinants ${\Delta _1} = {a_1}$ and ${\Delta %_k}\left( {k = 2,3,...,n} \right)$ are all positive, where the explicit expression of ${\Delta _k} $ is defined in \cite{frank1946zeros}.

%%%%%%%%%%%%%%%%%%3.First-Order Unstable Systems
%%%%%%%%%%%%%%%%%%

We also use repeatedly several classical results from the theory of algebraic geometry. In dealing with
the phase maximization problem, we employ the Bilherz stability criterion,
which concerns stability of polynomials with complex coefficients. We quote this result from  \cite{frank1946zeros}, \cite{marden1949geometry}.

\textit{Lemma A.2 (Bilherz Criterion)} \cite{frank1946zeros,marden1949geometry}:
Consider the complex polynomial
\begin{equation}\label{Bilherz}
(a_0+jb_0)s^n+(a_1+jb_1)s^{n-1}+\cdots+(a_n+jb_n),~~a_0+jb_0\neq 0.
\end{equation}
Define the associated $(2i-1)\times (2i-1)$ Bilherz submatrices
$$
\Delta_i=\left[\begin{array}{cccccccc}
a_{1} & a_3 & \cdots & a_{2i-1} & -b_2 & -b_4 & \cdots & -b_{2i-2} \\
a_{0} & a_2 & \cdots & a_{2i-2} & -b_1 & -b_3 & \cdots & -b_{2i-3} \\
\cdots & \cdots & \cdots & \cdots & \cdots & \cdots & \cdots & \cdots \\
0 & 0 & \cdots & a_i & 0 & 0 & \cdots & -b_{i-1} \\
0 & b_2 & \cdots & b_{2i-2} & a_{1} & a_3 & \cdots & a_{2i-3} \\
0 & b_1 & \cdots & b_{2i-3} & a_0 & a_2 & \cdots & a_{2i-4} \\
\cdots & \cdots & \cdots & \cdots & \cdots & \cdots & \cdots & \cdots \\
0 & 0 & \cdots & b_i & 0 & 0 & \cdots & a_{i-1}\end{array}\right].
$$
Then, the polynomial (\ref{Bilherz}) is stable if and only if $\det (\Delta_i)$ are positive for all $i=1,~2,~\cdots,~n$.

Next, the {\em Descartes' Rule of Signs} for polynomials
%and the {\em Sylvester resultant} for polynomial pairs
is routinely used in our development.

%\begin{lemma}\label{lemma:Descartes rule}
%\textit{The Descartes Rule}
\textit{Lemma A.3 (Descartes' Rule of Signs)} \cite{borwein1995polynomials,anderson1998descartes}:
Consider the $n$th-order polynomial
$$
{\displaystyle P(x)=a_{n}x^{n}+a_{n-1}x^{n-1}+\cdots +a_{1}x+a_{0}},~~~~a_n\neq 0.
$$
The number of positive roots (multiplicities counted) of the polynomial $P(x)$ is either equal to
the number of sign differences between the consecutive nonzero coefficients of $P(x)$, or is less than it by an even number.
%\end{lemma}
%
%%\begin{lemma}\label{lemma:Sylvester}
%\textit{Lemma 2.4 (The Sylvester resultant)} \cite{apostol1970resultants}:
%Consider the $m$th-order polynomial
%$$
%F\left( y \right) = {f_0} + {f_1}y + {f_2}{y^2} + ... + {f_m}{y^m}, ~~~ {f_m} \ne 0,
%$$
%and the $n$th-order polynomial
%$$
%Q\left( y \right) = {q_0} + {q_1}y + {q_2}{y^2} + ... + {q_n}{y^n}, ~~~~ {q_n} \ne 0.
%$$
%The The Sylvester resultant associated to $F\left( y \right) $ and $Q\left( y \right)$ is constructed as
%\begin{equation} \nonumber
%{\rm{Syl}}\left( {F,Q,y} \right) = \left[ {\begin{array}{*{20}{c}}
%{{f_m}}&{{f_{m - 1}}}& \cdots &{{f_0}}&{}&{}\\
%{}& \ddots &{}&{}& \ddots &{}\\
%{}&{}&{{f_m}}&{{f_{m - 1}}}& \cdots &{{f_0}}\\
%{{q_n}}&{{q_{n - 1}}}& \cdots &{{q_0}}&{}&{}\\
%{}& \ddots &{}&{}& \ddots &{}\\
%{}&{}&{{q_n}}&{{q_{n - 1}}}& \cdots &{{q_0}}
%\end{array}} \right]
%\end{equation}
%Then the two polynomials $F\left( y \right) $ and $Q\left( y \right)$ have a root $y$ in common if and only if ${\rm{det}}\left( {{\rm{Syl}}\left( {F,Q,y} \right)} \right) = 0$.
%\textit{The Sylvester resultant}
%\end{lemma}

Finally, we use repeatedly the following properties of arctangent function \cite{boros2005sums}:
%\begin{lemma}\label{lemma:arctan}

\textit{Lemma A.4} Suppose that $\xi\geq 0$ and $\eta\geq 0$. Then,
\begin{itemize}

\item[(i).]
$\tan^{-1}\xi$ is monotonically increasing with $\xi$.

%\item[(ii)]
%$$
%\frac{\xi}{1+\xi^2}\leq \tan^{-1}\xi)\leq \min \left\lbrace\xi,\,\frac{\pi}{2} \right\rbrace .
%$$
\item[(ii).]
$$
\tan^{-1}\xi + \tan^{-1}\eta =\left\lbrace \begin{aligned}
&\tan^{-1}\frac{\xi+\eta}{1-\xi\eta},~~~~~~~~~~\xi\eta\leq 1,
\\
&\pi-\tan^{-1}\frac{\xi+\eta}{\xi\eta-1},~~~~~\xi\eta>1.
\end{aligned} \right.
$$
\item[(iii).]
$$
\tan^{-1}\xi - \tan^{-1}\eta=\tan^{-1}\frac{\xi-\eta}{1+\xi\eta}.
$$
%\item[(vi)]
\end{itemize}
%\end{lemma}
Here we restrict $\tan^{-1}(\cdot)$ to its principal values in $\left[-\pi/2,~\pi/2 \right] $.

%%%%%%%%%%%%%%%%%
%%%%%%%%%%%%%%%%%
\section{Proof of Theorem 3.1} \label{Theorem 3.1-Proof}
\setcounter{equation}{0}
\setcounter{subsection}{0}
\renewcommand{\theequation}{B.\arabic{equation}}
\renewcommand{\thesubsection}{B.\arabic{subsection}}
For ease of presentation, we divide the proof in three different cases. In each case, we first determine the maximal gain margin %for ${K_{\PID}}\left( s \right)$ to stabilize ${\mathscr{P}_\alpha }$
and next the maximal phase margin. % for ${K_{\PID}}\left( s \right)$ stabilizing ${\mathcal{Q}_\theta }$.

\textbf{ Case 1}: ${\beta_0} > 0$ and ${\beta_1} < 0$. To stabilize ${\mathscr{P}_\mu}$, the first step is to determine the feasible ranges of the PI coefficients $\left( {{k_p},{k_i}} \right)$.
With the PI controller, the closed-loop characteristic equation $1 +\alpha P\left( s \right){K_{\PI}\left( s \right)} = 0$ of the system is given by
\begin{equation} \label{1-10}
\left( {1 + \alpha {\beta_0}{k_p}} \right){s^2} + \left( {\alpha {\beta_1}{k_p} + \alpha {\beta_0}{k_i} - p} \right)s + \alpha {\beta_1}{k_i} = 0.
\end{equation}
%that also satisfies Assumption 2.1.
%It involves the closed-loop characteristic equation of ${\mathscr{P}_\alpha }$, i.e., $1 +\alpha P\left( s \right){K_{PI}\left( s \right)} = 0$. Consequently, it gives
Define the sets
\begin{equation} \nonumber
{\kern - 2pt} \Omega _\alpha ^ +  = \left\{ {\left( {{k_p},{k_i}} \right):\alpha {\beta _0}{k_p} >  \!-\! 1,\alpha \left( {{\beta _1}{k_p} \!+\! {\beta _0}{k_i}} \right) > p, {\beta _1}{k_i} > 0} \right\},
\end{equation}
and
\begin{equation} \nonumber
{\kern - 2pt}  \Omega _\alpha ^ -  = \left\{ {\left( {{k_p},{k_i}} \right):\alpha {\beta _0}{k_p} <  \!- \!1,\alpha \left( {{\beta _1}{k_p} \!+\! {\beta _0}{k_i}} \right) < p,{\beta _1}{k_i} < 0} \right\}.
\end{equation}
%Note that $L\left( \infty  \right) = {b_0}\left| {{k_p}} \right| < 1$ is necessary to fulfill Assumption 2.1.
%Invoking Routh-Hurwitz criterion gives that  $K_{PI}\left( s \right)$ is capable of stabilizing ${\mathscr{P}_\alpha }$ if and only if all the coefficients have the same sign for this second-order polynomials. We first consider all the coefficients in (\ref{1-10}) have positive sign. Define the set
Note that $\Omega _\alpha^ +  \cap \Omega _\alpha^ -  = \emptyset $. For $\alpha  = 1$, it follows from  Routh-Hurwitz criterion that $K_{\PI}\left( s \right)$ stabilizes $P\left( s \right)$ if and only if  $\left( {{k_p},{k_i}} \right) \in {\Omega _1^ +} $ or $\left( {{k_p},{k_i}} \right) \in {\Omega _1^ -} $.
%Clearly,
%\begin{equation} \nonumber
%k_\M^{\PI}{\rm{ = }}\max \left\{ {k_{PI}^ + ,k_{PI}^ - } \right\},
%\end{equation}
%where
%\begin{equation} \label{1-11-1}
%\begin{array}{l}
%{\kern -5pt} k_{PI}^ +  = \sup \left\{ {\alpha  > 1:\alpha {\beta _0}{k_p} >  - 1,\alpha {\beta _1}{k_p} + \alpha {\beta _0}{k_i} > p,} \right.\\
%{\kern 148pt} {k_i} < 0,\left. {\left( {{k_p},{k_i}} \right) \in \Omega _1^ + } \right\},
%\end{array}
%\end{equation}
%\begin{equation} \label{1-11-2}
%\begin{array}{l}
%{\kern -5pt} k_{PI}^ -  = \sup \left\{ {\alpha  > 1:\alpha {\beta _0}{k_p} <  - 1,\alpha {\beta _1}{k_p} + \alpha {\beta _0}{k_i} < p,} \right.\\
%{\kern 148pt} {k_i} > 0,\left. {\left( {{k_p},{k_i}} \right) \in \Omega _1^ - } \right\}.
%\end{array}
%\end{equation}
Define
$$
k_{\PI}^ +  = \sup \left\{ {\alpha  \ge 1:\left( {{k_p},{k_i}} \right) \in \Omega _\alpha ^ + } \right\},
$$
and
$$
k_{\PI}^ -  = \sup \left\{ {\alpha  \ge 1:\left( {{k_p},{k_i}} \right) \in \Omega _\alpha ^ - } \right\}.
$$
Consider first that  $\left( {{k_p},{k_i}} \right) \in {\Omega _1^ +} $. It follows that $- 1/{\beta _0} < {k_p} < \left( {{\beta _0}{k_i} - p} \right)/\left| {{\beta_1}} \right|$, ${{k_i} < 0}$, which implies that $\left| {{\beta_1}} \right|/{\beta_0} > p$. Hence, we are led to
%Beyond that, it is worth noting from the set ${\Omega _1^ + }$ that $- 1/{\beta _0} < {k_p} < \left( {{\beta _0}{k_i} - p} \right)/\left| {{\beta_1}} \right|\left( {{k_i} < 0} \right)$, which gives rise to $\left| {{\beta_1}} \right|/{\beta_0} > p$.
%A necessary and sufficient condition when $\alpha  = 1$ leads to $\left( {{k_p},{k_i}} \right) \in {\Omega _1} $ for ${K_{PI}}\left( s \right)$ stabilizing ${\mathscr{P}_\alpha }$. Beyond that, it is worth noting from the set ${\Omega _1}$ that $- \frac{1}{{{b_0}}} < {k_p} < \frac{{{b_0}{k_i} - p}}{{\left| {{b_1}} \right|}}\left( {{k_i} < 0} \right)$, which gives rise to $\left| {{b_1}} \right|/{b_0} > p$. Hence, we are led to
%\begin{equation} \label{1-11}
%\begin{array}{l}
%{\kern -5pt} \displaystyle {k_\M^{\PI} \ge \sup \left\{ {\alpha  > 1:\alpha {\beta_0}{k_p} >  - 1,\alpha {\beta_1}{k_p} + \alpha {\beta_0}{k_i} > p,} \right.}\\
%{\kern 145pt}  \displaystyle { {k_i} < 0,\left. {\left( {{k_p},{k_i}} \right) \in \Omega _1^ + } \right\}}\\
%{\kern 0 pt}  \displaystyle {= \sup \left\{ {\alpha  > 1: - \frac{1}{{{\beta _0}}} < \alpha {k_p} < \frac{{\alpha {\beta _0}{k_i} \!-\! p}}{{\left| {{\beta _1}} \right|}},\left( {{k_p},{k_i}} \right) \in \Omega _1^ + } \right\}.}
%\end{array}
%\end{equation}
\begin{flalign} \label{1-11}
%\begin{array}{l}
& {\kern -6pt} \displaystyle k_{\PI}^ +  \le \sup \left\{ {\alpha  \ge 1: \!-\! \frac{1}{{{\beta _0}}} < \alpha {k_p} < \frac{{\alpha {\beta _0}{k_i} \!-\! p}}{{\left| {{\beta _1}} \right|}},\left( {{k_p},{k_i}} \right) \in \Omega _1^ + } \right\} \nonumber\\
& {\kern 10 pt}  = \left| {{\beta _1}} \right|/\left( {{\beta _0}p} \right).
%\end{array}
\end{flalign}
%It is clear that $\alpha $ achieves its maximum at $k_p^* =  - p/\left| {{\beta_1}} \right|,k_i^* = 0$.
It is clear that this upper bound can be attained asymptotically by the pair $\left( {k_p^*,k_i^*} \right)$ in the closure of $\Omega _\alpha ^ + $, with  $k_p^* =  - p/\left| {{\beta_1}} \right|,k_i^* = 0$ selected.  Hence, $k_{\PI}^ + = \left| {{\beta _1}} \right|/\left( {{\beta _0}p} \right)$. Analogously,  for $\left( {{k_p},{k_i}} \right) \in {\Omega _1^ -} $, we note that ${\beta _0}p > \left| {{\beta _1}} \right|$ and then find that $k_{\PI}^ - = \left( {{\beta _0}p} \right)/\left| {{\beta _1}} \right|$. It follows that
$$
k_\M^{\PI} = k_\M^\p = \max \left\{ {k_{\PI}^ + ,k_{\PI}^ - } \right\}.
$$
This establishes the expression of $k_\M^{\PI}$ and $k_\M^{\p}$ in (\ref{2-9}).
%\begin{equation} \label{1-11}
%\begin{array}{*{20}{l}}
%{\kern -10pt}{k_\M^{\PI} = \sup \left\{ {\alpha  > 1:\alpha {b_0}{k_p} >  - 1,\alpha {b_1}{k_p} + \alpha {b_0}{k_i}} \right.}\\
%{\kern 130pt}{\left. { > p,\left| {\left( {{k_p},{k_i}} \right) \in {\Omega _1}} \right.} \right\}}\\
%{\kern 7pt}{ = \sup \left\{ {\alpha >1 :\alpha \frac{{{b_0}\left| {{k_i}} \right| + p}}{{\left| {{b_1}} \right|}} < \frac{1}{{{b_0}}},\left| {\left( {{k_p},{k_i}} \right) \in {\Omega _1}} \right.} \right\}}\\
%{\kern  7pt}{ = \frac{{\left| {{b_1}} \right|/{b_0}}}{p}.}
%\end{array}
%\end{equation}
%
%It follows from (\ref{1-11}) that $- \frac{1}{{{b_0}}} < \alpha {k_p} < \frac{{\alpha {b_0}{k_i} - p}}{{\left| {{b_1}} \right|}}$. Clearly, $\alpha $ achieve its maximum at $k_p^* =  - \frac{p}{{\left| {{b_1}} \right|}},k_i^* = 0$ and thus we obtain $k_\M^{\PI} = \frac{{\left| {{b_1}} \right|/{b_0}}}{p}$. Analogously, we have $k_\M^\p = \sup \left\{ {\alpha  > 1:\alpha {b_0}{k_p} >  - 1, - \frac{1}{{{b_0}}} < {k_p} <  - \frac{p}{{\left| {{b_1}} \right|}}} \right\} = \frac{{\left| {{b_1}} \right|/{b_0}}}{p}$ for ${K_{P}}\left( s \right)$ to stabilize ${\mathscr{P}_\alpha }$. In the case that all the coefficients in (\ref{1-10}) have negative sign, we are then led to $k_\M^{\PI} = k_\M^\p = \frac{p}{{\left| {{b_1}} \right|/{b_0}}}$ in a similar fashion, where $\frac{p}{{\left| {{b_1}} \right|/{b_0}}} > 1$.

%We now focus on the maximal phase margin and examine $1 + {e^{ - j\theta }}P\left( s \right){K_{PI}}\left( s \right) = 0$, which results in

To derive the maximal phase margin $\theta _\M^{\PI}$, we examine  the characteristic equation $1 + {e^{ - j\theta }}P\left( s \right){K_{\PI}}\left( s \right) = 0$, which is given by
\begin{equation} \label{1-12}
\begin{array}{l}
{\kern -13pt}\left( {1 + {e^{ - j\theta }}{\beta_0}{k_p}} \right){s^2} + \left( {{e^{ - j\theta }}{\beta_1}{k_p} + {e^{ - j\theta }}{\beta_0}{k_i} - p} \right)s\\
{\kern 144pt} + {e^{ - j\theta }}{\beta_1}{k_i} = 0.
\end{array}
%\left( {1 \! +\! {e^{ -\!  j\theta }}{\beta _0}{k_p}} \right){s^2} \! + \!  \left( {{e^{ -\!  j\theta }}{\beta _1}{k_p} \! +\!  {e^{ - \! j\theta }}{\beta _0}{k_i} \! - \! p} \right)s \! + \! {e^{ -\!  j\theta }}{\beta _1}{k_i} = 0.
\end{equation}
Likewise, for $\theta = 0$, ${K_{\PI}}\left( s \right)$ stabilizes $P(s)$ if and only if  $\left( {{k_p},{k_i}} \right) \in {\Omega _1^ +} $ or $\left( {{k_p},{k_i}} \right) \in {\Omega _1^ -} $. Furthermore, in view of the Bilherz criterion (Lemma A.2), $K_{\PI}\left( s \right)$ stabilizes ${\mathcal{Q}_\nu}$, that is, the zeros of the complex polynomial in (\ref{1-12}) are in the open left-half of the complex plane if and only if the determinants ${\Delta _1} = {a_1} > 0$ and ${\Delta _2} = a_1^2{a_2} + {a_1}{b_1}{b_2} - b_2^2 >0$, which are found explicitly as
\begin{equation}\nonumber
{\kern -17pt} {\Delta _1} = \frac{{\left( {\left( { \!- {\beta _0}p + {\beta _1}} \right){k_p} + {\beta _0}{k_i}} \right)\cos \theta  + {\beta _0}{\beta _1}k_p^2 + \beta _0^2{k_p}{k_i} \!- \!p}}{{{{\left| {1 + {e^{ - j\theta }}{\beta _0}{k_p}} \right|}^2}}},
\end{equation}
\begin{equation}\nonumber
\begin{array}{l}
\displaystyle {\kern -5pt} {\Delta _2} = b_1^2\frac{{{\beta _1}{k_i}\cos \theta  + {\beta _0}{\beta _1}{k_p}{k_i}}}{{{{\left| {1 + {e^{ - j\theta }}{\beta _0}{k_p}} \right|}^2}}} - \frac{{\beta _1^2k_i^2{{\sin }^2}\theta }}{{{{\left| {1 + {e^{ - j\theta }}{\beta _0}{k_p}} \right|}^4}}} + \\
\displaystyle {\kern 98pt}  {b_1}\frac{{\left( {{\beta _0}p + {\beta _1}} \right){\beta _1}{k_p}{k_i} + {\beta _0}{\beta _1}k_i^2}}{{{{\left| {1 + {e^{ - j\theta }}{\beta _0}{k_p}} \right|}^4}}}{\sin ^2}\theta .
\end{array}
\end{equation}
%Note that the strict inequalities are revised to make the infimum and the minimum value equivalent, that is, ${\Delta _1} \ge 0$ and ${\Delta _2} \ge 0$.
Define similarly
\begin{equation}\nonumber
\begin{array}{l}
{\kern -5pt}  \theta _{\PI}^ +  = \max \left\{ {\sup \left\{ {\theta  \ge 0:\left( {{k_p},{k_i}} \right) \in \Omega _1 ^ + } \right\},} \right.\\
{\kern 115pt}  \left. { - \inf \left\{ {\theta  < 0:\left( {{k_p},{k_i}} \right) \in \Omega _1 ^ + } \right\}} \right\},
\end{array}
\end{equation}
and
\begin{equation}\nonumber
\begin{array}{l}
{\kern -5pt}  \theta _{\PI}^ -  = \max \left\{ {\sup \left\{ {\theta  \ge 0:\left( {{k_p},{k_i}} \right) \in \Omega _1 ^ - } \right\},} \right.\\
{\kern 115pt}  \left. { - \inf \left\{ {\theta  < 0:\left( {{k_p},{k_i}} \right) \in \Omega _1 ^ - } \right\}} \right\}.
\end{array}
\end{equation}
We first consider that $\left( {{k_p},{k_i}} \right) \in {\Omega _1^ +} $. %and hence,
%For ${K_{PI}}\left( s \right)$ to stabilize $P(s)$ when $\theta  = 0$, it yields $\left( {{k_p},{k_i}} \right) \in {\Omega_1} $ when all the coefficients in (\ref{1-12}) have positive sign. It follows from Lemma 2.1 that all the zeros of (\ref{1-12}) are in the open left-half of the complex plane for the second-order polynomial with complex coefficients to be stable i.e., $K_{PI}\left( s \right)$ stabilizes ${\mathcal{Q}_\theta }$, if and only if the determinants ${\Delta _1} = {b_1}$ and ${\Delta _2} = b_1^2{b_2} + {b_1}{c_1}{c_2} - c_2^2$ are all positive. Hence, we obtain
It is clear that
\begin{equation} \label{1-13}
%\begin{flalign} \label{1-13}
\begin{array}{*{20}{l}}
&{\kern -12pt}{\theta _{\PI}^ + \le \sup \left\{ {\theta  \ge 0:\cos \theta  > f\left( {{k_p},{k_i}} \right),{k_i} > \left( {p\! -\! ({\beta_1}/{\beta_0})} \right){k_p},\;} \right.}  \\  %\nonumber
&{\kern 172pt}{\left. {\left( {{k_p},{k_i}} \right) \in \Omega _1^ + } \right\},}
\end{array}
%\end{flalign}
\end{equation}
where the function $f\left( {{k_p},{k_i}} \right)$ is denoted by
\begin{equation} \nonumber
f\left( {{k_p},{k_i}} \right) = \frac{{ - {\beta_0}{\beta_1}k_p^2 - \beta_0^2{k_p}{k_i} + p}}{{\left( { - {\beta_0}p + {\beta_1}} \right){k_p} + {\beta_0}{k_i}}}.
\end{equation}
Here we obtain the upper bound in (\ref{1-13}) by weakening the Bilherz condition to $\Delta_1>0$ alone.
Since ${\cos \theta }$ is monotonically decreasing for $\theta  \in \left[ {0,\pi } \right]$, the maximal $\theta$ such that the constraint holds corresponds to the minimum of $f\left( {{k_p},{k_i}} \right)$.
Taking the first-order partial derivatives of $f\left( {{k_p},{k_i}} \right) $ with respect to ${k_p}$ and ${k_i}$, we have
%It is noted that the function ${\cos \theta }$ is a monotonically decreasing for $\theta  \in \left[ {0,\pi } \right]$. In other words, determining the maximal phase margin is equivalent to evaluate the minimal value of ${\cos \theta }$. Hence, we can recast the phase margin maximization in (\ref{1-13}) as a minimization problem of $f\left( {{k_p},{k_i}} \right) $.
%Towards this end, we examine the first-order partial derivatives of $f\left( {{k_p},{k_i}} \right) $ with respect to ${k_p}$ and ${k_i}$ built on the constrained conditions of the set ${\Omega _1}$, that is
%{\raggedright
%Taking the first-order partial derivatives of $f\left( {{k_p},{k_i}} \right) $ with respect to ${k_p}$ and ${k_i}$, we have
%}
\begin{equation} \label{1-14}
{\kern -8 pt}\frac{{\partial f}}{{\partial {k_p}}} \!=\! \frac{{{\beta_0}{\beta_1}\left( {{\beta_0}p \! -\! {\beta_1}} \right)k_p^2 \!-\! 2\beta_0^2{\beta_1}{k_p}{k_i} \!-\! \beta_0^3k_i^2 \!+\! p\left( {{\beta_0}p \!-\! {\beta_1}} \right)}}{{{{\left( {\left( { - {\beta_0}p + {\beta_1}} \right){k_p} + {\beta_0}{k_i}} \right)}^2}}},
\end{equation}
\begin{equation} \label{1-15}
{\kern -120pt}\frac{{\partial f}}{{\partial {k_i}}} = \frac{{{\beta_0}p\left( {\beta_0^2k_p^2 - 1} \right)}}{{{{\left( {\left( { - {\beta_0}p + {\beta_1}} \right){k_p} + {\beta_0}{k_i}} \right)}^2}}}.
\end{equation}
Note that for any $\left( {{k_p},{k_i}} \right) \in {\Omega _1^ +} $, ${\beta_0}\left| {{k_p}} \right| < 1$ and ${k_i} < 0$. Accordingly, $f\left( {{k_p},{k_i}} \right) $ is monotonically decreasing with ${k_i}$ and achieves its minimum at $k_i^* = 0$. Setting $k_i = 0$ and ${{\partial f}}/{{\partial {k_p}}} = 0$ in (\ref{1-14}), we obtain $k_p^* =  - \sqrt {p/\left( {{\beta_0}\left| {{\beta_1}} \right|} \right)} $. It is easy to verify that $f\left( {{k_p},0} \right)$ achieves the minimum at ${k_p} = k_p^*$. Beyond that, we assert that ${\Delta _2} \ge 0$ holds by submitting $k_p^*$ and $ k_i^*$ into ${\Delta _2}$. In view of (\ref{1-13}), the upper bound of $\theta _\M^{\PI}$ can be achieved and we are thus led to
\begin{equation}\nonumber
\cos \theta _\M^{\PI} = \cos \theta _\M^\p = f\left( { \!- \sqrt {{p}/{({{\beta_0}\left| {{\beta_1}} \right|}})} ,0} \right) = \frac{{2\sqrt {{\beta_0}\left| {{\beta_1}} \right|/p} }}{{{\beta_0} \!+\! (\left| {{\beta_1}} \right|/p})}.
\end{equation}
In a similar manner, for $\left( {{k_p},{k_i}} \right) \in {\Omega _1^ -} $, we also obtain that
$$\cos \theta _\M^{\PI} = \cos \theta _\M^\p = \frac{{2\sqrt {{\beta_0}\left| {{\beta_1}} \right|/p} }}{{{\beta_0} + (\left| {{\beta_1}} \right|/p})}.$$
Consequently, we prove the expression of $\theta _\M^{\PI} $ and $\theta _\M^{\p} $ in (\ref{2-9}).

\textbf{ Case 2}: ${\beta_0} > 0$ and ${\beta_1}  >  0$. %We may apply a similar analysis
%invoke similar techniques in the case
%to determine the maximal gain and phase margins.
In this case, it suffices to consider ${K_\p}\left( s \right)$ alone, which stabilizes $P\left( s \right)$ if $1 + {\beta _0}{k_p} > 0,{\beta _1}{k_p} > p$. It follows that
\begin{equation}\nonumber
\begin{array}{l}
{\kern -5pt}k_\M^\p \ge \sup \left\{ {\alpha  > 1:\alpha {\beta _0}{k_p} >  - 1,\alpha {\beta _1}{k_p} > p,} \right.\\
{\kern 150pt} \left. {{\beta _0}{k_p} >  - 1,{\beta _1}{k_p} > p} \right\}\\
{\kern 12pt}= \infty.
\end{array}
\end{equation}
Hence, $k_\M^{\PI} = k_\M^{\p} =\infty $. Similarly, we have
\begin{equation}\nonumber
\begin{array}{l}
{\kern -5 pt}  \displaystyle \theta _\M^\p \ge \sup \left\{ {\theta  > 0:\cos \theta _M^{P} \ge } \right.\frac{{ - {\beta _0}{\beta _1}k_p^2 + p}}{{\left( { - {\beta _0}p + {\beta _1}} \right){k_p}}},\\
{\kern 150 pt} \displaystyle  \left. {{\beta _0}{k_p} >  - 1,{\beta _1}{k_p} > p} \right\}\\
\displaystyle  {\kern 10 pt}  = \pi ,
\end{array}
\end{equation}
which gives rise to $\theta _\M^{\PI} = \theta _\M^\p = \pi $.
%\begin{equation}\nonumber
%k_\M^{\PI} = \sup \left\{ {\alpha  > 1:\alpha \left( {{b_1}{k_p} + {b_0}{k_i}} \right) > p,\left( {{k_p},{k_i}} \right) \in {\Omega _2}} \right\} = \infty ,
%\end{equation}
%\begin{equation}\nonumber
%{\kern -112pt}\cos \theta _\M^{\PI} > \frac{{ - {b_0}{b_1}k_p^2 - b_0^2{k_p}{k_i} + p}}{{\left( { - {b_0}p + {b_1}} \right){k_p} + {b_0}{k_i}}}.
%\end{equation}
%where ${\Omega _2} = \left\{ {\left( {{k_p},{k_i},{k_d}} \right):{b_1}{k_p} > p,{b_1}{k_i} > 0,{b_1}{k_d} >  - 1} \right\}.$ The right hand side of the last inequality arrive the minimum $-1$ at $k_p^* = \frac{1}{{{b_0}}},k_i^* = 0$. We are thus led to $k_\M^{\PI} = k_\M^\p = \infty$ and $\theta _\M^{\PI} = \theta _\M^\p = \pi $. Beyond that, it is clear that the case that all the coefficients of closed-loop characteristic equation have negative sign shares the same results. This establishes (\ref{2-9}).

\textbf{ Case 3}: ${\beta_0} = 0$. In this case, since $P( s )$ is strictly proper, a full PID controller can be used. The proof follows analogously as in Case 2. $\hfill\blacksquare$
%Likewise, we characterize the maximal gain and phase margins in this case. Without loss of generality, we assume that ${b_1}>0$. The proof as well as the results shares the analogous technique to those of the case ${b_0} > 0$ and hence are omitted. In particular, for P controllers to stabilize ${\mathcal{Q}_\theta }$, it gives rise to
%\begin{equation}\nonumber
%\theta _\M^\p = \sup \left\{ {\theta  > 0:{b_1}{k_p}\cos \theta  > p,{k_p} > p/{b_1}} \right\} = \frac{\pi }{2}.
%\end{equation}
%This completes the proof. $\hfill\blacksquare$

\section{Proof of Theorem 4.2} \label{Theorem 4.2-Proof}
\setcounter{equation}{0}
\setcounter{subsection}{0}
\renewcommand{\theequation}{C.\arabic{equation}}
\renewcommand{\thesubsection}{C.\arabic{subsection}}
%We first prove the theorem for the case of two real unstable poles, i.e., $p_1>0$, $p_2>0$.
%The proof for the case of complex conjugate poles can be established analogously and hence is omitted.
%
%For clarity, we present this proof in two cases. The first case is that $p_1$, $p_2$ are %two unstable
%real poles,  After that, we then consider that $p_1$, $p_2$ are %$a pair of conjugate complex poles,
%a complex conjugate pair, where  ${p_1} = \sigma  + j\nu$, ${p_2} = \sigma  - j\nu$, $\sigma > 0$.
%
%\textbf{ Case 1}: $p_1>0$, $p_2>0$. \noindent
At outset, it is worth pointing out that while Theorem 4.2 concerns only real poles, much of its proof
in this appendix, however, does not differentiate between real and complex poles and hence carries over
directly to the latter case, i.e., Theorem 4.3. We begin by defining, for $\alpha\geq 1$, the sets
%$\Xi _\alpha^ +$, $\Xi _\alpha^ -$
\begin{equation}\nonumber
\begin{array}{*{20}{l}}
{\kern -5pt} {\Xi _\alpha ^ +  = \left\{ {\left( {{k_p},{k_i},{k_d}} \right):} \right.\alpha {k_p} < \frac{{\alpha {k_i} + {p_1}{p_2}}}{z},\alpha {k_i} < 0, - 1 < \alpha {k_d}}\\
{\kern 40pt}  { < \frac{{\alpha {k_p} - {p_1} - {p_2}}}{z},\left( {\alpha {k_p} - z\alpha {k_d} - {p_1} - {p_2}} \right)\left( { - z\alpha {k_p}} \right.}\\
{\kern 94pt} {\left. {\left. { + \alpha {k_i} + {p_1}{p_2}} \right) >  - z\left( {1 + \alpha {k_d}} \right)\alpha {k_i}} \right\},}
\end{array}
\end{equation}
%\begin{equation}\nonumber
%\begin{array}{*{20}{l}}
%{\kern -5pt}{\Xi _1^ +  = \left\{ {\left( {{k_p},{k_i},{k_d}} \right):{k_p} < \frac{{{k_i} + {p_1}{p_2}}}{z},{k_i} < 0,} \right.}\\
%{\kern 145pt}{\left. { - 1 < {k_d} < \frac{{{k_p} - {p_1} - {p_2}}}{z}} \right\},}
%\end{array}
%\end{equation}
and
\begin{equation}\nonumber
\begin{array}{*{20}{l}}
{\kern -5pt}  {\Xi _\alpha ^ -  = \left\{ {\left( {{k_p},{k_i},{k_d}} \right):} \right.\alpha {k_p} > \frac{{\alpha {k_i} + {p_1}{p_2}}}{z},\frac{{\alpha {k_p} - {p_1} - {p_2}}}{z} < \alpha {k_d}}\\
{\kern 40pt} { <  - 1,\alpha {k_i} > 0,\left( {\alpha {k_p} - z\alpha {k_d} - {p_1} - {p_2}} \right)\left( { - z\alpha {k_p}} \right.}\\
{\kern 94pt} {\left. {\left. { + \alpha {k_i} + {p_1}{p_2}} \right) >  - z\left( {1 + \alpha {k_d}} \right)\alpha {k_i}} \right\}.}
\end{array}
\end{equation}
By Routh-Hurwitz criterion, the PID controller ${K_{\PID}}(s)$
%\left( s \right) = {k_p} + ({k_i}/s) + {k_d}s$
can stabilize ${\mathscr{P}_\mu}$ for all $\alpha\in [1,~\mu)$ if and only if
$\left( {{k_p},{k_i},{k_d}} \right) \in \Xi _\alpha ^ + $ or $\left( {{k_p},{k_i},{k_d}} \right) \in \Xi _\alpha ^ - $ for all $\alpha\in [1,~\mu)$. Moreover, for ${K_{\PID}}\left( s \right)$ to stabilize $P(s)$,
it is necessary that $\left( {{k_p},{k_i},{k_d}} \right) \in \Xi _1^ + $ or $\left( {{k_p},{k_i},{k_d}} \right) \in \Xi _1^ - $,
which requires, respectively,
\begin{eqnarray*}
p_1+p_2 &<& \frac{p_1p_2}{z}+z, \\
p_1+p_2 &>& \frac{p_1p_2}{z}+z.
\end{eqnarray*}
These conditions can be rewritten as
$({p_1} - z)({p_2} - z) > 0$ and $({p_1} - z)({p_2} - z) < 0$.
Since they are mutually exclusive, it is necessary to address the cases separately. Furthermore, the
condition $({p_1} - z)({p_2} - z) > 0$ implies that ${p_1} > z,{p_2} > z$, or ${p_1} < z,{p_2} < z$.
Consider the gain margin problem. It follows that
%
%
% In this vein, it follows from ${k_p}$ in ${\Xi _1^ +}$ that ${{p_1} + {p_2} < ({{{p_1}{p_2}}}/{z}) + z}$ and ${{p_1} + {p_2} > ({{{p_1}{p_2}}}/{z}) + z}$ for the case %$\left( {{k_p},{k_i},{k_d}} \right) \in \Xi _1^ - $ due to the fact that $k_i<0$.
%Using the Routh-Hurwitz criterion, a necessary and sufficient condition when ${K_{\PID}}\left( s \right)$ stabilizes $P(s)$ gives rise to $\left( {{k_p},{k_i},{k_d}} \right) %\in {\Xi _1^ +}$ or $\left( {{k_p},{k_i},{k_d}} \right) \in {\Xi _1^ -}$. It follows from ${k_p}$ in ${\Xi _1^ +}$ that ${{p_1} + {p_2} < ({{{p_1}{p_2}}}/{z}) + z}$ is %necessary and ${{p_1} + {p_2} > ({{{p_1}{p_2}}}/{z}) + z}$ for the case $\left( {{k_p},{k_i},{k_d}} \right) \in \Xi _1^ - $.
%Note that $\Xi _1^ +  \cap \Xi _1^ -  = \emptyset $. Clearly, we have
\begin{equation} \nonumber
k_\M^{\PID} = \left\{ \begin{array}{l}
k_{\PID}^ + , \;\;\;  {\rm{if}}  \; ({p_1} - z)({p_2} - z) > 0,\\
k_{\PID}^ - ,  \;\;\;  {\rm{if}}  \; ({p_1} - z)({p_2} - z) < 0,
\end{array} \right.
\end{equation}
where
\begin{eqnarray}
%\begin{array}{l}
k_{\PID}^ +  &=& \sup \left\{ {\alpha  \ge 1:\left( {{k_p},{k_i},{k_d}} \right) \in \Xi _\alpha ^ + } \right\}, \label{Theo4.2-1}\\
k_{\PID}^ -  &=& \sup \left\{ {\alpha  \ge 1:\left( {{k_p},{k_i},{k_d}} \right) \in \Xi _\alpha ^ - } \right\}. \label{Theo4.2-2}
%{\kern  7pt} \le \sup \left\{ {\alpha  > 1:\alpha \left( {{k_p} - z{k_d}} \right) < {p_1} + {p_2},} \right.
%\end{array}\\
%{\kern 128pt} {\left. {\left( {{k_p},{k_i},{k_d}} \right) \in \Xi _1^ - } \right\}.}
\end{eqnarray}
We first seek to determine $k_{\PID}^ +$. Toward this end,  we note that for any $\left( {{k_p},{k_i},{k_d}} \right) \in \Xi _\alpha ^ +$,
it is necessary that
\begin{equation}\label{ineqa1}
\begin{aligned}
& \alpha k_d > -1, \\
&\alpha k_i<0, \\
& \alpha k_d z+p_1+p_2 < \alpha k_p < \frac{\alpha k_i+p_1p_2}{z}<\frac{p_1p_2}{z}.
\end{aligned}
\end{equation}
Under the condition ${p_1} > z,{p_2} > z$, %${{p_1} + {p_2} \le  p_1p_2/z}$, then necessarily $p_1\geq z$, and $p_2\geq z$.
it follows from (\ref{ineqa1}) that
\begin{equation}\label{interval-kp}
0<p_1+p_2-z<\alpha k_p<\frac{p_1p_2}{z}.
\end{equation}
As such,
\begin{eqnarray*}
k_{\PID}^ + &\leq & \sup \left\{\alpha  \ge 1:~ p_1+p_2-z<\alpha k_p<\frac{p_1p_2}{z}\right\} \\
&=& \frac{p_1p_2}{z\left(p_1+ p_2 - z\right)},
\end{eqnarray*}
where the supremum is attained at $k_p^*=p_1+p_2-z$.
%$$
%k_p^* = \frac{{z(z - {p_1} - {p_2})}}{{{p_1}{p_2}}}.
%$$
Note, however, that the upper bound can be achieved asymptotically by the triple
$({k_p^*,k_i^*,k_d^*})$ in the closure of $\Xi _\alpha ^ + $, with this $k^*_p$, and additionally,
$k_d^* = -z(p_1+p_2-z)/(p_1p_2)$, $k_i^* = 0$, thus establishing that
$$
k_{\PID}^ +   = \frac{p_1p_2}{z\left(p_1+ p_2 - z\right)}.
$$
On the other hand, if ${p_1}< z,{p_2} < z$, clearly, the inequalities in (\ref{ineqa1}) imply that
\begin{equation}\label{interval-kd}
-1<\alpha k_d <\frac{p_1p_2-z(p_1+p_2)}{z^2}<0.
\end{equation}
Hence,
\begin{eqnarray*}
k_{\PID}^ + &\leq & \sup \left\{\alpha  \ge 1:~ -1<\alpha k_d <\frac{p_1p_2-z(p_1+p_2)}{z^2}\right\} \\
&=& \frac{{{z^2}}}{{z\left( {{p_1} + {p_2}} \right) - {p_1}{p_2}}}.
\end{eqnarray*}
Here the supremum is attained at
$$
k_d^* = \frac{{\left( {{p_1}{p_2} - } \right.z\left. {\left( {{p_1} + {p_2}} \right)} \right)}}{{{z^2}}}.
$$
%$k^*_d =\left(p_1p_2-z(p_1+p_2)\right)/z^2$.
Similarly, the upper bound can be achieved asymptotically by $k^*_d$ so selected, and additionally,  $\alpha k_p^* =   {{{p_1}{p_2}}}/{z}$, $k_i^* = 0$. The triple $( {k_p^*,k_i^*,k_d^*})$, likewise, lies in the closure of $\Xi _\alpha ^ + $. Hence, under the condition ${p_1}< z,{p_2} < z$,
%${{p_1} + {p_2} > p_1p_2/z}$,
we find that
$$
k_{\PID}^ +   = \frac{{{z^2}}}{{z\left( {{p_1} + {p_2}} \right) - {p_1}{p_2}}}.
$$
We have thus established (\ref{2-45-1}). For the case $({p_1} - z)({p_2} - z) < 0$, the derivation of $k_{\PID}^-$ is analogous, which results in (\ref{2-45-1a}).

We next establish the expressions for the phase margin, whose calculation, however, is considerably more
involved and the proof is longwinding.
As such, for it to be more accessible, we divide the proof in several steps, each of which begins
with a subheading highlighting the key technical background.

%%%%% Step 1:
\noindent  \textbf{ Step 1: Finding the critical phase with a set of fixed PID coefficients ${\left( {{k_p},{k_i},{k_d}} \right)}$}.
We denote similarly
\begin{equation} \nonumber
\theta _\M^{\PID} = \left\{ {\begin{array}{*{20}{l}}
{\theta^+_{\PID},\;\;\; {\rm{if}}\; ({p_1} - z)({p_2} - z) > 0,}\\
{\theta^-_{\PID},\;\;\;  {\rm{if}}\; ({p_1} - z)({p_2} - z) < 0,}
\end{array}} \right.
\end{equation}
where $\theta^+_{\PID}$ and $\theta^-_{\PID}$ represent the maximal phase margin in the two cases, respectively.
We first attempt to find $\theta_{\PID}^ +$.
For this purpose, we begin by examining the open-loop frequency response
\begin{equation}\nonumber
L\left( {j\omega } \right) = \frac{{j\omega-z }}{{\left( { j\omega-p_1 } \right)\left( {j\omega-p_2 } \right)}}\left( {{k_p} + \frac{{{k_i}}}{{j\omega }} + j{k_d}\omega } \right),
\end{equation}
whose magnitude is given by
\begin{equation}\label{T4.2-kpid0}
{\left| {L\left( {j\omega } \right)} \right|^2} = \frac{\left( {{\omega ^2} + {z^2}} \right)}{{\left( {{\omega ^2} + p_1^2} \right)\left( {{\omega ^2} + p_2^2} \right)}}
\left(k_p^2 +\left(k_d\omega-\frac{k_i}{\omega}\right)^2\right).
\end{equation}
Denote
$$
G(\omega)=\frac{\left( {{\omega ^2} + {z^2}} \right)}{{\left( {{\omega ^2} + p_1^2} \right)\left( {{\omega ^2} + p_2^2} \right)}}.
$$
Then the crossover frequencies ${\omega} > 0$ such that ${\left| {L\left( {j{\omega}} \right)} \right|^2} = 1$ can be
solved from the condition
\begin{equation}\label{cross1}
\left(k_p^2 +\left(k_d\omega-\frac{k_i}{\omega}\right)^2\right)G(\omega)=1,
\end{equation}
or equivalently, from the polynomial equation
\begin{equation}\label{T4.2-kpid1}
\begin{array}{l}
{\kern -18 pt} \left( {k_d^2 - 1} \right){\omega ^6} + \left( {k_p^2 - 2{k_i}{k_d} + {z^2}k_d^2 - p_1^2 - p_2^2} \right){\omega ^4} + \\
{\kern 12 pt}   \left( {{z^2}\left( {k_p^2 - 2{k_i}{k_d}} \right) + k_i^2 - p_1^2p_2^2} \right){\omega ^2} + {z^2}k_i^2 = 0.
\end{array}
\end{equation}
At each crossover frequency $\omega$, we have
\begin{equation}\nonumber
\sphericalangle L\left( {j{\omega}} \right) = \pi  + \sum\limits_{i = 1}^2 {{{\tan }^{ - 1}}\frac{{{\omega}}}{{{p_i}}}}  - {\tan ^{ - 1}}\frac{{{\omega}}}{z} + {\tan ^{ - 1}}\frac{{{k_d}{\omega} - \frac{{{k_i}}}{{{\omega}}}}}{{{k_p}}}.
\end{equation}
%Since $L\left( {j{\omega_l}} \right) = {e^{j\angle L\left( {j{\omega_l}} \right)}}$, we can match the phase of $L\left( {j{\omega_l}} \right)$, by setting
Define
\begin{equation}\label{T4.2-kpid2}
%{\theta _l} = {\tan ^{ - 1}}\frac{{{\omega _l}}}{{{p_1}}} + {\tan ^{ - 1}}\frac{{{\omega _l}}}{{{p_2}}} - {\tan ^{ - 1}}\frac{{{\omega _l}}}{z} + {\tan ^{ - 1}}\frac{{{k_d}{\omega _l} - \frac{{{k_i}}}{{{\omega _l}}}}}{{{k_p}}}.
\phi = \sum\limits_{i = 1}^2 {{{\tan }^{ - 1}}\frac{{{\omega}}}{{{p_i}}}}  - {\tan ^{ - 1}}\frac{{{\omega}}}{z} + {\tan ^{ - 1}}\frac{{{k_d}{\omega} - \frac{{{k_i}}}{{{\omega}}}}}{{{k_p}}}.
\end{equation}
%for some ${\theta _l}$.
Evidently,
$1 + e^{ - j\theta}L\left( {j\omega } \right) = 0$ whenever $\theta=\phi$. On the other hand, $1 + e^{ - j\theta}L\left( {j\omega } \right) \neq 0$ for all $\theta <\phi$ if $\phi\geq 0$, and for all $\theta>-\phi$ if $\phi<0$. Hence, with a given $\left( {{k_p},{k_i},{k_d}} \right)\in \Xi _1^ +$,
which means that $K_{\PID}(s)$ stabilizes $P\left( s \right)$,
the phase margin is $\phi$ if $\phi\geq 0$ and $-\phi$ if $\phi<0$. Rewrite $\phi=\phi \left( {{k_p},{k_i},{k_d}} \right)$. It follows that
the maximal phase margin achievable by $K_{\PID}(s)$ is given by the maximum of
\begin{equation}\nonumber
\begin{array}{l}
{\kern -4 pt} \phi^ +  = \sup \left\{ \phi\left( {{k_p},{k_i},{k_d}} \right) \ge 0: \right. ~\omega \mbox{ is a solution to } (\ref{cross1}), \\
{\kern 168 pt} \left. {\left( {{k_p},{k_i},{k_d}} \right) \in \Xi _1^ + } \right\},
\end{array}
\end{equation}
and
\begin{equation}\nonumber
\begin{array}{l}
{\kern -2 pt} \phi^ -  =  - \inf \left\{ {{\phi}\left( {{k_p},{k_i},{k_d}} \right) < 0:} \right.~\omega \mbox{ is a solution to } (\ref{cross1}),  \\
{\kern 170 pt}  \left. {\left( {{k_p},{k_i},{k_d}} \right) \in \Xi _1^ + } \right\},
\end{array}
\end{equation}
computed over all crossover frequencies that meet the equation (\ref{cross1}), or equivalently (\ref{T4.2-kpid1}); that is,
\begin{equation}\label{maxmin}
\theta _\M^{\PID} = \max \left\{ {{\phi ^ + },\;{\phi ^ - }} \right\}.
\end{equation}
Consequently, the computation of the maximal phase margin amounts to solving these two optimization problems.

%where ${\theta _l}$ is given by (\ref{T4.2-kpid2}). Clearly, we have
%\begin{equation}\nonumber
%\theta _{\Xi _1^ + }^{\PID}{\rm{ = }}\mathop {{\rm{min}}}\limits_l \left\{ {\theta _{\Xi _1^ + }^ + ,\theta _{\Xi _1^ + }^ - } \right\}.
%\end{equation}

%%%%% Step 2:
\noindent  \textbf{ Step 2:  Phase maximization reformulated as a nonlinear programming problem.}
One critical step in the subsequent proof is to recast the phase margin maximization problem, i.e., the computation of $\phi^+$ and $\phi^-$,
as a constrained nonlinear programming problem. Toward this end, it is clear that $\phi^-$ can be determined by solving the minimization problem
%
%We examine the maximization problem
%It follows from Theorems 3.1 and 3.2 that we can recast the phase margin maximization problem as a constrained nonlinear programming problem. We examine the maximization problem
\begin{equation}\nonumber
\begin{array}{l}
{\kern 70 pt}  \min \;\;\;f\left( {{k_p},{k_i},{k_d},{\omega}} \right)\\
{\rm{s}}.{\rm{t}}. {\kern 6 pt}  {\rm{     }}{g_1}\left( {{k_p},{k_i},{k_d},{\omega }} \right) =  - 1 - {k_d} \le 0,\\
{\kern 20 pt} {g_2}\left( {{k_p},{k_i},{k_d},{\omega}} \right) =  - {k_p} + z{k_d} + {p_1} + {p_2} \le 0,\\
{\kern 20 pt}{g_3}\left( {{k_p},{k_i},{k_d},{\omega}} \right) = z{k_p} - {k_i} - {p_1}{p_2} \le 0,\\
{\kern 20 pt}{g_4}\left( {{k_p},{k_i},{k_d},{\omega}} \right) = z{k_i} \le 0,\\
{\kern 20 pt}{g_5}\left( {{k_p},{k_i},{k_d},{\omega}} \right) =  - \left( {{k_p} - z{k_d} - {p_1} - {p_2}} \right)\left( { - z{k_p} + } \right.\\
\left. {{\kern 105 pt} {k_i} + {p_1}{p_2}} \right) - z\left( {1 + {k_d}} \right){k_i} \le 0,\\
{\kern 20 pt} h\left( {{k_p},{k_i},{k_d},{\omega }} \right) = {\left( {k_p^2 + \left( {{k_d}\omega  - \frac{{{k_i}}}{\omega }} \right)^2} \right)}G\left( {{\omega}} \right) - 1 = 0,
\end{array}
\end{equation}
where the objective function is given by $f\left( {{k_p},{k_i},{k_d},{\omega }} \right)=\phi\left( {{k_p},{k_i},{k_d}} \right)$, that is,
%
%$G\left( {{\omega _l}} \right) = \left( {\omega _l^2 + {z^2}} \right)/\left( {\left( {\omega _l^2 + p_1^2} \right)\left( {\omega _l^2 + p_2^2} \right)} \right)$, and
\begin{equation}\nonumber
\begin{array}{l}
{\kern -5pt} f\left( {{k_p},{k_i},{k_d},{\omega}} \right) =  \sum\limits_{i = 1}^2 {{{\tan }^{ - 1}}\frac{{{\omega}}}{{{p_i}}}}-{\tan ^{ - 1}}\frac{{{\omega}}}{z} + {\tan ^{ - 1}}\frac{{{k_d}{\omega} - \frac{{{k_i}}}{{{\omega}}}}}{{{k_p}}}.
\end{array}
\end{equation}
It is worth noting that the first five inequality constraints represent the set ${\Xi _1^ + }$, and the strict inequalities are used to make the infimum and the minimum value equivalent. The last equality constraint characterizes the crossover frequency $\omega> 0$.
%Noting, however, that in this formulation, the independent variable $\omega _l $ is limited by the equality constraint. As a consequence, the objective function and the %constraints are all continuously differentiable.
To tackle this nonlinear programming problem, we invoke the KKT condition and examine the first-order condition
\begin{equation}\nonumber
\begin{array}{*{20}{l}}
\begin{array}{l}
{\kern -5 pt}  \nabla f\left( {k_p^*,k_i^*,k_d^*,\omega^*} \right) + \sum\limits_{i = 1}^5 {{\lambda _i}\nabla {g_i}\left( {k_p^*,k_i^*,k_d^*,\omega ^*} \right)} \\
{\kern  110 pt} + {\lambda _6}\nabla h\left( {k_p^*,k_i^*,k_d^*,\omega ^*} \right) = 0,
\end{array}\\
{\kern -0 pt}  {\sum\limits_{i = 1}^5 {{\lambda _i}{g_i}\left( {k_p^*,k_i^*,k_d^*,\omega^*} \right)}  = 0,}\\
{\kern -0 pt}  {{\lambda _i} \ge 0,i \in \left\{ {1,2,3,4,5} \right\}.}
\end{array}
\end{equation}
The key idea is to show that to meet the first-order necessary condition of optimality, the optimal PID coefficients must be such that
$k_i^* = 0$. %and $k_d^* = -1$
To proceed, we 
%
%
%In particular, $k_i^* = 0$, then $\theta _\M^{\PID} = \theta _\M^{\PD}$. Under this circumstance, we first
take the partial derivatives of
%In an analogous manner, we invoke Karush–Kuhn–Tucker condition to solve the nonlinear programming above and examine the partial derivatives of
$f\left( {{k_p},{k_i},{k_d},{\omega}} \right)$, $h\left( {{k_p},{k_i},{k_d},{\omega}} \right)$, and ${g_i}\left( {{k_p},{k_i},{k_d},{\omega}} \right),~i = 1,2,3,4,5$ with respect to ${k_p}$, ${k_i}$, ${k_d}$, which gives rise to the following equations:
\begin{equation}\label{T4.2-kpid3}
\begin{array}{l}
\displaystyle {\kern -9pt} \frac{{\partial }}{{\partial {k_p}}}:\frac{1}{{1 + {{\rm H}^2}}}\frac{{k_d^*\omega^* \!- \frac{{k_i^*}}{{\omega ^*}}}}{{k_p^{*2}}} \!- {\lambda _2} + {\lambda _3}z + {\lambda _5}\left( {2z}  \! \right.k_p^* \!- k_i^* \!- {z^2}k_d^* \\
{\kern 58pt} \left. {- z\left( {{p_1} + {p_2}} \right) - {p_1}{p_2}} \right) + {\lambda _6} \cdot 2k_p^*G\left( {\omega ^*} \right) = 0,{\rm{   }}
\end{array}
\end{equation}
\begin{equation}\label{T4.2-kpid4}
\begin{array}{l}
\displaystyle {\kern -8pt}  \frac{{\partial }}{{\partial {k_i}}}:\frac{1}{{1 + {{\rm H}^2}}}\frac{1}{{\omega ^*k_p^*}} - {\lambda _3} + {\lambda _4}z + {\lambda _5}\left( { - k_p^* + {p_1} + {p_2} - z} \right)\\
{\kern 108pt}  + {\lambda _6} \cdot 2\left( {\frac{{k_i^*}}{{\omega ^{*2}}} - k_d^*} \right)G\left( {\omega ^*} \right) = 0,{\rm{   }}
\end{array}
\end{equation}
\begin{equation}\label{T4.2-kpid5}
\begin{array}{l}
\displaystyle {\kern -8pt}  \frac{{\partial }}{{\partial {k_d}}}:\frac{1}{{1 + {{\rm H}^2}}}\frac{{ - \omega ^*}}{{k_p^*}} - {\lambda _1} + {\lambda _2}z + {\lambda _5}\left( { - {z^2}k_p^* + {p_1}{p_2}z} \right) + \\
{\kern 108pt} {\lambda _6} \cdot 2\left( { - k_i^* + k_d^*\omega^{*2}} \right)G\left( {\omega^*} \right) = 0.
\end{array}
\end{equation}
Here ${\rm H} = \frac{{k_d^*\omega ^* - \frac{{k_i^*}}{{\omega ^*}}}}{{k_p^*}}$.
%We first consider the case $ {p_1} > z, {p_2} > z$.
For these equations, we first claim that it suffices to consider $\lambda_5=0$. Indeed, if $\lambda_5>0$, then the constraint ${g_5}\left( {{k_p^*},{k_i^*},{k_d^*},{\omega^* }} \right) \leq 0$ is active, that is,
\begin{equation}\label{T4.2-kpid7+1}
%{\kern -12 pt} \left( {k_p^* \!-\! zk_d^* \!-\! {p_1}\! -\! {p_2}} \right)\left( { - zk_p^* + } \right.\left. {k_i^* + {p_1}{p_2}} \right) =  -\! z\left( {1 + k_d^*} \right)k_i^*.\\
(k_p^* - zk_d^* - {p_1} - {p_2})( - zk_p^* + k_i^* + {p_1}{p_2}) =  - z(1 + k_d^*)k_i^*.
\end{equation}
From the Routh table associated with the closed-loop characteristic equation $1 + L\left( {s} \right) = 0$, this means that the characteristic equation
has a pair of imaginary roots $\pm j\omega^*$, as the solution to the equation
\begin{equation}\label{imaginary2}
\left( {k_p^* - zk_d^* - {p_1} - {p_2}} \right)\left( {j\omega^{*}} \right)^2 - zk_i^* = 0.
\end{equation}
%In other words,
%\begin{equation}\label{imaginary1}
%\omega^{*2}=\frac{zk^*_i}{-k_p^* +zk_d^* +{p_1} +{p_2}}.
%\end{equation}
As a result, $1 + L\left( {j\omega^{*} } \right) = 0$, %the system is unstable even when $\theta=0$
thus suggesting that $\phi\left( {{k^*_p},{k^*_i},{k^*_d}} \right)=0$. Indeed, this can be seen by
solving the equations (\ref{T4.2-kpid7+1}) and (\ref{imaginary2}), which yields
\begin{equation}\label{imaginary3}
\omega^{*2} = \frac{z\left( {{p_1}{p_2} - zk_p^*} \right)}{{k_p^* - {p_1} - {p_2} + z}}.
\end{equation}
From (\ref{T4.2-kpid7+1}) and (\ref{imaginary3}), we find that
$$
k^*_d \omega^{*2}-k^*_i=\frac{\left( {k_p^* - {p_1} - {p_2}} \right)\left( {{p_1}{p_2} - zk_p^*} \right)}{{k_p^* - {p_1} - {p_2} + z}}.
$$
Thus, for $k_p^* \geq {p_1} + {p_2}$, $k^*_d \omega^{*2}-k^*_i\geq 0$. Note also that as a function of $k^*_p$,
$p_1+p_2-z<k^*_p<p_1p_2/z$, $\omega^{*2}$ decreases monotonically with $k^*_p$, and hence from (\ref{imaginary3}) for $k_p^* > {p_1} + {p_2}$,
$\omega^{*2}<p_1p_2$. Using Lemma A.4, this allows us to obtain, for $k_p^* > {p_1} + {p_2}$,
$$
\sum\limits_{i = 1}^2 {{{\tan }^{ - 1}}\frac{{{\omega^*}}}{{{p_i}}}}=\tan^{-1}\frac{(p_1+p_2)\omega^*}{p_1p_2-\omega^{*2}}.
$$
Furthermore, by means of (\ref{imaginary2}) and (\ref{imaginary3}), and by invoking Lemma A.4, it follows similarly that
$$
{\tan ^{ - 1}}\frac{{{k^*_d}{\omega^*} - \frac{{{k^*_i}}}{{{\omega^*}}}}}{{{k^*_p}}}-{\tan ^{ - 1}}\frac{{{\omega^*}}}{z}=-\tan^{-1}\frac{(p_1+p_2)\omega^*}{p_1p_2-\omega^{*2}}.
$$
Hence, for $k_p^* > {p_1} + {p_2}$, $\phi\left( {{k^*_p},{k^*_i},{k^*_d}} \right)=0$. This invalidates that  $\lambda_5 > 0$ is a meaningful solution. Thus, we consider $\lambda_5 = 0$. For $p_1+p_2-z<k^*_p\leq p_1+p_2$, the same assertion holds analogously. We omit the proof for brevity.

Similarly, it suffices to consider $\lambda_2=0$. For if otherwise, then it is necessary that ${g_2}\left(\! {{k^*_p},{k^*_i},{k^*_d},{\omega^*}} \! \right)=0$. Together with the constraints ${g_4}\left(\! {{k^*_p},{k^*_i},{k^*_d},{\omega^*}} \!\right)\leq 0$ and ${g_5}\left( \!{{k^*_p},{k^*_i},{k^*_d},{\omega^*}}\! \right)\leq 0$, this mandates that $k^*_i=0$. It follows from (\ref{T4.2-kpid1}) that this equation has a solution $\omega^*=0$. Indeed, if $\omega^* > 0$, then it can be shown analogously that $\phi\left( {{k^*_p},{k^*_i},{k^*_d}} \right)=0$. For  $k^*_i=0$, $\omega^*=0$, we claim that $\phi\left( {{k^*_p},{k^*_i},{k^*_d}} \right)=0$ as well. Towards this end, we write the closed-loop characteristic polynomial as
$$
\begin{array}{l}
{\kern  -5  pt}  d(s,{k_i}) = (1 + {k_d}){s^3} + ({k_p} - z{k_d} - {p_1} - {p_2}){s^2}\\
 {\kern 130 pt}  + ({p_1}{p_2} - z{k_p} + {k_i})s - z{k_i}.
\end{array}
$$
When ${k_i} \to k_i^* = 0$, $d(s,{k_i})$ has a root $s \to {s^*} = 0$, i.e., $d({s^*},k_i^*) = 0$. From the eigenvalue perturbation result in \cite{chen2017stability} (Theorem 4), it follows that in the neighborhood  of $s^*$, $d(s,{k_i})$ has a root that can be expanded in the eigenvalue series
$$
s = \frac{z}{{{p_1}{p_2} - z{k_p}}}{k_i} + o(\left| {{k_i}} \right|).
$$
This gives rise to
$$
\frac{{{k_i}}}{s} = \frac{{{p_1}{p_2} - z{k_p}}}{z} + o(1).
$$
Therefore, when ${k_i} \to k_i^* = 0$, we have ${k_i}/s \to ({p_1}{p_2} - z{k_p})/z$, and in turn, $L(s) \to  - 1$. Consequently, $\sphericalangle L(j\omega ) =  - \pi $ when $k^*_i=0$, $\omega^*=0$, and thus $\phi\left( {{k^*_p},{k^*_i},{k^*_d}} \right)=0$.
%Next, we claim that it suffices to consider $\lambda_2=\lambda_3=0$ as well. For if otherwise, then it is necessary that ${g_2}\left(\! {{k^*_p},{k^*_i},{k^*_d},{\omega^*}} \! \right)=0$ or ${g_3}\left(\! {{k^*_p},{k^*_i},{k^*_d},{\omega^*}} \!\right)=0$. Together with the constraints ${g_4}\left(\! {{k^*_p},{k^*_i},{k^*_d},{\omega^*}} \!\right)\leq 0$ and ${g_5}\left( \!{{k^*_p},{k^*_i},{k^*_d},{\omega^*}}\! \right)\leq 0$, this mandates that $k^*_i=0$, thus activating the constraint ${g_5}\left( \!{{k^*_p},{k^*_i},{k^*_d},{\omega^*}} \!\right)\leq 0$. From the above proof, we arrive at the same conclusion that $\phi\left( \!{{k^*_p},{k^*_i},{k^*_d}} \! \right)=0$.

We now substitute $\lambda_2=\lambda_5=0$ into (\ref{T4.2-kpid3}), (\ref{T4.2-kpid4}), and (\ref{T4.2-kpid5}). It follows from (\ref{T4.2-kpid4}) and (\ref{T4.2-kpid5}) that
\begin{equation}\label{T4.2-kpid6}
\begin{array}{l}
{\kern -11pt} \!- {\lambda _1} - {\lambda _3}{\omega ^{*2}} + {\lambda _4}z{\omega ^{*2}} = 0,
\end{array}
\end{equation}
and from (\ref{T4.2-kpid3}) and (\ref{T4.2-kpid5}) that
\begin{equation}\label{T4.2-kpid7}
\begin{array}{l}
{\kern -2 pt}  {\omega ^*} + {\lambda _1}k_p^* + {\lambda _3}z(k_d^*{\omega ^{*2}} - k_i^*) = 0.
\end{array}
\end{equation}
These two equations yield possible solutions
\begin{equation}\label{Cases}
\left\{ \begin{array}{l}
{\rm{(i)}}~~\,{\lambda _1} > 0,\;{\lambda _4} > 0,~k^*_d<0\\
{\rm{(ii)}}~~{\lambda _1} = 0,\;{\lambda _3} > 0,\;{\lambda _4} > 0,~k^*_d<0
\end{array} \right.
\end{equation}
and ${\lambda _1} = 0$, $\omega^*=0$, $k^*_i=0$. It follows similarly that $\phi\left( {{k^*_p},{k^*_i},{k^*_d}} \right)=0$ when $\omega^*=0$, $k^*_i=0$. Otherwise, it holds that $k^*_d=-1$, $k^*_i=0$ in Case (i), and $k_p^* = {p_1}{p_2}/z$, $k^*_i=0$ in Case (ii).

\noindent  \textbf{ Step 3:  Determining $\theta^+_{\PID}$ for the case of ${p_1} > z,~{p_2} > z$.}
Consider now the case ${p_1} > z,{p_2} > z$. This step itself will be divided into two parts, dealing with
Case (i) and Case (ii) of (\ref{Cases}), respectively.

\noindent  \textbf{ Step 3a:  Minimizing $\hat \theta \left( {{k_p}} \right)$.}
In Case (i) of (\ref{Cases}), the function to be optimized
becomes $\phi\left(k_p,0,-1\right)$, and the equation (\ref{T4.2-kpid1}) reduces to
\begin{equation}\label{crossoverkp}
\left( {k_p^2 + {z^2} - p_1^2 - p_2^2} \right){\omega^{2}} + \left( {{z^2}k_p^2 - p_1^2p_2^2} \right) = 0.
\end{equation}
Then a necessary condition for the PID controller
${K_{\PID}}(s)$ to stabilize the plant $P(s)$, i.e., $e^{-j\theta}P(s)$ with $\theta=0$, is
\begin{equation}\label{boundary1}
0 < {p_1} + {p_2} - z < {k_p} < {p_1}{p_2}/z,
\end{equation}
that is, the condition (\ref{interval-kp}) with $\alpha=0$. It thus suffices to restrict $k_p$ to the interval given in (\ref{boundary1}).
Within this range, the equation (\ref{T4.2-kpid1}) results in the
unique positive solution $\hat{\omega}(k_p)$, given by (\ref{hattheta2}).
In what follows we prove $\hat\theta(k_p) < 0$ for $0<p_1+p_2-z<k_p<p_1p_2/z$,
where $\hat\theta (k_p)=\phi\left(k_p,0,-1\right)$ is defined by (\ref{hattheta1}). This being the case,
to determine $\phi^-$ for ${p_1} > z,{p_2} > z$,
%$p_1+p_2 \le p_1p_2/z$,
 it is necessary to
solve the univariate %maximization or
minimization problem of $\hat\theta\left(k_p\right)$. %and $\tilde \theta \left( {{k_d}} \right)$.
To facilitate the solution, we study the properties of $\hat\theta\left(k_p\right)$ for $0<p_1+p_2-z<k_p<p_1p_2/z$. %and further $\tilde \theta \left( {{k_d}} \right)$.
%For convenience, we write explicitly
%$$
%{{\omega}}}{z} - {\tan ^{ - 1}}\frac{\omega}{k_p}
%$$
%and
%$$
%\omega(k_p)=\sqrt {\frac{{p_1^2p_2^2 - {z^2}k_p^2}}{k_p^2 + {z^2} - p_1^2 - p_2^2}}.
%$$
A straightforward calculation shows that
\begin{equation}\label{T4.2-kpid11}
\frac{d\hat{\omega}(k_p)}{dk_p} = \frac{{{k_p}\left( {{z^2}\left( {p_1^2 + p_2^2} \right) - {z^4} - p_1^2p_2^2} \right)}}{{{{\left( {k_p^2 + {z^2} - p_1^2 - p_2^2} \right)}^2} \hat{\omega}(k_p) }},
\end{equation}
where $\hat{\omega}(k_p)$ is given in (\ref{hattheta2}).
%We analyze $\hat \theta \left( {{k_p}} \right)$ for $k_p$ in two disjoint intervals.
%\medskip\noindent {\textbf{Negative $\hat\theta\left(k_p\right)$: $0 < p_1+p_2-z<k_p<p_1p_2/z$.}}
% It follows from Step 2 that in this case, it holds that $ {p_1} > z, {p_2} > z$.
%In its full generality, i
It is immediately clear from (\ref{T4.2-kpid11}) that $\hat{\omega}(k_p)$ is monotonically decreasing.
%We first consider that ${k_p} > 0$ and then evaluate the value of $\hat \theta \left( {{k_p}} \right)$
%at ${k_p} = {p_1} + {p_2} - z$. For
At ${k_p} = {p_1} \!+\! {p_2} \!-\!z$, we have
$$
\hat{\omega}^2(p_1+p_2-z)=\bar{\omega}^2=\frac{p_1p_2+z(p_1+p_2-z)}{2}.
$$
%$, where
%$\bar{\omega}$ is given by (\ref{omegabar}).
% \left( {{p_1}{p_2} + z\left( {{p_1} + {p_2} - z} \right)} \right)/2$,
%\begin{equation}\nonumber
%\tilde \omega _1^2 = \frac{1}{2}\left( {{p_1}{p_2} + z\left( {{p_1} + {p_2} - z} \right)} \right),
%\end{equation}
This leads to ${p_1}{p_2} -  \bar{\omega}^2 =\left( {{p_1}{p_2} - z\left( {{p_1} + {p_2} - z} \right)} \right)/2 > 0$.
In other words, for all $k_p$ such that $p_1+p_2-z<k_p <p_1p_2/z$, $\hat{\omega}^2/{p_1}{p_2} < 1$. By Lemma A.4, it follows that
\begin{equation}\nonumber
{\tan ^{ - 1}}\frac{{{{ \hat\omega }}}}{{{p_1}}} + {\tan ^{ - 1}}\frac{{{{ \hat\omega }}}}{{{p_2}}} = {\tan ^{ - 1}}\frac{{\left( {{p_1} + {p_2}} \right){{ \hat\omega }}}}{{{p_1}{p_2} - \hat\omega ^2}}.
\end{equation}
Similarly,
\begin{equation}\nonumber
{\kern 0 pt}  {\tan ^{ - 1}}\frac{{{{\hat\omega }}}}{z} + {\tan ^{ - 1}}\frac{{{{ \hat\omega }}}}{{{k_p}}} = \left\{ \begin{array}{l}
\displaystyle  {\tan ^{ - 1}}\frac{{\left( {{k_p} + z} \right){{ \hat\omega }}}}{{{k_p}z -  \hat\omega ^2}}, \;\;\;\;\;\;\;\;\,  \frac{{ \hat\omega ^2}}{{{k_p}z}} \le 1,\\
\displaystyle  \pi  - {\tan ^{ - 1}}\frac{{\left( {{k_p} + z} \right){{\hat\omega }}}}{{ \hat\omega ^2 - {k_p}z}}, \;\; \frac{{ \hat\omega ^2}}{{{k_p}z}} > 1.
\end{array} \right.
\end{equation}
Hence, if $\hat{\omega}^2 > {k_p}z$, we have
\begin{equation}\nonumber
\hat \theta \left( {{k_p}} \right) =  - \pi  + {\tan ^{ - 1}}\frac{{({p_1} + {p_2})\hat\omega }}{{{p_1}{p_2} - {\hat\omega ^2}}} + {\tan ^{ - 1}}\frac{{({k_p} + z)\hat\omega }}{{{\hat\omega ^2} - {k_p}z}} < 0.
\end{equation}
Otherwise, if $\hat{\omega}^2 \le {k_p}z$, then
\begin{equation}\nonumber
\hat \theta \left( {{k_p}} \right) = {\tan ^{ - 1}}\frac{{\left( {{p_1} + {p_2}} \right){{ \hat\omega }}}}{{{p_1}{p_2} -  \hat\omega ^2}} - {\tan ^{ - 1}}\frac{{\left( k_p+z \right){{ \hat\omega }}}}{{{k_p}z -  \hat\omega^2}}.
\end{equation}
Since $p_1+p_2-z<k_p<p_1p_2/z$, we have
$$
\hat \theta \left( {{k_p}} \right) < {\tan ^{ - 1}}\frac{{\left( {{p_1} + {p_2}} \right){{\hat\omega }}}}{{{p_1}{p_2} -\hat\omega ^2}} - {\tan ^{ - 1}}\frac{{\left( k_p+z \right){{\hat\omega }}}}{{{k_p}z - \hat\omega^2}}<0
$$
As a consequence, we have shown that under the condition $0<p_1+p_2-z<k_p <p_1p_2/z$, $\hat \theta \left( {{k_p}} \right)<0$.
%On the other hand,
%in the second scenario, i.e., under the condition $p_1+p_2-z\leq 0<k_p <p_1p_2/z$, since $p_1+p_2\leq z$, we also have $d\omega(k_p)/dk_p \leq 0$ and
%henceforth $\omega^2(k_p)\leq \omega^2(0)\leq p_1p_2/2$. Mimicking the steps above, this too leads to $\hat \theta \left( {{k_p}} \right)<0$. To summarize,
%whenever $0<k_p <p_1p_2/z$, we have $\hat \theta \left( {{k_p}} \right)<0$.
Note in particular that at ${k_p} = {p_1}{p_2}/z$, $\hat\omega(k_p) = 0$ and hence
$\hat \theta \left( {{p_1}{p_2}/z} \right) = 0$.

%%%%% Step 4:
%\noindent  \textbf{ Step 4:  Minimizing $\hat \theta \left( {{k_p}} \right)$.} %$\phi _{{k_p}}^ - $
%and $\phi _{{k_p}}^ + $.} %and $\phi _{{k_d}}^ + $ ($\phi _{{k_d}}^ - $). }
We now attempt to find the minimum of $\hat \theta \left( {{k_p}} \right)$  by exploiting its differentiability. We show that this can be accomplished by solving the polynomial equation
(\ref{poly1}).
%
%
%In the sequel we show that
%$\hat \theta \left( {{k_p}} \right)$ is a monotone or a pseudo-convex function, and we find accordingly the global
%minimum of $\hat \theta \left( {{k_p}} \right)$ in the range of $p_1+p_2-z<k_p<p_1p_2/z$.
To begin with, we take the derivative of $\hat \theta \left( {{k_p}} \right)$ and obtain
\begin{equation}\label{T4.2-kpid10}
\frac{{d\hat \theta \left( {{k_p}} \right)}}{{d{k_p}}} = \displaystyle \frac{{\frac{{ \hat\omega '}}{{{p_1}}}}}{{1 + \frac{{{{ \hat\omega }^2}}}{{p_1^2}}}} + \frac{{\frac{{\hat\omega '}}{{{p_2}}}}}{{1 + \frac{{{{ \hat\omega }^2}}}{{p_2^2}}}} - \frac{{\frac{{ \hat\omega '}}{z}}}{{1 + \frac{{{{ \hat\omega }^2}}}{{{z^2}}}}} - \frac{{\frac{{{k_p} \hat\omega ' -  \hat\omega }}{{k_p^2}}}}{{1 + \frac{{{{ \hat\omega }^2}}}{{k_p^2}}}},
\end{equation}
where $\hat\omega '$ is the derivative of $\hat\omega(k_p)$, given in (\ref{T4.2-kpid11}).
%It is worth pointing out that (\ref{T4.2-kpid10}) holds for both
%positive and negative $k_p$.
%Since $\omega(p_1p_2/z)=0$, it follows that
%
%
%
%Note that ${p_1} + {p_2} - z < {k_p} < {p_1}{p_2}/z$. This gives rise to ${z^4} + p_1^2p_2^2 > {z^2}\left( {p_1^2 + p_2^2} \right)$. Consequently, we obtain that ${\tilde %\omega '} < 0$ at ${k_p} = {p_1}{p_2}/z$. Additionally, noting that $\tilde \omega  = 0$ for ${k_p} = {p_1}{p_2}/z$, we may rewrite (\ref{T4.2-kpid10}) as
%
%
%
%\begin{equation}\nonumber
%\frac{{d\hat \theta \left( {{p_1p_2/z}} \right)}}{{d{k_p}}} = \left( {\frac{1}{{{p_1}}} + \frac{1}{{{p_2}}} - \frac{1}{z} - \frac{1}{{{k_p}}}} \right)\omega'(p_1p_2/z).
%\end{equation}
%It is easy to see that
%By invoking the inequality ${k_p} < {p_1}{p_2}/z$ and noting that $ - {p_1}{p_2}/{k_p} <  - z$, it is not a difficult exercise to show that
%\begin{equation}\nonumber
%\begin{array}{l}
%\displaystyle \frac{1}{{{p_1}}} + \frac{1}{{{p_2}}} - \frac{1}{z} - \frac{1}{{{k_p}}}<
%\frac{1}{{{p_1}{p_2}}}\left( {{p_1} + {p_2} - z - \frac{{{p_1}{p_2}}}{{{z}}}} \right)< 0.
%\end{array}
%\end{equation}
%Hence, $\hat \theta '\left( {{k_p}} \right) > 0$ at ${k_p} = {p_1}{p_2}/z$. Furthermore, noting
Note that at $k_i^* = 0$, $k_d^* =  - 1$,
\begin{equation}\label{T4.2-kpid13}
\left( {{{\hat\omega }^2} + p_1^2} \right)\left( {{{\hat\omega }^2} + p_2^2} \right) = \left( {{{\hat\omega }^2} + {z^2}} \right)\left( {{{\hat\omega }^2} + k_p^2} \right),
\end{equation}
and by means of (\ref{T4.2-kpid11}), we obtain the equation (\ref{T4.2-kpid12}) at the top of the next page.
\begin{figure*}[bht]
%\hrulefill\\
\begin{equation}\label{T4.2-kpid12}
\displaystyle{
\frac{{d\hat \theta \left( {{k_p}} \right)}}{{d k_p }}
= \frac{\frac{(p_1+p_2-z-k_p)\hat\omega^2 +p_1p_2(p_1+p_2)-k_p z(k_p+z)}{(k^2_p+z^2-p^2_1-p^2_2)^2}k_p\left(z^2(p^2_1+p^2_2)-z^4-p^2_1p^2_2\right)+\hat\omega^2(\hat\omega^2+z^2)}
{{\hat\omega \left( {{{\hat\omega }^2} + p_1^2} \right)\left( {{{\hat\omega }^2} + p_2^2} \right)}}
}
\end{equation}
\hrulefill
\end{figure*}
By substituting (\ref{hattheta2}) into (\ref{T4.2-kpid12}), we find that
\begin{eqnarray}\label{polyequation-derivative}
\frac{{d\hat \theta \left( {{k_p}} \right)}}{{d k_p }} &=& \hat{W}(k_p)\left[k_p\left(k_p+z-p_1-p_2\right)\hat\omega^2\right. \nonumber \\
\hfill & \hfill & \left. ~+\left(k^3_p z-k_p p_1p_2(p_1+p_2)+p^2_1p^2_2\right)\right], ~ \label{T4.2-kpid12a}
\end{eqnarray}
where
\begin{equation}\label{T4.2-kpid12b}
 \hat{W}(k_p)=\frac{p^2_1p^2_2+z^4-z^2(p^2_1+p^2_2)}{\hat\omega (\hat\omega^2+p^2_1)(\hat\omega^2+p^2_2)(k^2_p+z^2-p^2_1-p^2_2)^2}>0.
\end{equation}
It thus follows that $\hat\theta'(k_p)=0$ if and only if
\begin{equation}\label{alge1}
k_p\left(k_p+z-p_1-p_2\right)\hat\omega^2+k^3_p z-k_p p_1p_2(p_1+p_2)+p^2_1p^2_2=0.
\end{equation}
Substituting (\ref{hattheta2}) into (\ref{alge1}) gives rise to the polynomial equation
\begin{equation}\label{kp-fifthorder-poly}
\begin{array}{l}
{\kern -9 pt} zk_p^5 - {z^2}k_p^4 + z\left( {({p_1} + {p_2})(z - {p_1} - {p_2} - ({p_1}{p_2}/z)) + 2{p_1}{p_2}} \right)\cdot\\
{\kern 12 pt}  k_p^3+ 2p_1^2p_2^2k_p^2 + \left( {{p_1}{p_2}({p_1} + {p_2})(p_1^2 + p_2^2 - {p_1}{p_2}) + } \right.z{p_1}{p_2} \cdot \\
{\kern 32 pt} \left. {({p_1}{p_2} - z({p_1} + {p_2}))} \right){k_p} + p_1^2p_2^2\left( {{z^2} - p_1^2 - p_2^2} \right) = 0.
\end{array}
\end{equation}
%We next seek to establish the uniqueness of  such a ${{\bar k}_p}$ that $\hat \theta '({{\bar k}_p}) = 0$.
%prove that there exists  a unique ${{\bar k}_p}$ such that $\hat \theta '({{\bar k}_p}) = 0$.
According to the Descartes' rule (Lemma A.3), the polynomial equation (\ref{kp-fifthorder-poly}) admits at most three positive solutions.
We claim, however, that only one positive solution lies in the interval $\left(p_1+p_2-z,~p_1p_2/z\right)$. To see this, denote
%However, whether these solutions fall into the interval $\left(p_1+p_2-z,~p_1p_2/z\right)$ cannot be ascertained in general. %Under certain circumstance,
%In fact, we claim that it can be determined by analysing the derivative of $\hat \theta \left( {{k_p}} \right)$ and the equation (\ref{kp-fifthorder-poly}).
%We take the derivative of $\hat \theta \left( {{k_p}} \right)$ and by a straightforward, albeit tedious calculation, we are led to
%$$
%\begin{array}{l}
%{\kern -5 pt} \hat \theta '\left( {{k_p}} \right) = W({k_p})\left( {{k_p}\left( {{k_p} + z - {p_1} - {p_2}} \right){\omega ^2}} \right.\\
%{\kern 90pt} \left. {\; + \left( {k_p^3z - {k_p}{p_1}{p_2}({p_1} + {p_2}) + p_1^2p_2^2} \right)} \right),
%\end{array}
%$$
%where $W({k_p}) > 0 $ is given in (\ref{T4.2-kpid12b}). Define
$$
D\left( {{k_p}} \right) = k_p^3z - {k_p}{p_1}{p_2}({p_1} + {p_2}) + p_1^2p_2^2.
$$
It is clear that ${k_p}\left( {{k_p} + z - {p_1} - {p_2}} \right){\hat\omega ^2} = 0$ at ${k_p} = {p_1} + {p_2} - z$ and at ${k_p} = {p_1}{p_2}/z$; the latter follows since at ${k_p} = {p_1}{p_2}/z$, $\hat\omega  = 0$. In view of (\ref{polyequation-derivative}), this suggests that at both ${k_p} = {p_1} + {p_2} - z$ and ${k_p} = {p_1}{p_2}/z$,
$\hat \theta '\left( {{k_p}} \right)=\hat{W}(k_p)D\left( {{k_p}} \right)$.
At ${k_p} = {p_1} \!+\! {p_2} \!- \!z$, we have
\begin{equation} \nonumber
\begin{array}{l}
{\kern -5 pt}  D\left( {{k_p}} \right) = {k_p}\left( {k_p^2z - {p_1}{p_2}({p_1} + {p_2} - z)} \right) + {p_1}{p_2}\left( {{p_1}{p_2} - z{k_p}} \right)\\
{\kern 23pt}  = \left( {k_p^2 - {p_1}{p_2}} \right)\left( {{k_p}z - {p_1}{p_2}} \right).
\end{array}
\end{equation}
Note that $k_p^2 - {p_1}{p_2} > 0$ since  $p_1 > z, p_2 > z$, which results in $D\left( {{k_p}} \right) < 0$ at ${k_p} = {p_1} \!+\! {p_2} \!- \!z$. As such,
we assert from (\ref{polyequation-derivative}) that $\hat \theta '\left( {{k_p}} \right) < 0$ at ${k_p} = {p_1} + {p_2} - z$. At ${k_p} = {p_1}{p_2}/z$, we have $D\left( {{k_p}} \right) = zk_p^2\left( {{k_p} - ({p_1} + {p_2} - z)} \right) > 0$; that is, $\hat \theta '\left( {{k_p}} \right) > 0$ at ${k_p} = {p_1}{p_2}/z$. These facts  imply  that there must exist at least one $\hat{k}_p\in \left(p_1\!+\!p_2\!-\!z,~p_1p_2/z\right)$ to the equation (\ref{kp-fifthorder-poly}).
Furthermore, it can be verified that the polynomial in (\ref{kp-fifthorder-poly}) can be factorized as
% It follows from (\ref{kp-fifthorder-poly}) that the equation can be factorized as
$$
z({k_p} - {p_1})({k_p} - {p_2})M({k_p}),
$$
where $M({k_p}) = k_p^3 + {c_{2}}k_p^2 + {c_{1}}{k_p} + {c_{0}}$ is given by (\ref{poly1}). This can be seen by verifying directly that the equation (\ref{alge1}) has solutions at $k_p=p_1$ and $k_p=p_2$.
Evidently, ${c_2} > 0,{c_1} < 0,{c_0} < 0$. As such, by Descartes' rule, $M({k_p}) = 0$ admits one and only one positive solution ${{\hat k}_p}$ in $\left(p_1\!+\!p_2\!-\!z,~p_1p_2/z\right)$. In summary, by now we have proved that there exists  a unique $\hat{k}_p\in \left(p_1\!+\!p_2\!-\!z,~p_1p_2/z\right)$ such that $\hat \theta '({{\hat k}_p}) = 0$. Hence, the minimum of $\hat \theta(k_p)$ is $\hat \theta({{\hat k}_p})$.
%Noting also that $\hat \theta ({k_p}) < 0$, it follows that for ${p_1} + {p_2} - z < {k_p} < {p_1}{p_2}/z$,
%$$
%\inf \hat \theta ({k_p}) = \min \{ \hat \theta ({{\hat k}_p}),\hat \theta ({p_1} + {p_2} - z)\},
%$$
%where ${{\hat k}_p}$ is a positive solution to the equation $M({k_p}) = 0$.
%Noting also that ${k_p} > {p_1} + {p_2} - z > {p_1},{p_2}$, this leads to that solving the equation reduces to find the solutions for $M({k_p}) = 0$, i.e., the equation (\ref{poly1}). Furthermore, there exist at most one positive solution $\bar{k}_p\in \left(p_1\!+\!p_2\!-\!z,~p_1p_2/z\right)$ to (\ref{poly1}).
%In summary, by now we have proved that there exists  a unique $\hat{k}_p\in \left(p_1\!+\!p_2\!-\!z,~p_1p_2/z\right)$ such that $\hat \theta '({{\hat k}_p}) = 0$.
%Hence, the minimum of $\hat \theta(k_p)$ is $\hat \theta({{\hat k}_p})$.
%$\theta _{\PID}^ +  = \phi^ -  = |\hat \theta ({{\hat k}_p})|.$ This completes the proof for the first expression in (\ref{2-45+2}).

%%%%% Step 5:
\noindent  \textbf{ Step 3b: Maximizing $\tilde \theta ({k_d})$.}
In parallel, under the same condition ${p_1} > z$, ${p_2} > z$,
%$p_1+p_2\leq p_1p_2/z$,
it remains to analyze Case (ii) of (\ref{Cases}). In this case, the function to be optimized is $\phi\left(p_1p_2/z,0,k_d\right)$, over the range $-1<k_d<({p_1}{p_2}/{z^2}) - (({p_1} + {p_2})/z)$. The equation (\ref{T4.2-kpid1}) now reduces to
\begin{equation}\label{crossoverkp}
\left(k^2_d-1\right){\omega^{2}} + \left( {{z^2}k_d^2 +\frac{p^2_1p^2_2}{z^2}- p_1^2-p_2^2} \right) = 0,
\end{equation}
which results in the solution $\tilde{\omega}(k_d)$ given in (\ref{hattheta-kd2}), and the function $\tilde{\theta}(k_d)$ given in (\ref{hattheta-kd1}). Evidently,
$\tilde{\theta}(-1)=0$. Calculating similarly the derivative of $\tilde{\theta}(k_d)$, we are led to
\begin{eqnarray}\label{polyequation-derivative-kd}
\frac{{d\tilde \theta ({k_d})}}{{d{k_d}}} &=& \tilde{W}({k_d})\left[ {(({p_1} + {p_2}){k_d} - zk_d^3 + ({p_1}{p_2}/z)){\tilde\omega ^2}} \right. \nonumber \\
\hfill & \hfill & \left. ~{ + {p_1}{p_2}{k_d}\left( z{k_d} + {p_1} + {p_2} - ({p_1}{p_2}/z) \right)}\right], ~ \label{T4.2-kpid12a}
\end{eqnarray}
where
\begin{equation}\label{T4.2-kpid12b-kd}
\tilde{W}({k_d}) = \frac{{{{({p_1}{p_2}/z)}^2} + {z^2} - p_1^2 - p_2^2}}{{\tilde\omega ({\tilde\omega ^2} + p_1^2)({\tilde\omega ^2} + p_2^2){{(1 - k_d^2)}^2}}} > 0.
\end{equation}
Under the condition $z \le {p_1}{p_2}/({p_1} + {p_2})$, it is clear that $({p_1}{p_2}/z) - {p_1} - {p_2} \ge 0$. Note from (\ref{Cases}) that ${k_d} < 0$. It is easy to realize that
$\tilde{\theta}'(k_d)>0$ for $k_d\in (-1,~0)$;
in other words, $\tilde{\theta}(k_d)$ is monotonically increasing in $(-1,~0)$, which,
together with the fact that $\tilde{\theta}(-1)=0$, suggests that $\tilde{\theta}(k_d)>0$ over $(-1,~0)$ and achieves its maximum at $k_d^* = 0$. We claim, however, that this leads to a contradiction since by (\ref{T4.2-kpid7}), it
necessitates ${\lambda _1} > 0,~k_p^* < 0$, which contradicts to the solution $k_p^*=p_1p_2/z$ in Case (ii)
of (\ref{Cases}).
%
%submitting $k_d^* = 0$ into (\ref{T4.2-kpid7}),
%
%
%this means that ${\lambda _1} > 0,k_p^* < 0$. On the other hand, for $z > {p_1}{p_2}/({p_1} + {p_2})$, we find that the %maximal phase is achieved at $k_d^* = {p_1}{p_2}/{z^2} - ({p_1} + {p_2})/z$ since $\tilde{\theta}'(k_d)>0$. It follows, %however, from the expression of $d(s,{k_i})$ that either $d(s,{k_i})$ has a root at $s = 0$ (i.e., ${\omega ^*} = 0$), or %$1 + k_d^* = 0$. In each case, it amounts to a contradiction, or $\phi\left(p_1p_2/z,0,k_d\right) = 0$. Hence,
As such, for ${p_1} > z, {p_2} > z$, Case (ii) can be ruled out.
In summary, we now ascertain that for ${p_1} > z, {p_2} > z$,
%${p_1} + {p_2} \le {p_1}{p_2}/z$,
\begin{eqnarray*}
{\phi ^ - } &=& -\inf\left\{\phi\left(k_p,0,-1\right) < 0:~\hat{\omega} ({k_p})  \mbox{ is given by (\ref{hattheta2})}, \right. \\
\hfill & \hfill & ~~~~~~~~~~~~~~~~~~~~~~~~\left. p_1+p_2-z<k_p<p_1p_2/z\right\} \\
\hfill &=& -\inf\{\hat\theta(k_p) < 0:~p_1+p_2-z<k_p<p_1p_2/z\} \\
\hfill &=& -\hat \theta ({{\hat k}_p}),
%\hfill &=& -\min \{ \hat \theta ({{\hat k}_p}),~\hat \theta ({p_1} + {p_2} - z)\},
\end{eqnarray*}
and $\theta _{\PID}^ +  = \phi ^ - =|\hat \theta ({{\hat k}_p})|$.
This completes Step 3 and in turn the proof for the first expression in (\ref{phase-realpoles}).

\noindent  \textbf{ Step 4:  Determining $\theta^+_{\PID}$ for the case of ${p_1} < z, ~{p_2} < z$.}
The case of ${p_1} < z$, ${p_2} < z$ is dual to that of ${p_1} > z$, ${p_2} > z$, and hence can be analyzed
analogously.
%In this case, a necessary condition for the PID controller
%${K_{\PID}}(s)$ to stabilize the plant $P(s)$ is
%\begin{equation}\label{boundary2}
%- 1 < {k_d}\; < \frac{p_1p_2 - z(p_1 + p_2)}{z^2} < 0,
%\end{equation}
%namely, the condition (\ref{interval-kd}) with $\alpha=0$.
By mimicking Steps 3a-3b, we find that Case (i) of (\ref{Cases}) can be ruled out.
%
%
%To see this, we first consider ${p_1} + {p_2} - z < 0$. For negative $k_p$, i.e., ${k_p} \in ({p_1} + {p_2} - z,0]$, (or %${k_p} \in ( - {p_1}{p_2}/z,0]$), the maximal phase margin is achieved at $k_p^* = 0$ due to monotonicity. Likewise, for %${k_p} \in [0,{p_1}{p_2}/z)$, the maximal phase margin is also achieved at $k_p^* = 0$. From (\ref{T4.2-kpid7}), however, %this necessitates that ${\lambda _3} > 0$, which in turn implies that $k_p^* = {p_1}{p_2}/z$, and hence leads to %contradiction. On the other hand, for the case ${p_1} + {p_2} - z \ge 0$, we first establish the monotonicity. It follows %that
%\begin{equation} \nonumber
%\begin{array}{l}
%{\kern -5 pt}  D({k_p}) = {k_p}(zk_p^2 - {p_1}{p_2}({p_1} + {p_2} - z)) + {p_1}{p_2}({p_1}{p_2} - z{k_p})\\
%{\kern 22 pt}  \ge {p_1}{p_2}({p_1}{p_2} - z{k_p})(2z - {p_1} - {p_2})/z\\
%{\kern 22 pt}  > 0
%\end{array}
%\end{equation}
%since ${p_1} < z$, ${p_2} < z$. In view of the monotonicity, the maximum of $\phi\left(k_p,0,-1\right)$ is achieved at %$k_p^* = {p_1} \!+\! {p_2}\!-\! z$. From the expression of $d(s,{k_i})$, this suggests that $d(s,{k_i})$ has a root at %${\omega ^*} = 0$, or ${p_1} {p_2}\!- \!z{k_p^*}= 0$, which results in a contradiction, or $\phi\left(k_p,0,-1\right) = 0$.
In Case (ii) of (\ref{Cases}), $\tilde \theta \left( {{k_d}} \right)=\phi ({p_1}{p_2}/z,0,{k_d})$, the function to be optimized, can be shown similarly to be a positive function for $ - 1 < {k_d}\; < ({p_1}{p_2} - z({p_1} + {p_2}))/z^2 < 0$.  Thus,
%$\varphi _{{k_d}}^ + $.
%Similarly, we analyze the case ${p_1} + {p_2} > {p_1}{p_2}/z$. The proof  follows analogously to the steps above; for this reason, we shall only provide a sketch of proof. %We begin with showing that
\begin{equation}\nonumber
\begin{array}{l}
{\kern -6 pt}  \theta^{+}_{\PID}  = \sup \{ \phi ({p_1}{p_2}/z,0,{k_d}) \ge 0:~\tilde{\omega} ({k_d})~\mbox{is given by (\ref{hattheta-kd2}),}\\
{\kern 90 pt} - 1 < {k_d} < ({p_1}{p_2} - z({p_1} + {p_2}))/z^2~\}\\
{\kern 12 pt}  {\rm{ = }}\sup \{ \tilde \theta ({k_d}) \ge 0: - 1 < {k_d}\; < ({p_1}{p_2} - z({p_1} + {p_2}))/z^2~\}.
\end{array}
\end{equation}
%where ${\tilde \theta \left( {{k_d}} \right)}$ is given in (\ref{hattheta-kd1}) and $\tilde{\omega}(k_d)$ by %(\ref{hattheta-kd2}).
With $\tilde{\theta}'(k_d)$ calculated in (\ref{polyequation-derivative-kd}), it can be shown similarly that $\tilde{\theta}'(k_d)=0$ if and only if
\begin{equation}\label{kd-fifthorder-poly}
- {z^2}\left({k_d} + \frac{p_1}{z}\right)\left({k_d} + \frac{p_2}{z}\right)N({k_d}) = 0,
\end{equation}
where $N({k_d}) = k_d^3 + {d_{2}}k_d^2 + {d_{1}}{k_d} + {d_{0}}$ is given by (\ref{kd-poly1}), and $N(k_d)$ has a unique negative root $\tilde{k}_d$
in the interval $(-1,~({p_1}{p_2} - z({p_1} + {p_2}))/z^2)$. Consequently, in the case ${p_1} < z,~ {p_2} < z$, $\theta^{+}_{\PID}= \tilde \theta ({{\tilde k}_d})$.
This completes the proof for the second expression in (\ref{phase-realpoles}), and hence the proof for the entire case of $({p_1} - z)({p_2} - z) > 0$.

\noindent  \textbf{ Step 5: Determining $\theta^-_{\PID}$.}
The final step in our proof is to establish (\ref{phase-realpoles1a}), that is, to calculate
$\theta^-_{\PID}$ in the case of $({p_1} - z)({p_2} - z) < 0$.
This step follows in exactly the same manner and hence is omitted. The entire proof is now completed. $\hfill\blacksquare$

\section{Proof of Theorem 4.3} \label{Theorem 4.3-Proof}
\setcounter{equation}{0}
\setcounter{subsection}{0}
\renewcommand{\theequation}{D.\arabic{equation}}
\renewcommand{\thesubsection}{D.\arabic{subsection}}

In the case of complex poles $p_1=\sigma+j\nu$ and $p_2=\sigma-j\nu$,
the proof shares the essential spirit of that for Theorem 4.2; indeed, one can see that much of the proof
for Theorem 4.2 carries over with no differentiation between real or complex poles.
For this reason, we only provide a sketch of proof for Theorem 4.3. To begin with,
we recognize that the feasible PID parameter set ${\Xi _\alpha ^ +}$ and ${\Xi _\alpha ^ -}$,
defined in Appendix C, remain in the same forms. The critical phase is found as
\begin{equation}\label{T4.2-kpid-complex1}
\phi  = {\tan ^{ - 1}}\frac{{\omega  + \nu }}{\sigma } + {\tan ^{ - 1}}\frac{{\omega  - \nu }}{\sigma } - {\tan ^{ - 1}}\frac{\omega }{z} + {\tan ^{ - 1}}\frac{{{k_d}\omega  - \frac{{{k_i}}}{\omega }}}{{{k_p}}}.
\end{equation}
Note that for ${\omega ^2} \le {\sigma ^2} + {\nu ^2} = {p_1}{p_2}$,
$$
\begin{array}{l}
\displaystyle {\kern -5 pt} {\tan ^{ - 1}}\frac{{\omega  + \nu }}{\sigma } + {\tan ^{ - 1}}\frac{{\omega  - \nu }}{\sigma } = {\tan ^{ - 1}}\frac{{2\sigma \omega }}{{{\sigma ^2} + {\nu ^2} - {\omega ^2}}}\\
\displaystyle {\kern 110 pt}  = {\tan ^{ - 1}}\frac{{\left( {{p_1} + {p_2}} \right)\omega }}{{{p_1}{p_2} - {\omega ^2}}},
\end{array}
$$
while for  ${\omega ^2} > {\sigma ^2} + {\nu ^2} = {p_1}{p_2}$,
$$
{\tan ^{ - 1}}\frac{{\omega  + \nu }}{\sigma } + {\tan ^{ - 1}}\frac{{\omega  - \nu }}{\sigma } = \pi  - {\tan ^{ - 1}}\frac{{\left( {{p_1} + {p_2}} \right)\omega }}{{{p_1}{p_2} - {\omega ^2}}}.
$$
In other words, $\phi $ is in the same form of (\ref{T4.2-kpid2}).
%Using exactly the same steps leads to the same answer for $\theta _\M^{\PID}$.
The crossover frequencies are determined in the form as in (\ref{cross1}), with
$$
\begin{array}{l}
\displaystyle {\kern  0 pt} G(\omega ) = \frac{{\left( {{\omega ^2} + {z^2}} \right)}}{{\left( {{{\left( {\omega  - \nu } \right)}^2} + {\sigma ^2}} \right)\left( {{{\left( {\omega  + \nu } \right)}^2} + {\sigma ^2}} \right)}}\\
\displaystyle {\kern 22 pt} = \frac{{\left( {{\omega ^2} + {z^2}} \right)}}{{\left( {{\omega ^2} + p_1^2} \right)\left( {{\omega ^2} + p_2^2} \right)}}.
\end{array}
$$
As such, the ensuing proof follows in exactly the same manner.  $\hfill\blacksquare$

\section{Proof of Theorem 5.1} \label{Theorem 5.1-Proof}
\setcounter{equation}{0}
\setcounter{subsection}{0}
\renewcommand{\theequation}{E.\arabic{equation}}
\renewcommand{\thesubsection}{E.\arabic{subsection}}
The proof is also similar to that for Theorem 4.2, and hence we shall mainly focus on the key differences. Analogously, as noted in the proof of Theorem 4.2,
it suffices to address real unstable poles, i.e., $p_1>0$, $p_2>0$. For $\alpha\geq 1$, define the set
\begin{equation}\nonumber
\begin{array}{l}
{\kern -6 pt} \Psi_\alpha  = \displaystyle \left\{ {\left( {{k_p},{k_i}} \right):} \right.{p_1} + {p_2} < \alpha{k_p} < \frac{{\alpha{k_i} + {p_1}{p_2}}}{z},\alpha{k_i} < 0,\\
{\kern 30 pt}  \displaystyle \left. {\left( {\alpha{k_p} - {p_1} - {p_2}} \right)\left( { - z\alpha{k_p} + \alpha{k_i} + {p_1}{p_2}} \right) >  - z\alpha{k_i}} \right\}.
\end{array}
\end{equation}
Then it follows similarly that
%It follows from the Routh-Hurwitz criterion that the PI controller
%${K_{PI}}\left( s \right) = {k_p} + (k_i/s)$ stabilizes ${\mathscr{P}_\alpha }$ for
%all $\alpha\in [1,~\mu)$ if and only if $(k_p,~k_i)\in \Psi_\alpha$ for all $\alpha\in [1,~\mu)$.
%Furthermore, for ${K_{PI}}\left( s \right)$ to stabilize $P(s)$, it is necessary that
%$(k_p,~k_i)\in \Psi_1$, which necessitates that ${p_1} + {p_2} < {k_p} < \left( {{k_i} + {p_1}{p_2}} \right)/z$ and
%in turn ${p_1} + {p_2} < {k_p} < p_1p_2/z$, whereas the latter follows from the fact that
%$k_i<0$. Hence,
\begin{eqnarray*}
k_\M^{\PI} &=& \sup \left\{\alpha  \geq 1:~(k_p,~k_i)\in \Psi_\alpha\right\} \\
& \leq & \sup \left\{ {\alpha  > 1:~\alpha \left( {z{k_p} - {k_i}} \right) < {p_1}{p_2},\left( {{k_p},{k_i}} \right) \in \Psi_1 } \right\}.
\end{eqnarray*}
Since for any $\left( {{k_p},{k_i}} \right) \in \Psi_1$, $k_i<0$, $k_p>p_1+p_2$, the supremum is achieved
at $k_p^* = {p_1} + {p_2},k_i^* = 0$; that is,
$$
k_\M^{\PI}=k_\M^\p\leq \frac{{{p_1}{p_2}}}{{z\left( {{p_1} + {p_2}} \right)}}.
$$
The upper bound, however, is attained since $\left( {{k^\ast_p},{k^\ast_i}} \right) \in \Psi_\alpha$ for
$\alpha=p_1p_2/\left(z(p_1+p_2)\right)$. This proves (\ref{2-44}).

%It follows that to maximize $\alpha $, it is equivalent to minimize $\left( {z{k_p} - {k_i}} \right)$. In light of the fact that ${k_p} > {p_1} + {p_2},{k_i} < 0$, we find %that $k_\M^{\PI}$ achieves its maximum at  On the other hand, note that $\left( {{k_p^*},{k_i^*}} \right) \in \Psi $ by submitting $k_p^* = {p_1} + {p_2}$, $k_i^* = 0$ into %$\Psi $. We are thus led to $k_\M^{\PI} = k_\M^\p = $.

%%%%%%%%%%%%%%%%%%%%Kp
To establish the maximal phase margin, consider similarly
%Clearly, for ${K_{PI}}\left( s \right)$ to stabilize $P(s)$, it gives rise to ${p_1} + {p_2} < {k_p} < \left( {{k_i} + {p_1}{p_2}} \right)/z$
%and ${p_1} < z$, ${p_2} < z$.
the open-loop frequency response
\begin{equation}\nonumber
L\left( {j\omega } \right) = \frac{{ j\omega-z }}{{\left( { j\omega-p_1 } \right)\left( {  j\omega- {p_2} } \right)}}\left( {{k_p} - j\frac{{{k_i}}}{\omega }} \right).
\end{equation}
The crossover frequencies of $L\left( {j\omega } \right)$ are found to be the solutions to the equation
\begin{equation}\label{T4.2-kp1}
{\kern -6 pt} {\omega ^6} + \left( {p_1^2 + p_2^2 - k_p^2} \right){\omega ^4} + \left( {p_1^2p_2^2 - {z^2}k_p^2 \!-\! k_i^2} \right){\omega ^2} \!- \!{z^2}k_i^2 = 0.
\end{equation}
Mimicking the proof of Theorem 4.2, we find that
\begin{equation}\nonumber
\theta _\M^{\PI} = \min \left\{ {\theta _{\PI}^ + ,\;\theta _{\PI}^ - } \right\},
\end{equation}
where
\begin{equation}\nonumber
\theta _{\PI}^ +  = \sup \left\{ {{\theta \left( \omega  \right)} \! > \!0:~ \! \omega \mbox{ is a solution to (\ref{T4.2-kp1})}, \left( {{k_p},{k_i}} \right) \! \in  \! \Psi_1 } \right\},
\end{equation}
\begin{equation}\nonumber
\theta _{\PI}^ -  =  \!- \! \inf \left\{ {{\theta \left( \omega  \right)} \! <  \! 0:~ \! \omega \mbox{ is a solution to (\ref{T4.2-kp1})}, \left( {{k_p},{k_i}} \right) \! \in  \! \Psi_1 } \right\},
\end{equation}
and
\begin{equation}\label{T4.2-kp2}
{\theta \left( \omega  \right)} = {\tan ^{ - 1}}\frac{{{\omega }}}{{{p_1}}} + {\tan ^{ - 1}}\frac{{{\omega }}}{{{p_2}}} - {\tan ^{ - 1}}\frac{{{\omega }}}{z} - {\tan ^{ - 1}}\frac{{ {k_i}}}{{{\omega}{k_p}}}.
\end{equation}
By reformulating the optimization of $\theta$ in terms of a nonlinear programming problem, and by
employing the KKT condition, we find that both $\theta^+_{\PI}$ and $\theta^-_{\PI}$
are achieved at
$k_i^* = 0$. As such, $\theta _\M^{\PI} = \theta _\M^\p$ and the equations (\ref{T4.2-kp1}) and (\ref{T4.2-kp2}) reduce to
\begin{equation}\label{T4.2-kp1+1}
{\omega ^4} + \left( {p_1^2 + p_2^2 - k_p^2} \right){\omega ^2} +  {p_1^2p_2^2 - {z^2}k_p^2} = 0,
\end{equation}
and
\begin{equation}\label{T4.2-kp2+2}
{\theta \left( \omega  \right)} = {\tan ^{ - 1}}\frac{{{\omega }}}{{{p_1}}} + {\tan ^{ - 1}}\frac{{{\omega }}}{{{p_2}}} - {\tan ^{ - 1}}\frac{{{\omega }}}{z}.
\end{equation}
%Accordingly, $\theta _{\PI}^ +$ and $\theta _{\PI}^ -$ reduce to
%\begin{equation}\nonumber
%\begin{array}{l}
%{\kern -5 pt}  \theta _{\PI}^ +  = \sup \left\{ {\theta \left( \omega  \right)  > 0:} \right.~\omega \mbox{ is a solution to (\ref{T4.2-kp1+1})},  \\
%{\kern 140 pt}  \left. {{p_1} + {p_2} < {k_p} < {p_1}{p_2}/z} \right\},
%\end{array}
%\end{equation}
%and
%\begin{equation}\nonumber
%\begin{array}{l}
%{\kern -4 pt}  \theta _{\PI}^ -  = -\inf \left\{ {\theta  \left( \omega  \right) < 0:} \right.~\omega \mbox{ is a solution to (\ref{T4.2-kp1+1})},  \\
%{\kern 140 pt}  \left. {{p_1} + {p_2} < {k_p} < {p_1}{p_2}/z} \right\},
%\end{array}
%\end{equation}
%where $\theta \left( \omega  \right)$ is now given by (\ref{T4.2-kp2+2}).
Note that the solutions to the equation (\ref{T4.2-kp1+1}) are given by
\begin{equation}\nonumber
\omega_1^2 = \frac{{\left( {k_p^2 \!- p_1^2 \!- p_2^2} \right) + \sqrt {{{\left( {k_p^2 \!- p_1^2 \!- p_2^2} \right)}^2} + 4\left( {k_p^2{z^2} \!- p_1^2p_2^2} \right)} }}{2},
\end{equation}
and
\begin{equation}\nonumber
\omega_2^2 = \frac{{\left( {k_p^2 \!- p_1^2 \!- p_2^2} \right) - \sqrt {{{\left( {k_p^2 \!- p_1^2 \!- p_2^2} \right)}^2} + 4\left( {k_p^2{z^2} \!- p_1^2p_2^2} \right)} }}{2}.
\end{equation}
It is straightforward to show that $\omega_1^2$ is monotonically increasing with $k_p$, and that $\omega_2^2$ decreases monotonically with $k_p$.
For ${p_1} + {p_2} < {k_p} < p_1p_2/z$, it then holds that
\begin{eqnarray}
\omega^2_1 & > & p_1p_2 + z(p_1+p_2), \label{crossover1} \\
\omega^2_2 & < & p_1p_2 - z(p_1+p_2). \label{crossover2}
\end{eqnarray}
Denote by ${\theta _i}\left( {{\omega _i}} \right) = \theta \left( {{\omega _i}} \right)$, %$\theta =\theta_i$ for $\omega=\omega_i$,
$i=1,~2$. In light of Lemma A.4 and (\ref{crossover1}) and (\ref{crossover2}), we have
\begin{eqnarray*}
%\theta_1 &=& {\tan ^{ - 1}}\frac{{{\omega_1 }}}{{{p_1}}} + {\tan ^{ - 1}}\frac{{{\omega_1 }}}{{{p_2}}} - {\tan ^{ - 1}}\frac{{{\omega_1 }}}{z} \\
%&=& \pi - \tan^{-1}\frac{(p_1+p_2)\omega_1}{\omega^2_1-p_1p_2}-{\tan ^{ - 1}}\frac{{{\omega_1 }}}{z} \\
%& > & 0,
{\theta _1 \left( {{\omega _1}} \right) } &=& \pi  - {\tan ^{ - 1}}\frac{{({p_1} + {p_2}){\omega _1}}}{{\omega _1^2 - {p_1}{p_2}}} - {\tan ^{ - 1}}\frac{{{\omega _1}}}{z} > 0,
\end{eqnarray*}
and
\begin{eqnarray*}
%\theta_2 &=& {\tan ^{ - 1}}\frac{{{\omega_2 }}}{{{p_1}}} + {\tan ^{ - 1}}\frac{{{\omega_2 }}}{{{p_2}}} - {\tan ^{ - 1}}\frac{{{\omega_2 }}}{z} \\
%&=& \tan^{-1}\frac{(p_1+p_2)\omega_2}{p_1p_2-\omega^2_2}-{\tan ^{ - 1}}\frac{{{\omega_2 }}}{z} \\
%& < & 0.
{\theta _2 \left( {{\omega _2}} \right) }& =& {\tan ^{ - 1}}\frac{{({p_1} + {p_2}){\omega _2}}}{{{p_1}{p_2} - \omega _2^2}} - {\tan ^{ - 1}}\frac{{{\omega _2}}}{z} < 0.
\end{eqnarray*}
As a result, to determine $\theta _\M^{\PI}$, we maximize $\theta_1\left( {{\omega _1}} \right)$ and minimize $\theta_2\left( {{\omega _2}} \right)$.
Taking the derivative of $\theta_i\left( {{\omega _i}} \right)$ with respect to $\omega_i$ %leads to
gives rise to
\begin{eqnarray}
\frac{{d{{\theta_i(\omega_i) }}}}{{d{\omega_i}}} &=& \frac{{{p_1}}}{{\omega_i^2 + p_1^2}} + \frac{{{p_2}}}{{\omega_i^2 + p_2^2}} - \frac{z}{{\omega_i^2 + {z^2}}} \nonumber  \\
&=& \frac{q(\omega_i)}{(\omega_i^2+p^2_1)(\omega_i^2+p^2_2)(\omega_i^2+z^2)}, \label{derivative2}
\end{eqnarray}
%Taking the derivative of $\theta_i\left( {{\omega _i}} \right)$ with respect to $\omega_i$ gives rise to the equations (\ref{derivative2}) and (\ref{derivative1}) with $\omega  = {\omega _i},i = 1,2$, given in the proof of Theorem 4.2.
where
\begin{equation}\label{derivative1}
\begin{array}{l}
{\kern -10 pt}   q(\omega_i) = \left( {{p_1} + {p_2} - z} \right){\omega_i ^4} + \left( {\left( {{p_1} + {p_2}} \right)\left( {{z^2} + {p_1}{p_2}} \right) - } \right.\\
{\kern 20 pt} \left. {z\left( {p_1^2 + p_2^2} \right)} \right){\omega_i ^2} + {p_1}{p_2}z\left( {z\left( {{p_1} + {p_2}} \right) - {p_1}{p_2}} \right).
\end{array}
\end{equation}
%\begin{eqnarray}
%\frac{{d{\theta }}}{{d{\omega}}} &=& \frac{{{p_1}}}{{\omega^2 + p_1^2}} + \frac{{{p_2}}}{{\omega^2 + p_2^2}} - \frac{z}{{\omega^2 + {z^2}}} \nonumber  \\
%&=& \frac{q(\omega)}{(\omega^2+p^2_1)(\omega^2+p^2_2)(\omega^2+z^2)}, \label{derivative2}
%\end{eqnarray}
%where
%\begin{equation}\label{derivative1}
%\begin{array}{l}
%{\kern -10 pt}   q(\omega ) = \left( {{p_1} + {p_2} - z} \right){\omega ^4} + \left( {\left( {{p_1} + {p_2}} \right)\left( {{z^2} + {p_1}{p_2}} \right) - } \right.\\
%{\kern 20 pt} \left. {z\left( {p_1^2 + p_2^2} \right)} \right){\omega ^2} + {p_1}{p_2}z\left( {z\left( {{p_1} + {p_2}} \right) - {p_1}{p_2}} \right).
%\end{array}
%\end{equation}
%\begin{equation}\label{derivative1}
%\begin{aligned}
%&q(\omega)=\left( {{p_1} + {p_2} - z} \right)\omega^4 + \left(\left( {{p_1} + {p_2}} \right)\left( {{z^2} + {p_1}{p_2}} \right)
%- z\left( {p_1^2 + p_2^2} \right) \right)\omega^2
%\\
%&~~~~~~~~~~~~~+ {p_1}{p_2}z\left( z\left( {{p_1} + {p_2}} \right) - {p_1}{p_2} \right).
%\end{aligned}
%\end{equation}
%\begin{aligned}
%&q(\omega)=\left( {{p_1} + {p_2} - z} \right)\omega^4 + \left(\left( {{p_1} + {p_2}} \right)\left( {{z^2} + {p_1}{p_2}} \right)
%- z\left( {p_1^2 + p_2^2} \right) \right)\omega^2
%\\
%&~~~~~~~~~~~~~+ {p_1}{p_2}z\left( z\left( {{p_1} + {p_2}} \right) - {p_1}{p_2} \right).
%\end{aligned}
%\end{equation}
Setting $d\theta_i/d\omega_i=0$ yields $q(\omega_i)=0$.
%with $q(\omega_i )$ given in (\ref{derivative1}).
%\begin{equation}\label{derivative1}
%\begin{aligned}
%&\left( {{p_1} + {p_2} - z} \right)\omega^4 + \left(\left( {{p_1} + {p_2}} \right)\left( {{z^2} + {p_1}{p_2}} \right)
%- z\left( {p_1^2 + p_2^2} \right) \right)\omega^2
%\\
%&~~~~~~~~~~~~~+ {p_1}{p_2}z\left( z\left( {{p_1} + {p_2}} \right) - {p_1}{p_2} \right) \right)=0.
%\end{aligned}
%\end{equation}
From the inequality $p_1+p_2<p_1p_2/z$, it follows that $p_1>z$, $p_2>z$. These inequalities together suggest that
$$
\begin{aligned}
& {{p_1} + {p_2} - z}>0, \\
& \left( {{p_1} + {p_2}} \right)\left( {{z^2} + {p_1}{p_2}} \right)
- z\left( {p_1^2 + p_2^2} \right)>0, \\
& {p_1}{p_2}z\left( z\left( {{p_1} + {p_2}} \right) - {p_1}{p_2} \right)  <0.
\end{aligned}
$$
According to the Descartes' Rule (Lemma A.3), we conclude that the polynomial equation $q(\omega_i)=0$ admits one and only one
positive solution $\omega^2_0$. Assume that the solution meets the condition $\omega^2_0\geq p_1p_2 - z(p_1+p_2)$.
Substituting $\omega_0$ into (\ref{derivative1}), we arrive at the condition
$$
\begin{aligned}
&\left( {{p_1} + {p_2} - z} \right)\omega^4_0 + \left(\left( {{p_1} + {p_2}} \right)\left( {{z^2} + {p_1}{p_2}} \right)
- z\left( {p_1^2 + p_2^2} \right) \right)\omega^2_0
\\
&~~~~~~~~~~~~~+ {p_1}{p_2}z\left( z\left( {{p_1} + {p_2}} \right) - {p_1}{p_2} \right)>0,
\end{aligned}
$$
thus leading to contradiction.
This means that any positive root of $q(\omega_i)$ must be such that $\omega^2_0< p_1p_2 - z(p_1+p_2)$,
implying that $d\theta_1/d\omega_1 \neq 0$, and that $\theta _{\PI}^ +$ is achieved
at the endpoints $k_p=p_1+p_2$ or $k_p=p_1p_2/z$. In particular, since $q(\omega_1)>0$,
$\theta_1\left( {{\omega _1}} \right)$ is monotonically increasing. Therefore, $\theta _{\PI}^ +$ is achieved
at $k_p=p_1p_2/z$, i.e., $\theta _{\PI}^ +  = {\theta _1}\left( {{\omega _1}\left( {{p_1}{p_2}/z} \right)} \right)$. On the other hand, since $\omega^2_2 \in [0,~p_1p_2 - z(p_1+p_2))$,
we assert that $d\theta_2/d\omega_2= 0$ at $\omega_0$, which can be found explicitly by solving the
positive root of $q(\omega_2)$.
% given by (\ref{PI-phase-analytical}).
Note that at the endpoints $k_p=p_1+p_2$ or $k_p=p_1p_2/z$, $\theta_2\left( {{\omega _2}} \right)=0$. Since
$\omega^2_0$ is the sole positive root, the infimum of $\theta_2$ is achieved at $\omega_0$, i.e, $\theta _{\PI}^ -  = {\theta _2}\left( {{\omega _0}} \right)$.
%which is given in (\ref{PI-phase-analytical}).

It remains to show that $\theta _{\PI}^+\geq \theta _{\PI}^-$. Towards this end, we examine
${\theta _1\left( {{\omega _1}} \right)} + { {\theta _2\left( {{\omega _2}} \right)}}$. Since
${\theta _1\left( {{\omega _1}} \right)}>0$ and ${\theta _2\left( {{\omega _2}} \right)} < 0$, it follows that $\theta _{\PI}^+\geq \theta _{\PI}^-$ provided that $ {\theta _1\left( {{\omega _1}} \right)} + {\theta _2\left( {{\omega _2}} \right)} \geq 0$. Note that
\begin{equation} \label{theta1}
{\kern -10 pt}{\theta _1\left( {{\omega _1}} \right)} \!+ \!{\theta _2\left( {{\omega _2}} \right)} = \sum\limits_{i = 1}^2 \left({{{\tan }^{ - 1}}\frac{{{\omega _i}}}{{{p_1}}} \!+\! {{\tan }^{ - 1}}\frac{{{\omega _i}}}{{{p_2}}}\!-\! {{\tan }^{ - 1}}\frac{{{\omega _i}}}{z}}\right) .
\end{equation}
Note also from (\ref{T4.2-kp1+1}) that
\begin{eqnarray}
{\omega _1}{\omega _2} &=& \sqrt {p_1^2p_2^2 - {z^2}k_p^2},  \label{omega1 and 2-time} \\
{\omega _1} + {\omega _2} &=& \sqrt {k_p^2 - p_1^2 - p_2^2 + 2\sqrt {p_1^2p_2^2 - {z^2}k_p^2} }.\label{omega1 and 2-plus}
\end{eqnarray}
Furthermore, from the latter equality, we find that
\begin{equation} \nonumber
\frac{{d\left( {{\omega _1} + {\omega _2}} \right)}}{{d{k_p}}} = \frac{1}{{{\omega _1} + {\omega _2}}}\left( 1 - \frac{{{z^2}}}{{\omega _1}{\omega _2}} \right){k_p}.
\end{equation}
It is thus clear that ${\omega _1}{\omega _2}$ is monotonically decreasing with ${k_p}$, while
${{\omega _1} + {\omega _2}}$ is monotonically increasing with $k_p$ whenever ${\omega _1}{\omega _2}\geq z^2$
and monotonically decreasing with $k_p$ whenever ${\omega _1}{\omega _2}\leq z^2$.
Without loss of generality, assume that ${p_2} \le {p_1}$. We now evaluate ${\theta _1\left( {{\omega _1}} \right)} + {\theta _2\left( {{\omega _2}} \right)}$ in four different cases.

\medskip\noindent
{\bf Case 1: ${\omega _1}{\omega _2}\leq z^2< p^2_2\leq p^2_1$. }
It is immediate from Lemma A.4 that
\begin{equation} \label{inequa1}
{\kern -4pt}{\theta _1\left( {{\omega _1}} \right)} + {\theta _2\left( {{\omega _2}} \right)} = \sum\limits_{i = 1}^2 {{{\tan }^{ - 1}}\frac{{{p_i}\left( {{\omega _1} + {\omega _2}} \right)}}{{p_i^2- {\omega _1}{\omega _2}}}}   -    {\tan ^{ - 1}}\frac{{z\left( {{\omega _1} + {\omega _2}} \right)}}{{{z^2}- {\omega _1}{\omega _2}}}.
\end{equation}
From (\ref{omega1 and 2-plus}), it is clear that
$$
({\omega _1} + {\omega _2})^2\geq k_p^2 - p_1^2 - p_2^2> p_1p_2.
$$
This leads to
%gives rise to
$$
({\omega _1} + {\omega _2})^2p_1p_2> p^2_1p^2_2\geq (p^2_1-\omega_1\omega_2)(p^2_2-\omega_1\omega_2),
$$
or equivalently,
$$
\prod^2_{i=1} \frac{p_i ({\omega _1} + {\omega _2})}{p^2_i-\omega_1\omega_2}>1.
$$
Denote by
$$
\gamma_i = \frac{p_i ({\omega _1} + {\omega _2})}{p^2_i-\omega_1\omega_2},~~~~~i=1,~2.
$$
Then according to Lemma A.4, we have
\begin{equation}\label{inequa1a}
\sum^2_{i=1}{{{\tan }^{ - 1}}\frac{{{p_i}\left( {{\omega _1} + {\omega _2}} \right)}}{{p_i^2- {\omega _1}{\omega _2}}}}
=\pi- \tan^{-1}\frac{\gamma_1+\gamma_2}{\gamma_1 \gamma_2-1}.
\end{equation}
Together with (\ref{inequa1}), this leads to 
$$
{\theta _1\left(\omega _1\right)} + \theta _2\left(\omega _2\right)=
\pi- \tan^{-1}\frac{\gamma_1+\gamma_2}{\gamma_1 \gamma_2-1}-
{\tan ^{ - 1}}\frac{{z\left( {{\omega _1} + {\omega _2}} \right)}}{{{z^2}- {\omega _1}{\omega _2}}},
$$
and hence $\theta_1\left( {{\omega _1}} \right)+\theta_2\left( {{\omega _2}} \right)\geq 0$.
%
%Since ${\omega _1}{\omega _2}$ and ${{\omega _1} + {\omega _2}}$ are both
%monotonically decreasing with ${k_p}$ under this circumstance, for $i=1,~2$,
%$$
%{{{\tan }^{ - 1}}\frac{{{p_i}\left( {{\omega _1} + {\omega _2}} \right)}}{{p_i^2- {\omega _1}{\omega _2}}}}
%$$
%is monotonically decreasing with $k_p$, which attains it minimum at $k_p=p_1p_2/z$. Since at $k_p=p_1p_2/z$,
%\begin{eqnarray*}
%\omega_1 &=& \sqrt{\frac{p^2_1p^2_2}{z^2}-(p^2_1+p^2_2)}, \\
%\omega_2 &=& 0,
%\end{eqnarray*}
%we have
%$$
%{{{\tan }^{ - 1}}\frac{{{p_i}\left( {{\omega _1} + {\omega _2}} \right)}}{{p_i^2- {\omega _1}{\omega _2}}}}
%\geq \tan^{ - 1}\frac{\sqrt{(p_1p_2/z)^2-(p^2_1+p^2_2)}}{p_i}.
%$$
%Note however that
%$$
%p_1p_2<\frac{p^2_1p^2_2}{z^2}-(p^2_1+p^2_2).
%$$
%Hence, accordingly Lemma ,
%$$
%\sum^2_{i=1}\tan^{ - 1}\frac{\sqrt{(p_1p_2/z)^2-(p^2_1+p^2_2)}}{p_i}=\pi-
%\tan^{-1} \frac{(p_1+p_2)\sqrt{(p_1p_2/z)^2-(p^2_1+p^2_2)}}{(p_1p_2/z)^2-(p^2_1+p^2_2+p_1p_2)}.
%$$
%This consequently leads us to
%$$
%\sum^2_{i=1}{{{\tan }^{ - 1}}\frac{{{p_i}\left( {{\omega _1} + {\omega _2}} \right)}}{{p_i^2- {\omega _1}{\omega _2}}}}
%\geq \pi-
%\tan^{-1} \frac{(p_1+p_2)\sqrt{(p_1p_2/z)^2-(p^2_1+p^2_2)}}{(p_1p_2/z)^2-(p^2_1+p^2_2+p_1p_2)},
%$$
%and in view of (\ref{inequa1}), to ${\theta _1} + {\theta _2}\geq 0$.

\medskip\noindent
{\bf Case 2: $z^2<{\omega _1}{\omega _2}\leq p^2_2\leq p^2_1$. } In this case, (\ref{inequa1a}) remains to hold but 
$$
\sum\limits_{i = 1}^2{{\tan }^{ - 1}}\frac{{{\omega _i}}}{z}=\pi -
{\tan ^{ - 1}}\frac{{z\left( {{\omega _1} + {\omega _2}} \right)}}{{{\omega _1}{\omega _2}-{z^2}}},
$$
By substituting (\ref{inequa1a}), we have 
\begin{equation} \label{inequa2}
\begin{array}{l}
{\kern -12 pt} \displaystyle{ {\theta _1} \left( {{\omega _1}} \right) + {\theta _2}\left( {{\omega _2}} \right) = 
{\tan ^{ - 1}}\frac{{z({\omega _1} + {\omega _2})}}{{{\omega _1}{\omega _2} - {z^2}}} - {\tan ^{ - 1}}\frac{{{\gamma _1} + {\gamma _2}}}{{{\gamma _1}{\gamma _2} - 1}}.}
%\sum\limits_{i = 1}^2 {{{\tan }^{ - 1}}\frac{{{p_i}\left( {{\omega _1} + {\omega _2}} \right)}}{{p_i^2 - {\omega _1}{\omega _2}}}}  + \\
%{\kern 118 pt} \displaystyle {\tan ^{ - 1}}\frac{{z\left( {{\omega _1} + {\omega _2}} \right)}}{{{\omega _1}{\omega _2}-{z^2} }} - \pi \\
%{\tan ^{ - 1}}\frac{{z({\omega _1} + {\omega _2})}}{{{\omega _1}{\omega _2} - {z^2}}} - {\tan ^{ - 1}}\frac{{{\gamma _1} + {\gamma _2}}}{{{\gamma _1}{\gamma _2} - 1}}.
\end{array}
%{\kern -4 pt}  {{\theta _1}\left( \!{{\omega _1}} \!\right) \!+ \! {\theta _2}\left( \!{{\omega _2}} \!\right) \! = \! \sum\limits_{i = 1}^2 {{{\tan }^{ - 1}}\frac{{{p_i}\left( \! {{\omega _1} \! +  \! {\omega _2}} \! \right)}}{{p_i^2  \! - \! {\omega _1}{\omega _2}}}}  \! + \! {{\tan }^{ - 1}}\frac{{z\left( \! {{\omega _1} \! + \!  {\omega _2}} \! \right)}}{{{z^2} \! -  \! {\omega _1}{\omega _2}}} \! - \! \pi .}
\end{equation}
It follows from a straightforward, albeit tedious, calculation that 
$$
\frac{{z({\omega _1} + {\omega _2})}}{{{\omega _1}{\omega _2} - {z^2}}} \geq \frac{{{\gamma _1} + {\gamma _2}}}{{{\gamma _1}{\gamma _2} - 1}}.
$$
Therefore, $\theta_1\left( {{\omega _1}} \right)+\theta_2\left( {{\omega _2}} \right)\geq 0$.

\medskip\noindent
{\bf Case 3: $z^2< p^2_2<{\omega _1}{\omega _2}\leq p^2_1$. } Likewise, by invoking Lemma A.4, we obtain
\begin{eqnarray}
{\theta _1\left( {{\omega _1}} \right) } \!+\! {\theta _2\left( {{\omega _2}} \right) } \! \! &=& \!  \! {{{\tan }^{ - 1}}\frac{{{p_1}\left( {{\omega _1} \! + \! {\omega _2}} \right)}}{{p_1^2 \! - \!  {\omega _1}{\omega _2}}}}
-{{{\tan }^{ - 1}}\frac{{{p_2}\left( {{\omega _1} \! + \!  {\omega _2}} \right)}}{{{\omega _1}{\omega _2}} \! - \! p_2^2}} \nonumber \\
\hfill & \hfill & ~~~
{\kern -18 pt} +{\tan ^{ - 1}}\frac{{z\left( {{\omega _1} \! + \!  {\omega _2}} \right)}}{{{\omega _1}{\omega _2} \! - \! {z^2}}}. \label{inequa3}
\end{eqnarray}
Using (\ref{inequa3}), the proof can be carried out analogously as in Case 2, which also involves a tedious calculation, leading to 
the same conclusion that $\theta_1 \left( {{\omega _1}} \right) +\theta_2 \left( {{\omega _2}} \right) \geq 0$.
%As in Case 2, it follows that
%$$
%{\theta _1}\left( {{\omega _1}} \right) + {\theta _2}\left( {{\omega _2}} \right) = {\tan ^{ - 1}}\frac{{z({\omega _1} + {\omega _2})}}{{{\omega _1}{\omega _2} - {z^2}}} - {\tan ^{ - 1}}\frac{{{\gamma _1} + {\gamma _2}}}{{{\gamma _1}{\gamma _2} - 1}},
%$$
%which leads to the assertion $\theta_1 \left( {{\omega _1}} \right) +\theta_2 \left( {{\omega _2}} \right) \geq 0$.
%Since in this case (\ref{inequa-case2}) holds, it follows at once that $\theta_1(\omega _1)+\theta_2(\omega _2)\geq 0$.

\medskip\noindent
{\bf Case 4: $z^2< p^2_2\leq p^2_1<{\omega _1}{\omega _2}$. } From (\ref{omega1 and 2-time}), we have ${\omega _1}{\omega _2} < {p_1}{p_2} \le p_1^2$. 
As such, this case is not possible and hence is precluded. 
%
%We ascertain  that this case leads to contradiction since it follows from (\ref{omega1 and 2-time}) that ${\omega _1}{\omega _2} \le {p_1}{p_2} \le p_1^2$.
%As in Case 3, it follows that
%\begin{equation} \label{inequa4}
%\begin{array}{l}
%{\kern -15 pt} \displaystyle {\theta _1}\left( {{\omega _1}} \right) + {\theta _2}\left( {{\omega _2}} \right) = \pi  + {\tan ^{ - 1}}\frac{{z\left( {{\omega _1} + {\omega _2}} \right)}}{{{\omega _1}{\omega _2} - {z^2}}}\\
%{\kern 98 pt} \displaystyle  - \sum\limits_{i = 1}^2 {{{\tan }^{ - 1}}\frac{{{p_i}\left( {{\omega _1} + {\omega _2}} \right)}}{{{\omega _1}{\omega _2} - p_i^2}}} ,
%\end{array}
%{\kern -4 pt}  {{\theta _1}\left( \!{{\omega _1}} \!\right) \! + \! {\theta _2}\left( \!{{\omega _2}} \!\right) \! = \! \pi  \! + \! {{\tan }^{ - 1}}\frac{{z\left( \!{{\omega _1} \! + \! {\omega _2}} \!\right)}}{{{\omega _1}{\omega _2} \! - \! {z^2}}}  \! - \!  \sum\limits_{i = 1}^2 {{{\tan }^{ - 1}}\frac{{{p_i}\left( \!{{\omega _1} \! + \! {\omega _2}} \!\right)}}{{{\omega _1}{\omega _2} \! - \! p_i^2}}},}
%end{equation}
%which leads to the assertion $\theta_1 \left( {{\omega _1}} \right) +\theta_2 \left( {{\omega _2}} \right) \geq 0$.

\medskip\noindent
In summary, we have showed that for all $k_p$ such that $p_1+p_2<k_p<{p_1}{p_2}/z$, it is always true that
$\theta_1 \left( {{\omega _1}} \right) +\theta_2\left( {{\omega _2}} \right) \geq 0$. Consequently, $\theta_1 \left( {{\omega _1}} \right) \geq  -\theta_2\left( {{\omega _2}} \right)$, and hence $\theta _{\PI}^+\geq \theta _{\PI}^-$. This establishes (\ref{2-44-phase}) together with (\ref{PI-phase-analytical}).
Finally, to find the optimal $k_p$, we note that $\omega_0$ can be alternatively determined from the equation
$$
\frac{p_1}{\omega^2_0 + p_1^2} + \frac{p_2}{\omega^2_0 + p_2^2} = \frac{z}{\omega^2_0 + z^2},
$$
or equivalently,
$$
\frac{(p_1+p_2)(\omega^2_0+p_1p_2)}{(\omega^2_0 + p_1^2)(\omega^2_0 + p_2^2)} = \frac{z}{\omega^2_0 + z^2}.
$$
Since as a crossover frequency $\omega_0$ satisfies the relation
$$
\frac{k^2_p (\omega^2_0 + z^2)}{(\omega^2_0 + p_1^2)(\omega^2_0 + p_2^2)} = 1,
$$
the optimal proportional gain $k^\ast_p$ is obtained as
$$
k^\ast_p = \sqrt{\frac{p_1+p_2}{z}(\omega^2_0+p_1p_2)}.
$$
This establishes (\ref{optimal-kp}) and the proof is now completed. The proof for the case of complex conjugate poles can be established analogously and hence is omitted.
$\hfill\blacksquare$

\section*{Acknowledgement}
The authors would like to thank Dr. Islam Boussaada, Universit\'e Paris Saclay, CNRS-CentraleSupelec, France,
for helpful discussions.

\bibliographystyle{IEEEtran}
\bibliography{IEEEabrv,myieeebibfile}
%%%%%%%%%%%%{References}

\ifCLASSOPTIONcaptionsoff
  \newpage
\fi

\vspace{- 1.0 cm}
\begin{IEEEbiography}[{\includegraphics[width=1in,height=1.25in,clip,keepaspectratio]{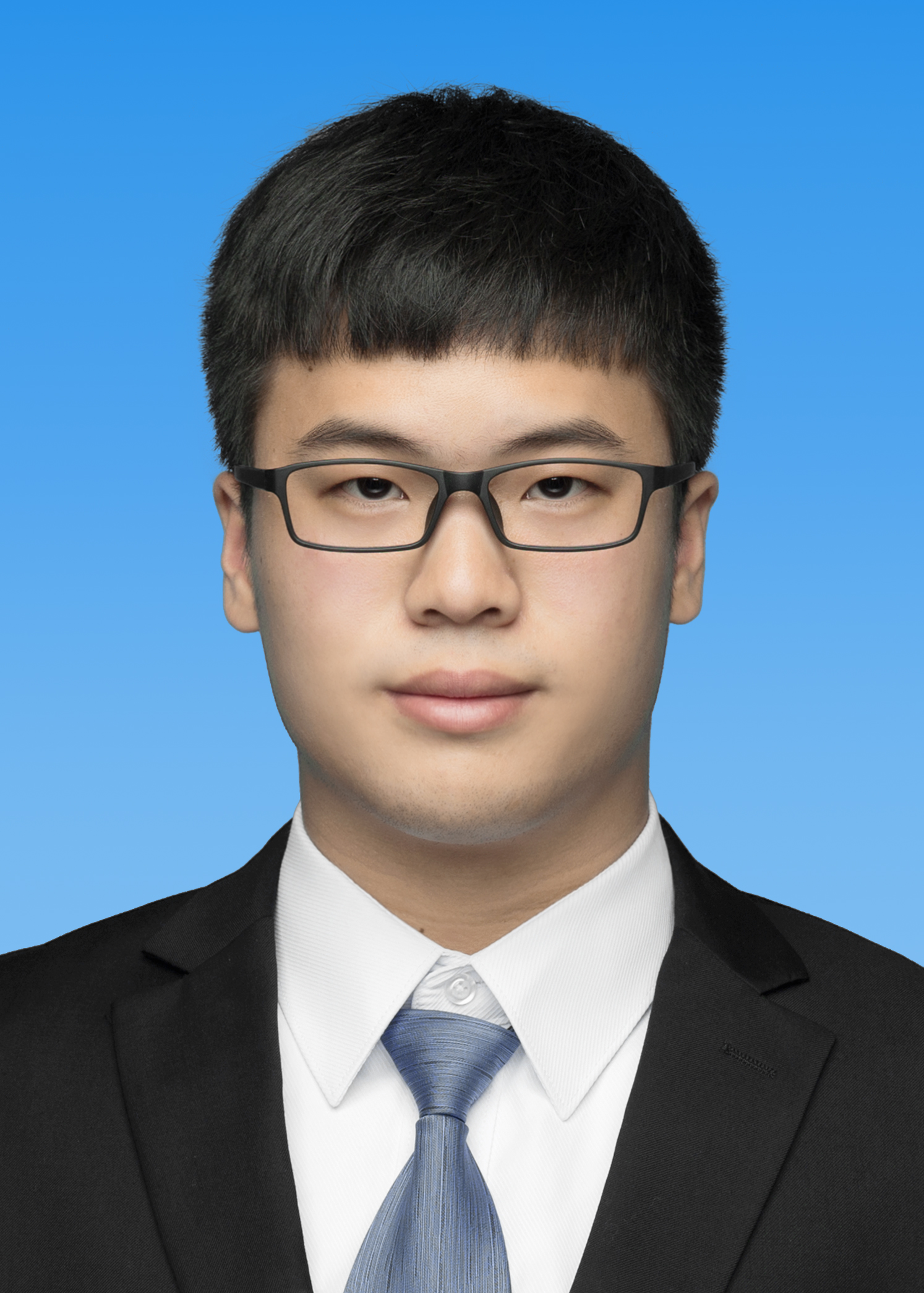}}]
{Qi Mao} received the B.S. degree in electrical engineering and automation from the Tianjin University of Science and Technology, Tianjin, China, in 2015, and the M.S. degree in the control science and engineering from the Tianjin University, Tianjin, China, in 2018, respectively. He is currently working toward the Ph.D. degree in the Department of Electrical Engineering, City University of Hong Kong, Kowloon, Hong Kong. His current research interests include PID control, time-delay systems, multiagent systems, and flight control.
\end{IEEEbiography}

\vspace{- 1.0cm}
\begin{IEEEbiography}[{\includegraphics[width=1in,height=1.25in,clip,keepaspectratio]{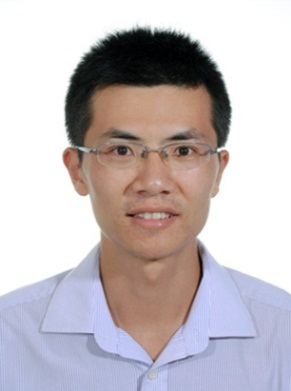}}]
{Yong Xu} (Member, IEEE) was born in Zhejiang Province, China, in 1983. He received the B.S.degree in information engineering from Nanchang Hangkong University, Nanchang, China, in 2007, the M.S. degree in control science and engineering from Hangzhou Dianzi University, Hangzhou, China, in 2010, and the Ph.D. degree in control science and engineering from Zhejiang University, Hangzhou, China, in 2014. He was a Visiting Student with the Department of Electronic and Computer Engineering, Hong Kong University of Science and Technology, Hong Kong, from June 2013 to November 2013, and was a Research Fellow from February 2018 to August 2018. He was honored Pearl River Young Scholars Program of Guangdong Province in 2017. He is currently a Professor with the School of Automation, Guangdong University of Technology, Guangzhou, China. His research interests include PID control, networked control systems, state estimation, and positive systems.
\end{IEEEbiography}

\vspace{- 1.0 cm}
\begin{IEEEbiography}[{\includegraphics[width=1in,height=1.25in,clip,keepaspectratio]{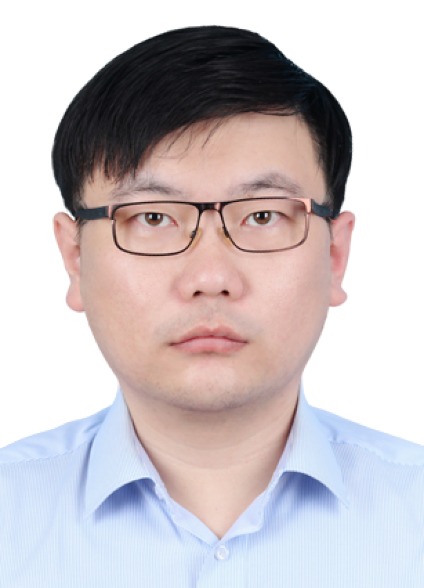}}]
{Jianqi Chen} received the B.Eng. degree from Zhejiang University, Hangzhou, China, in 2014, and the Ph.D. degree from the City University of Hong Kong, Hong Kong, in 2020. He is currently a Postdoctoral Researcher with the Department of Electrical Engineering, City University of Hong Kong. His research interests include PID control, time-delay systems, networked control, and cyber-physical systems. Dr. Chen was the recipient of the Guan Zhao-Zhi Best Paper Award at the 38th Chinese Control Conference in 2019.
\end{IEEEbiography}

%\begin{IEEEbiography}[{\includegraphics[width=1in,height=1.25in,clip,keepaspectratio]{Islam.pdf}}]
%{Islam Boussaada} received his Master in Mathematics from Carthage University, and an M.Sc. degree in Pure Mathematics from %University Paris 7 in 2004. In December 2008, he defended his Ph.D. degree in Mathematics from University of Rouen %Normandy. In June 2016, he received his HDR degree (French Habilitation) in Physics from University Paris %Saclay-Universit\'{e}. In 2010, IB was appointed for one year as a post-doctoral fellow in the control of time-delay %systems at L2S, Universit\'{e} Paris Saclay, CentraleSupelec-CNRS. From September 2011 until September 2017, he has been an %associate professor at IPSA and an associate researcher at MODESTY Team of L2S. Since October 2017, IB is appointed %associate researcher at DISCO Team and full professor at IPSA.  His research interests belong to the qualitative theory of %dynamical systems and its application in control problems.
%\end{IEEEbiography}

\vspace{- 1.0 cm}
\begin{IEEEbiography}[{\includegraphics[width=1in,height=1.25in,clip,keepaspectratio]{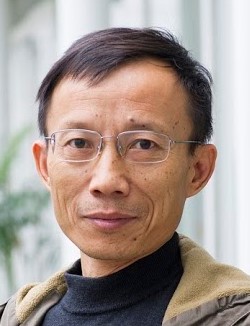}}]
{Jie Chen} ((S’87-M’89-SM’98-F’07)) received the B.S. degree in aerospace engineering from Northwestern Polytechnic University, Xi’an, China, in 1982, the M.S.E. degree in electrical engineering, the M.A. degree in mathematics, and the Ph.D. degree in electrical engineering all from the University of Michigan, Ann Arbor, MI, USA, in 1985, 1987, and 1990, respectively.

From 1994 to 2014, he was a Professor with the University of California, Riverside, CA, USA, and was a Professor and a Chair with the Department of Electrical Engineering from 2001 to 2006. He is a Chair Professor with the Department of Electrical Engineering, City University of Hong Kong, Hong Kong. His research interests include linear multi-variable systems theory, system identification, robust control, optimization, time-delay systems, networked control, and multi-agent systems.

Dr. Chen is a fellow of American Association for the Advancement of Science, a fellow of International Federation of Automatic Control, and an IEEE Distinguished Lecturer. He is currently an Associate Editor for SIAM Journal on Control and Optimization and International Journal of Robust and Nonlinear Control.
\end{IEEEbiography}

%\vspace{- 1.2 cm}
\begin{IEEEbiography}[{\includegraphics[width=1in,height=1.25in,clip,keepaspectratio]{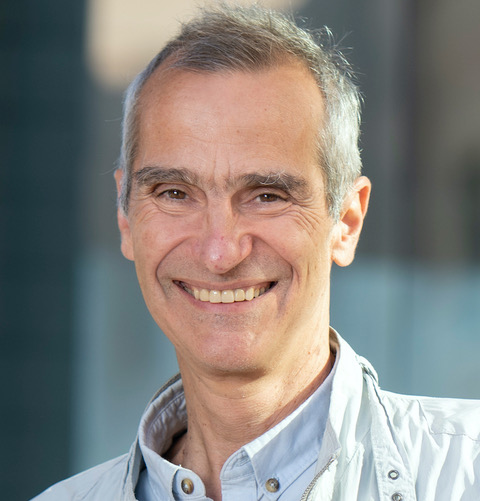}}]
{Tryphon T. Georgiou} (M’79-SM’99-F’00) is a UCI Distinguished Professor in mechanical and aerospace engineering at the University of California, Irvine. He was educated at the National Technical University of Athens, Greece (Diploma in Mechanical and Electrical Engineering, 1979), and the University of Florida, Gainesville (PhD 1983). Prior to joining the University of California, Irvine, he served on the faculty at the University of Minnesota, Iowa State University, and Florida Atlantic University.

Dr Georgiou has received the George S. Axelby Outstanding Paper award of the IEEE Control Systems Society for the years 1992, 1999, 2003, and 2017, he is a Fellow of the Institute of Electrical and Electronic Engineers (IEEE), a Fellow of the International Federation of Automatic Control (IFAC), a Fellow of the Society for Industrial and Applied Mathematics (SIAM), and a Foreign Member of the Royal Swedish Academy of Engineering Sciences (IVA).
\end{IEEEbiography}

%\newpage


%\vfill


\end{document}